\documentclass[draft]{amsart}
\usepackage{amsmath,amssymb, amsbsy}
 \usepackage{paralist}
\usepackage[dvips]{graphicx}
\newcommand{\supess}{\mathop{\rm ess\,sup}}
\newcommand{\infess}{\mathop{\rm ess\,inf}}
\newcommand{\R}{{\mathbb R}}
\newcommand{\N}{{\mathbb N}}

\newcommand{\weakly}{\rightharpoonup}
\renewcommand{\a }{\alpha }

\renewcommand{\d }{\delta }
\newcommand{\D }{\Delta }
\newcommand{\e }{\varepsilon }

\newcommand{\La }{\Lambda }
\newcommand{\n }{\nabla }

\newcommand{\Di}{{\mathcal D}^{1,2}(\R^N)}
\newcommand{\alchi}{\raisebox{1.7pt}{$\chi$}}

\newcommand{\media}{\mkern12mu\hbox{\vrule height4pt depth-3.2pt
    width6pt} \mkern-17.9mu\int}

\newenvironment{pf}{\noindent{\sc Proof}.\enspace}{\hfill\qed
\medskip}
\newenvironment{pfn}[1]{\noindent{\bf Proof of {#1}.\enspace}}{\hfill\qed
\medskip}
\newtheorem{Theorem}{Theorem}[section]
\newtheorem{Corollary}[Theorem]{Corollary}
\newtheorem{Lemma}[Theorem]{Lemma}
\newtheorem{Proposition}[Theorem]{Proposition}
\theoremstyle{definition}

\newtheorem{Example}[Theorem]{Example}

\newtheorem{remark}[Theorem]{Remark}
\begin{document}

\title[Schr\"odinger operators with inverse-square 
anisotropic 
potentials]{On Schr\"odinger operators with multisingular inverse-square 
anisotropic potentials}

\author[Veronica Felli \and Elsa M. Marchini \and Susanna
Terracini]{Veronica Felli \and Elsa M. Marchini \and Susanna
  Terracini}

\address{Universit\`a degli Studi di Milano Bicocca, Dipartimento di Matematica e
  Applicazioni, Via Cozzi 53, 20125 Milano, Italy.}  \email{{\tt
    veronica.felli@unimib.it, elsa.marchini@unimib.it,
    susanna.terracini@unimib.it}. }

\date{January 15, 2007} 

\thanks{Supported by Italy MIUR, national project
``Variational Methods and Nonlinear Differential Equations''.
}

\begin{abstract}
  \noindent We study positivity, localization of binding and essential
  self-adjointness properties of a class of Schr\"odinger operators
  with many anisotropic inverse square singularities, including the case of multiple dipole potentials. 
\end{abstract}

\maketitle

\section{Introduction and statement of the main results}\label{intro}

In this paper we analyze some basic spectral properties of 
Schr\"odinger operators associated with potentials possessing multiple anisotropic 
singularities of degree $-2$. The interest in such a class of operators arises in nonrelativistic molecular physics,
where the interaction between an electric charge  and the dipole moment of a
molecule con be described by an inverse square potential with 
an anisotropic coupling strength. 
More precisely,  the Schr\"odinger
operator acting on the wave function of an electron interacting with a polar
molecule (supposed to be point-like) can be written as
$$
H\!=\!-\frac{\hbar^2}{2m}\,\Delta+e\,\frac{x\cdot {\mathbf
    D}}{|x|^3}-E,
$$
where $e$ and $m$ denote respectively the charge and the mass of the
electron and ${\mathbf D}$ is the dipole moment of the molecule, see
\cite{leblond}. Therefore, in crystalline matter, the presence of many 
dipoles leads to consider multisingular Schr\"odinger operators of the form
\begin{equation}\label{eq:17}
-\D
-{\displaystyle{\sum_{i=1}^k}}
\dfrac{\lambda_i\,(x-a_i)\cdot{\mathbf d}_i}{|x-a_i|^3},
\end{equation}
where  $k\in\N$, $(a_1,\dots,a_k)\in \R^{k N}$, $N\geq 3$, $a_i\neq a_j$ for $i\neq j$,
 $(\lambda_1,\dots,\lambda_k)\in\R^k$, 
$({\mathbf d}_1,\dots,{\mathbf d}_k)\in\R^{k N}$, $\lambda_i>0$ and $|{\mathbf d}_i|=1$ for any
$i=1,\dots,k$.

Potentials of the form $\frac{\lambda \,(x\cdot{\mathbf d})}{|x|^3}$
are purely angular multiples of radial inverse-square functions; as such, they 
can be regarded as critical due to their lack of 
inclusion in  the Kato class: hence they are highly interesting
 from a mathematical point of view.
In addition, they share many, but not all, features with the isotropic
inverse square radial potentials:  in particular,  having  the same order of homogeneity, 
they satisfy a Hardy-type inequality (see, for instance \cite{terracini}).

A rich literature deals with Schr\"odinger equations and operators
with isotropic Hardy-type singular potentials, both in the case of one
pole, see e.g. \cite{AFP, egnell, FG, GP, Jan, RW, SM, terracini}, and
in that of multiple singularities, see \cite{caohan, chen, duyckaerts,
  esteban, FMT1, FT, FT2}.  In contrast, only a few papers deal with
the case of anisotropic potentials; in \cite{FMT2} the authors proved
an asymptotic formula for solutions to equation associated with
dipole-type Schr\"odinger operators near the singularity. This
asymptotic analysis will play a crucial role in the discussion of many
fundamental properties of Schr\"odinger operators of the form
\eqref{eq:17}, such as positivity, essential self-adjointness, and
spectral other properties, following the techniques developed in
\cite{FMT1} for Schr\"odinger operators with multipolar inverse-square
potentials.

A natural question is about the effect of the  
configuration of singularities and
the  orientations of dipoles on  the positivity of the associated Schr\"odinger operator. 
The quadratic form associated with the operator  
(\ref{eq:17}) is, denoting
\[
\{\mathcal{L},\mathcal{D},\mathcal{A}\}=\{\lambda_1,\dots,\lambda_k,{\mathbf
  d}_1,\dots,{\mathbf d}_k,a_1,\dots,a_k\},
\]
\begin{align*}
Q_{\mathcal{L},\mathcal{D},\mathcal{A}}:\Di\to\R,\qquad
Q_{\mathcal{L},\mathcal{D},\mathcal{A}}(u):=\int_{\R^N} |\n
u(x)|^2dx-{\displaystyle{\sum_{i=1}^k}}
\dfrac{\lambda_i\,(x-a_i)\cdot{\mathbf d}_i}{|x-a_i|^3}\,u^2(x)\,dx,
\end{align*}
where $\Di$ is the functional space given by  the
completion of $C^\infty_{\rm c}(\R^N)$ with respect to the Dirichlet
norm
$$
\|u\|_{\Di}:=\bigg(\int_{\R^N}|\n u(x)|^2\,dx\bigg)^{1/2}.
$$ 
We recall  that a quadratic form $Q:\Di\to\R$ is
said to be {\em positive definite} if 
$$
\inf_{\Di\setminus\{0\}}\frac{Q(u)}{\|u\|_{\Di}^2}>0.
$$
In the case of a simple dipole operator $-\D-\frac{\lambda
  \,(x\cdot{\mathbf d})}{|x|^3}$, the positivity only depends on the value
of $\lambda $ with respect to the threshold
\begin{equation}
\label{eq:bound}
\Lambda_N:=\sup_{u\in\Di\setminus\{0\}}\dfrac{{\displaystyle
{\int_{\R^N}{\dfrac{x\cdot{\mathbf
d}}{|x|^3}\,u^2(x)\,dx}}}}{{\displaystyle{\int_{\R^N}{|\n u(x)|^2\,dx}}}}.
\end{equation} 
We notice that, by rotation
invariance, $\Lambda_N$ does not depend on the unit vector 
${\mathbf d}$ and, by classical Hardy's 
inequality, $\Lambda_N<4/(N-2)^2$. In particular $\Lambda_N$ is the best constant 
in the following Hardy type inequality 
\begin{equation}\label{eq:hardytype}
\int_{\R^N}{\dfrac{x\cdot{\mathbf d}}{|x|^3}\,u^2(x)\,dx}\leq\Lambda_N
\int_{\R^N}|\n u(x)|^2\,dx\quad\text{for all }u\in\Di 
\text{ and for any unit vector }{\mathbf d}.
\end{equation}
Some numerical approximations of 
  $\Lambda_N$ can be found in \cite[Table 1]{FMT2}. 

It is easy to verify that the quadratic form associated to $-\D-\frac{\lambda
  \,(x\cdot{\mathbf d})}{|x|^3}$
is positive definite in $\Di$ if and only if $|\lambda |<\Lambda_N^{-1}$.
Furthermore, the positivity condition for a one dipole operator can be expressed 
as a condition on the first eigenvalue of the angular component of the operator 
on the unit sphere $\mathbb S^{N-1}$. Indeed, letting 
\begin{equation*}\label{firsteig}
\mu_1^\lambda=\min_{\psi\in H^1(\mathbb
S^{N-1})\setminus\{0\}}\frac{\int_{\mathbb S^{N-1}}|\n_{\mathbb
S^{N-1}}\psi(\theta)|^2\,dV(\theta)-\lambda \int_{\mathbb S^{N-1}}(\theta\cdot
{\mathbf d})\psi^2(\theta)\,dV(\theta)}{\int_{\mathbb S^{N-1}}\psi^2(\theta)\,dV(\theta)}
\end{equation*}
be the first eigenvalue of the
operator $-\D_{\mathbb S^{N-1}}-\lambda\,(\theta\cdot {\mathbf d})$ on $\mathbb S^{N-1}$,
it was proved in \cite[Lemma 2.5]{FMT2} that 
 the quadratic form associated to $-\D-\frac{\lambda
  \,(x\cdot{\mathbf d})}{|x|^3}$
 is positive definite if and only if $\mu_1^\lambda>-\big(\frac{N-2}2\big)^{2}$.

The analysis of the spectral properties of Schr\"odinger operators with multiple 
isotropic inverse square singularities performed in \cite{FMT1} highlighted
how the positivity of the associated quadratic form depends on the 
location and the strength of singularities. In the case of 
 multiple anisotropic singularities, the problem of positivity becomes a
more delicate issue, being the interaction between two dipoles strongly
affected by their mutual orientation. Unlike the isotropic case in which the interaction 
between two poles is either attractive or repulsive depending on the sign 
of coefficients, in the anisotropic one the constructive or destructive 
character of the interaction is determined by the mutual position and orientation.
As a consequence, in contrast with the isotropic
case, it is possible to orientate the dipoles in such a 
way that the interaction  is quite strong even if they are very far away from each other.

The following proposition yields  a
sufficient condition on the magnitudes for the quadratic form to be positive definite
for any localization and orientation of the dipoles.

\begin{Proposition}\label{p:posde}
  A sufficient condition for $Q_{\mathcal{L},\mathcal{D},\mathcal{A}}$ to be positive definite
  for every choice of $\mathcal{A}=\{a_1,a_2,\dots,a_k\}$ and $\mathcal{D}=\{{\mathbf
    d}_1,\dots,{\mathbf d}_k\}$ is that
$$
\sum_{i=1}^k\lambda_i<\Lambda_N^{-1}.
$$
Conversely, if $\sum_{i=1}^k\lambda_i>\Lambda_N^{-1}$ then there
exists a configuration of dipoles $\{ \mathcal{A},\mathcal{D}\}$ such
that $Q_{\mathcal{L},\mathcal{D},\mathcal{A}}$ is not positive
definite.
\end{Proposition}

In this paper we deal with a more general class of
Schr\"odinger operators with locally anisotropic  inverse-square
singularities 
including those with dipole-type potentials introduced in
(\ref{eq:17}). More precisely, we are interested in  operators 
with potentials exhibiting many singularities which are
locally  $L^{\infty}$-angular multiples of radial inverse-square potentials. 

For any $h\in L^{\infty}\big({\mathbb S}^{N-1}\big)$, let  
$\mu_1(h)$ be the first eigenvalue of the
operator $-\D_{\mathbb S^{N-1}}-h(\theta)$ on~$\mathbb S^{N-1}$,
i.e. 
$$
\mu_1(h)=\min_{\psi\in H^1(\mathbb
S^{N-1})\setminus\{0\}}\frac{\int_{\mathbb S^{N-1}}|\n_{\mathbb
S^{N-1}}\psi(\theta)|^2\,dV(\theta)-\int_{\mathbb S^{N-1}}h(\theta)
\psi^2(\theta)\,dV(\theta)}{\int_{\mathbb S^{N-1}}\psi^2(\theta)\,dV(\theta)}.
$$
We  recall  that  $\mu_1(h)$ is  simple and attained by  
a smooth positive eigenfunction $\psi_1^{h}$
such that $\min_{\mathbb S^{N-1}}\psi_1^{h}>0$. Moreover, 
 if, for some $\lambda\in\R$, $h(\theta)=\lambda$ for a.e.
  $\theta\in \mathbb S^{N-1}$, then  $\mu_1(h)=-\lambda$.
On the other hand,  if $h$ is not constant, then 
$$
-\supess_{\mathbb S^{N-1}}
  h<\mu_1(h)<-\media_{{\mathbb S}^{N-1}} h(\theta)\,dV(\theta),
$$
see \cite{FMT2}. 
The quadratic form associated to $-\D
-\frac{h(x/|x|)}{|x|^2}$
 is positive definite if and only if 
\begin{equation}\label{eq:lnh}
\Lambda_N(h):=\sup_{u\in\Di\setminus\{0\}}\dfrac{{\displaystyle
{\int_{\R^N}{\dfrac{h(x/|x|)}{|x|^2}\,u^2(x)\,dx}}}}{{\displaystyle{\int_{\R^N}{|\n
u(x)|^2\,dx}}}}<1,
\end{equation}
or, equivalently, if
and only if $\mu_1(h)>-(N-2)^2/4$, see \cite[Lemma 2.5]{FMT2}.
Furthermore, it is easy to verify that 
$$
\Lambda_N(h)\geq 0\quad\text{for all }h\in L^{\infty}\big({\mathbb S}^{N-1}\big)
$$
and
$$\Lambda_N(h)= 0 \quad\text{if and only if}\quad h\leq 0 \text{ a.e. on }
\mathbb S^{N-1}.
$$
A necessary condition on the angular 
coefficients  for positivity of
the quadratic form associated with multiple dipole-type potentials  for at least a
configuration of singularities is that each single dipole-type local subsystem is 
positive definite, as the following proposition clarifies.

\begin{Proposition}\label{t:mainresult}
  Let $h_1,\dots,h_k,h_\infty\in  L^{\infty}\big({\mathbb S}^{N-1}\big)$,
 $W\in L^{\frac N2}(\R^N)\cap L^{\infty}(\R^N)$, and $R,r_i\in\R^+$, $i=1,\dots,k$. 
  If there exists a configuration of poles
  $\{a_1,\dots,a_k\}$ such that $a_i\neq a_j$ for $i\neq j$ and the quadratic form
\begin{align*}
u\in\Di\mapsto &\int_{\R^N} |\n
u(x)|^2dx-{\displaystyle{\sum_{i=1}^k}}\int_{B(a_i,r_i)}
\dfrac{h_i\big(\frac{x-a_i}{|x-a_i|}\big)}{|x-a_i|^2}\,u^2(x)\,dx
\\
&\quad-\int_{\R^N\setminus B(0,R)}
\frac{h_\infty\big(\frac{x}{|x|}\big)}{|x|^2}\,u^2(x)\,dx-\int_{\R^N}
W(x)\,u^2(x)\,dx
\end{align*}
is positive definite, then 
\begin{equation}\label{eq:19}
\mu_1(h_i)>-\frac{(N-2)^2}4
 \quad \text{ for any }i=1,\dots,k,\infty.
\end{equation}
\end{Proposition}
By virtue of Proposition \ref{t:mainresult}, the following class 
of anisotropic multiple inverse square potentials provides a
suitable framework for the 
analysis of coercivity conditions for Schr\"odinger dipole-type operators:
\begin{align*}
{\mathcal V}:=\left\{
\begin{array}{ccc}
V(x)={\displaystyle{\sum_{i=1}^k\alchi_{B(a_i,r_i)}(x)\frac{
h_i\big(\frac{x-a_i}{|x-a_i|}\big)}{|x-a_i|^2}
      +\alchi_{\R^N\setminus B(0,R)}(x)
\frac{
h_\infty\big(\frac{x}{|x|}\big)}{|x|^2}+W(x)}}:\
  k\in\N,\\[20pt] r_i,R\in\R^+, 
  a_i\in\R^N,\ a_i\neq a_j\ \text{for }i\neq j,\ W\in L^{N/2}(\R^N)\cap L^{\infty}(\R^N),
  \\[10pt]
  h_i\in L^{\infty}\big({\mathbb S}^{N-1}\big),\ \mu_1(h_i)>-\frac{(N-2)^2}4
 \text{ for any $i=1,\dots,k,\infty$}
\end{array}
\right\}.
\end{align*}
If $V=\sum_{i=1}^k\frac{\lambda_i\,(x-a_i)\cdot{\mathbf
d}_i}{|x-a_i|^3}$, then $V\in\mathcal V$ with 
$h_i(\theta)=\lambda_i\theta\cdot {\mathbf d}_i$ and $h_\infty(\theta)
=\theta\cdot\big(\sum_{i=1}^k\lambda_i{\mathbf d}_i\big)$.

By Hardy's and Sobolev's inequalities, it follows that, for any
$V\in\mathcal V$, the first eigenvalue $\mu(V)$ of the operator
$-\D-V$ in $\Di$ is finite, namely 
\begin{equation}\label{eq:30}
\mu(V):=\inf_{u\in\Di\setminus\{0\}}\frac{\displaystyle\int_{\R^N}
\big(|\n u(x)|^2-V(x)u^2(x)\big)\,dx}
{\displaystyle\int_{\R^N}|\n u(x)|^2\,dx}\,>-\infty. 
\end{equation}
We notice that, in view of Sobolev-type embeddings in Lorentz spaces, 
(\ref{eq:30}) holds also for any potential $V$ lying in the Marcinkiewicz space
$L^{N/2,\infty}$.

If the potentials are supported in sufficiently small neighborhoods of 
singularities, then condition~(\ref{eq:19}) turns out to be also sufficient 
for positivity.

\begin{Lemma}{\bf [\,Shattering of singularities\,]}\label{l:sep}
Let  $a_1,a_2,\dots,a_k\in\R^N$, $a_i\neq
a_j$ for $i\neq j$, and $h_1,\dots,h_k,h_\infty
\in L^{\infty}\big({\mathbb S}^{N-1}\big)$ with
$\mu_1(h_i)>-(N-2)^2/4$, for $i=1,\dots,k,\infty$. Then  there exist
${\mathcal U}_1,\dots,{\mathcal U}_k,{\mathcal U}_\infty\subset\R^N$
such that ${\mathcal U}_i$ is a neighborhood of $a_i$ for every 
$i=1,\dots,k$, ${\mathcal U}_\infty$ is a neighborhood of $\infty$, and 
the quadratic form associated to the operator
\begin{align*}
-\Delta
-\sum_{i=1}^k\alchi_{{\mathcal U}_i}(x)
\frac{h_i\big(\frac{x-a_i}{|x-a_i|}\big)}{|x-a_i|^2}
-\alchi_{{\mathcal U}_\infty}(x)
\frac{h_\infty\big(\frac{x}{|x|}\big)}{|x|^2}
\end{align*}
is positive definite.
\end{Lemma}

An analogous  result will be proved also for 
 potentials with infinitely many  dipole-type
 singularities localized in sufficiently small neighborhoods of 
 equidistanced poles, see Lemma~\ref{l:ret}.

 Lemma \ref{l:sep} and Proposition \ref{t:mainresult} establish an
 equivalence between condition (\ref{eq:19}) and the property of being
 compact perturbations of positive operators, as stated in the
 following theorem.

\begin{Theorem}\label{l:semibounded}
For $h_1,\dots,h_k,h_\infty\in  L^{\infty}\big({\mathbb S}^{N-1}\big)$,
 $W\in L^{\frac N2}(\R^N)\cap L^{\infty}(\R^N)$, and $R,r_i\in\R^+$, $i=1,\dots,k$,
let 
\begin{equation}\label{eq:31}
V(x)=\sum_{i=1}^k\alchi_{B(a_i,r_i)}(x)\frac{
h_i\big(\frac{x-a_i}{|x-a_i|}\big)}{|x-a_i|^2}
      +\alchi_{\R^N\setminus B(0,R)}(x)
\frac{
h_\infty\big(\frac{x}{|x|}\big)}{|x|^2}+W(x).
\end{equation}
Then (\ref{eq:19}) is satisfied if and only if there  exists  $\widetilde W\!\in\!
L^{N/2}(\R^N)\cap L^{\infty}(\R^N)$ such that $\mu(V-\widetilde W)>0$. 
\end{Theorem}

By Theorem \ref{l:semibounded},  
Schr\"odinger operators with
potentials in ${\mathcal V}$ are semi-bounded in $L^2(\R^N)$, i.e.
$$
\nu_1(V):=\inf_{u\in
H^1(\R^N)\setminus\{0\}}\frac{\displaystyle\int_{\R^N}\big(|\n
u(x)|^2-V(x)u^2(x)\big)\,dx}{\displaystyle
\int_{\R^N}|u(x)|^2\,dx}\geq - \|\widetilde
W\|_{L^{\infty}(\R^N)}>-\infty,
$$
thus the class ${\mathcal V}$ provides a quite natural setting to
study the spectral properties of multisingular dipole Schr\"odinger
operators in $L^2(\R^N)$. Actually, condition (\ref{eq:19})
characterizing ${\mathcal V}$ is slightly stronger than
semi-boundedness; indeed the operator
$-\Delta-\big(\frac{N-2}2\big)^2$ provides an example of a
$L^2(\R^N)$ semi-bounded operator violating the strict inequality in (\ref{eq:19}). 
On the other hand, we notice that  any   semi-bounded operator in $L^2(\R^N)$ with  a potential  
of the form (\ref{eq:31}) satisfies a weaker condition than (\ref{eq:19}), namely  
$\mu_1(h_i)\geq -\frac{(N-2)^2}4$ for any 
$i=1,\dots,k,\infty$.

The analysis of stability of positivity of Schr\"odinger operators leads
to the problem of {\emph{localization of binding}}  raised by Sigal and 
Ouchinnokov \cite{SO}: if 
$-\Delta-V_1$ and $-\Delta-V_2$ are positive operators, is $-\Delta-V_1-V_2(\cdot-y)$
positive for $|y|$ large?
An affirmative answer to the above question can be found in 
 \cite{Simon} for  compactly supported
potentials and in  \cite{pinchover95} for potentials in  the Kato class.
 When dealing with potentials with an inverse square  singularity at infinity,
the problem becomes more delicate, due to the interaction of singularities
which overlap at infinity and  a {\emph{localization of binding}} type result
requires the additional assumption of some  control of the resulting 
singularity at infinity. Indeed, if, for $j=1,2$, 
$$
V_j=
\sum_{i=1}^{k_j}\alchi_{B(a_i^j,r_i^j)}(x)\frac{
h_i^j\big(\frac{x-a_i^j}{|x-a_i^j|}\big)}{|x-a_i^j|^2}
      +\alchi_{\R^N\setminus B(0,R_j)}(x)
\frac{
h_\infty^j\big(\frac{x}{|x|}\big)}{|x|^2}+W_j(x)\in{\mathcal V},
$$
then 
a necessary condition for positivity of $-\Delta-V_1-V_2(\cdot-y)$ for some $y$ is 
that 
\begin{equation}\label{eq:locbind1}
\mu_1(h_\infty^1+h_\infty^2)>
-\bigg(\frac{N-2}2\bigg)^{\!\!2},
\end{equation}   
see  Proposition \ref{p:nclb}. In
\cite{FMT1}, the authors proved that assumption (\ref{eq:locbind1}) is
also sufficient for {\emph{localization of binding}} when
singularities are locally isotropic, see Theorem~\ref{t:scattering}.
It is worthwhile noticing that a strong lack of isotropy can cause the failure of localization of
binding even under assumption (\ref{eq:locbind1}); in section
\ref{sec:localization-binding}, we will construct two anisotropic potentials in the class
${\mathcal V}$ satisfying (\ref{eq:locbind1}) for which no localization of
binding result holds true. 

\begin{Example}\label{ex:controes}
{\it For $N\geq 4$, there exist $V_1,V_2\in{\mathcal V}$ such that 
$\mu(V_1),\mu(V_2)>0$,  \eqref{eq:locbind1} holds, and for every $R>0$ there
exists $y_R\in\R^N$ such that $|y_R|>R$ and the quadratic form associated to 
the operator
$-\Delta-\left(V_1+V_2(\cdot-y_R)\right)$ is not positive semidefinite, i.e.
$\mu(V_1+V_2(\cdot-y_R))<0$.  }
\end{Example}

On the other 
hand, it is still possible to prove the following   {\emph{localization of binding}}
type result under a stronger control on the singularities at infinity than 
(\ref{eq:locbind1}).

\begin{Theorem}\label{t:localization} 
 Let 
\begin{align*}
&V_1(x)=
\sum_{i=1}^{k_1}\alchi_{B(a_i^1,r_i^1)}(x)\frac{
h_i^1\big(\frac{x-a_i^1}{|x-a_i^1|}\big)}{|x-a_i^1|^2}
      +\alchi_{\R^N\setminus B(0,R_1)}(x)
\frac{
h_\infty^1\big(\frac{x}{|x|}\big)}{|x|^2}+W_1(x)\in{\mathcal V},\\ 
&V_2(x)=\sum_{i=1}^{k_2}\alchi_{B(a_i^2,r_i^2)}(x)\frac{
h_i^2\big(\frac{x-a_i^2}{|x-a_i^2|}\big)}{|x-a_i^2|^2}
      +\alchi_{\R^N\setminus B(0,R_2)}(x)
\frac{
h_\infty^2\big(\frac{x}{|x|}\big)}{|x|^2}+W_2(x)\in{\mathcal V}.
\end{align*}
Assume that $\mu(V_1),\mu(V_2)>0$, and $\supess_{{\mathbb
    S}^{N-1}}(h_1^\infty)^+ +\supess_{{\mathbb S}^{N-1}}(h_2^\infty)^+
<(N-2)^2/4$.  Then, there exists $R>0$ such that
$\mu(V_1+V_2(\cdot-y))>0$ for every $y\in\R^N$ with $|y|\geq R$.
\end{Theorem}

A further key property of Schr\"odinger operators which turns out to
be very sensitive to the presence of singular terms is 
the {\em essential self-adjointness},
namely the existence of a unique self-adjoint extension. Semi-bounded
Schr\"odinger operators are essentially self-adjoint whenever the
potential is not too singular (see \cite{reedsimon}).  
On the other hand,  inverse square potentials exhibit a  quite strong 
singularity which makes the problem of  essential self-adjointness nontrivial.
We mention that essential self-adjointness in the case of Hardy type potentials was
discussed in \cite{kwss} for the one-pole case and in \cite{FMT1} 
for   many poles.  
The following theorem provides a necessary and sufficient condition 
on the magnitudes of dipole moments  for the essential self-adjointness of
multisingular dipole Schr\"odinger operators. An extension to the case
of infinitely many  dipole-type
 singularities distributed on reticular structures is contained in Theorem 
\ref{t:self-adjo-ret}.

\begin{Theorem}\label{t:self-adjo}
Let 
$$
V(x)=\displaystyle{\sum_{i=1}^k\alchi_{B(a_i,r_i)}(x)
\frac{h_i\big(\frac{x-a_i}{|x-a_i|}\big)}{|x-a_i|^2}
+\alchi_{\R^N\setminus B(0,R)}(x)
\frac{h_\infty\big(\frac{x}{|x|}\big)}{|x|^2}+W(x)}
\in{\mathcal V}.
$$
Then the Schr\"odinger operator $-\D-V$  is essentially self-adjoint
in $C^{\infty}_{\rm c}(\R^N\setminus\{a_1,\dots,a_k\})$ if and only if
$\mu_1(h_i)\geq-\big(\frac{N-2}{2}\big)^{\!2}+1$, 
for all $i=1,\dots,k$. 
\end{Theorem}

The proof of the above theorem is based on the
asymptotic analysis performed in  \cite{FMT2}, where the exact behavior 
 near the poles of solutions 
to Schr\"odinger equations with dipole-type singular
potentials is evaluated. 
From Theorem \ref{t:self-adjo}, it follows that, if
$V\in{\mathcal V}$ with $\mu_1(h_i)\geq-\big(\frac{N-2}{2}\big)^{\!2}+1$,
 for all
$i=1,\dots,k$, then the {\em {Friedrichs
extension}} $(-\D-V)^F$ defined as 
\begin{align}\label{eq:52}
D\big((-\D-V)^F\big)&=\{u\in H^1(\R^N):\-\D u-Vu\in L^2(\R^N)\},
\quad u\mapsto-\D u-Vu,
\end{align}
is the unique self-adjoint extension. On the other hand, if 
 $\mu_1(h_i)<-\big(\frac{N-2}{2}\big)^{\!2}+1$
 for some $i$, then $-\D-V$  admits  many self-adjoint extensions,
among which the Friedrichs extension is the only one whose domain is
included in $H^1(\R^N)$, namely it is the unique self-adjoint extension
to which we can associate a natural  quadratic form.

The paper is organized as follows. Section \ref{sec:proof-prop-refp:p}
contains the proofs of Propositions \ref{p:posde} and
\ref{t:mainresult}.  In section \ref{sec:bound} we prove a positivity
criterion in the spirit of the Allegretto-Piepenbrink theory, which is
used to prove Lemma \ref{l:sep} and its reticular version (see Lemma
\ref{l:ret}); a key tool in the proof of Lemma \ref{l:sep} is the
analysis of positivity of potentials obtained as juxtaposition of
potentials with different singularity rate, which is a nontrivial
issue due to the lack of isotropy, see Lemmas \ref{l:intercapedine},
\ref{l:sha1}, and \ref{l:sep_variante}.  In
section \ref{sec:pertinf}  we discuss the stability of positivity
with respect to perturbations of the potentials with  singularities
localized at dipolar-shaped neighborhoods either of a dipole or of
infinity. In section \ref{sec:localization-binding} we prove Theorem
\ref{t:localization} and show that condition \eqref{eq:locbind1} is no
sufficient for localization of binding by constructing a suitable example.
Section \ref{sec:semi-bound-essent} is devoted to the proof of Theorem \ref{t:self-adjo}
and of its reticular counterpart.

\medskip
\noindent
{\bf Notation. } We list below some notation used throughout the
paper.\par
\begin{itemize}
\item[-]$B(a,r)$ denotes the ball $\{x\in\R^N: |x-a|<r\}$ in $\R^N$ with
center at $a$ and radius $r$.
\item[-] $\R^+:=(0,+\infty)$ is the half line of positive real numbers.
\item[-] For any $A\subset \R^N$, $\alchi_{A}$ denotes the
characteristic function of $A$. 
\item[-] $S$ is the best constant in the Sobolev inequality
$S\|u\|_{L^{2^*}(\R^N)}^2\leq \|u\|_{{\mathcal D}^{1,2}(\R^N)}^2$.
\item[-] For all $t\in\R$,  $t^+:=\max\{t,0\}$ (respectively
$t^-:=\max\{-t,0\}$) denotes the positive (respectively negative)
part of $t$.
\item[-] For all functions $f:\ \R^N\to\R$, $\mathop{\rm supp}f$
denotes the support of $f$, i.e. the closure of the set of points
where $f$ is non zero.
\item[-] $\omega_N$ denotes the volume of the unit ball in $\R^N$.
\item[-] For any open set $\Omega\subset\R^N$, ${\mathcal D}'(\Omega)$
denotes the space of distributions in $\Omega$.
\item[-] For any $f\in L^{\infty}(A)$, with either $A={\mathbb S}^{N-1}$
or $A=\R^N$, 
  we denote the essential supremum of $f$ in $A$ as 
 $\supess_{A}f:=\inf\{\a\in\R: f(x)\leq \alpha \text{
    for a.e. }x\in A\}$, while the essential infimum
  of $f$ in $A$ is denoted as $\infess_{A}f:=-\supess_{A}(-f)$.
\end{itemize}

\section{Proof of Propositions \ref{p:posde} 
and  \ref{t:mainresult}}\label{sec:proof-prop-refp:p}

This section is devoted to the proof of Propositions \ref{p:posde}
and  \ref{t:mainresult}.

\bigskip\noindent

\begin{pfn}{Proposition \ref{p:posde}}
From \eqref{eq:bound}, it follows that
\begin{equation}\label{eq:3}
Q_{\mathcal{L},\mathcal{D},\mathcal{A}}(u)\geq
  \bigg(1-\Lambda_N\sum_{i=1}^k{\lambda_i}
\bigg) \int_{\R^N} |\n 
  u|^2\,dx. 
\end{equation}
Hence a sufficient condition for $Q_{\mathcal{L},\mathcal{D},\mathcal{A}}$ to be 
positive definite is that 
$$
\sum_{i=1}^k{\lambda_i}<\Lambda_N^{-1}.
$$
Assume now that $\sum_{i=1}^k{\lambda_i}>\Lambda_N^{-1}$ and fix    
${\mathbf d}\in \R^N$, $|{\mathbf d}|=1$.  From \eqref{eq:bound} and density
of $C^{\infty}_{\rm c}(\R^N)$ in $\Di$, there exists some function
$\phi\in C^{\infty}_{\rm c}(\R^N)$ such that
$$
\int_{\R^N}|\nabla\phi|^2\,dx
-\bigg(\sum_{i=1}^k\lambda_i\bigg)\int_{\R^N}\frac{x\cdot{\mathbf d}}{|x|^3}\,
\phi^2\,dx < 0.
$$
Let $a_1,\dots,a_k\in\R^N$  and set ${\mathcal
  A}=\{a_1,\dots,a_k\}\subset\R^N$, ${\mathcal
  D}=\{\mathbf d,\dots,\mathbf d\}$.
For any $\mu>0$, consider the function $\phi_{\mu}(x)=
\mu^{-\frac{N-2}2}\phi(x/\mu)$. A change of variable yields
$$
\int_{\R^N}|\nabla\phi_{\mu}|^2\,dx
-\sum_{i=1}^k\lambda_i\int_{\R^N}\frac{(x-a_i)\cdot{\mathbf d}}
{|x-a_i|^3}\,\phi_{\mu}^2\,dx=
\int_{\R^N}|\nabla\phi|^2\,dx
-\sum_{i=1}^k\lambda_i\int_{\R^N}\frac{x\cdot{\mathbf d}}
{|x|^3}\,\phi^2\Big(x+\frac{a_i}{\mu}\Big)\,dx
$$
for all $\mu>0$.
Letting $\mu\to\infty$, the Dominated Convergence Theorem 
yields
$$
\int_{\R^N}|\nabla\phi_{\mu}|^2\,dx
-\sum_{i=1}^k\lambda_i\int_{\R^N}\frac{(x-a_i)\cdot{\mathbf d}}
{|x-a_i|^3}\,\phi_{\mu}^2\,dx
\longrightarrow \int_{\R^N}|\nabla\phi|^2\,dx
-\bigg(\sum_{i=1}^k\lambda_i\bigg)\int_{\R^N}\frac{x\cdot{\mathbf d}}{|x|^3}\,
\phi^2\,dx<0
$$
therefore  
$Q_{\mathcal{L},\mathcal{D},\mathcal{A}}(\phi_{\mu})<0$ for $\mu$ sufficiently large,
thus proving the second part of Proposition \ref{p:posde}.
\end{pfn} 

\begin{pfn}{Proposition \ref{t:mainresult}}
Assume that, for some configuration $\{a_1,\dots, a_k\}$, for some $\e>0$,
and for any $u\in\Di$,
\begin{align}
\int_{\R^N}|\n u(x)|^2dx&-{\displaystyle{\sum_{i=1}^k}}\int_{B(a_i,r_i)}
\dfrac{h_i\big(\frac{x-a_i}{|x-a_i|}\big)}{|x-a_i|^2}\,u^2(x)\,dx\\
&-\int_{\R^N\setminus B(0,R)}
\frac{h_\infty\big(\frac{x}{|x|}\big)}{|x|^2}\,u^2(x)\,dx-\int_{\R^N}
W(x)\,u^2(x)\,dx
\geq\e\int_{\R^N}|\n u(x)|^2\,dx.\nonumber
\end{align}
Arguing by contradiction, suppose that, for some $i=1,\dots,k$,
$\mu_1(h_i)\leq-(N-2)^2/4$, or equivalently that
$\La_N(h_i)\geq1$, with $\La_N(h_i)$ as in \eqref{eq:lnh}. Let
$0<\delta<\e\Lambda_N(h_i)^{-1}$.  By \eqref{eq:lnh} and density
of $C^{\infty}_{\rm c}(\R^N)$ in $\Di$,
there exists $\phi\in C^{\infty}_{\rm
  c}(\R^N)$ such that
\begin{equation}\label{eq:prs}
  \int_{\R^N}|\n\phi(x)|^2dx-\int_{\R^N}
\dfrac{h_i\big(\frac{x}{|x|}\big)}{|x|^2}\,\phi^2(x)\,dx
  <\delta\int_{\R^N}\dfrac{h_i\big(\frac{x}{|x|}\big)}{|x|^2}\,\phi^2(x)\,dx.
\end{equation}
The rescaled function $\phi_{\mu}(x)=\mu^{-(N-2)/2}\phi(x/\mu)$ satisfies
\begin{align}
&\int_{\R^N}|\n \phi_{\mu}(x-a_i)|^2dx-{\displaystyle{\sum_{j=1}^k}}\int_{B(a_j,r_j)}
\dfrac{h_j\big(\frac{x-a_j}{|x-a_j|}\big)}{|x-a_j|^2}\,\phi_{\mu}^2(x-a_i)\,dx\\
&\quad-\int_{\R^N\setminus B(0,R)}
\frac{h_\infty\big(\frac{x}{|x|}\big)}{|x|^2}\,\phi_{\mu}^2(x-a_i)\,dx-\int_{\R^N}
W(x)\,\phi_{\mu}^2(x-a_i)\,dx\nonumber\\
&=\int_{\R^N}|\n\phi(x)|^2dx-\int_{B\big(0,\frac{r_i}\mu\big)}
\dfrac{h_i\big(\frac{x}{|x|}\big)}{|x|^2}\,\phi^2(x)\,dx
-\int_{\R^N\setminus B\big(\frac{-a_i}\mu,\frac{R}\mu\big)}
\frac{h_\infty\big(\frac{x+a_i/\mu}{|x+a_i/\mu|}\big)}{|x+a_i/\mu|^2}\,\phi^2(x)\,dx
\nonumber\\
&\quad-\sum_{j\neq i}\int_{B\big(\frac{a_j-a_i}\mu,\frac{r_j}\mu\big)}
\frac{h_j\big(\frac{x-(a_j-a_i)/\mu}{|x-(a_j-a_i)/\mu|}\big)}
{\big|x-(a_j-a_i)/\mu\big|^2}\,\phi^2(x)\,dx-
\mu^2\int_{\R^N}W(\mu x+a_i)\,\phi^2(x)\,dx\nonumber\\
&=\int_{\R^N}|\n\phi(x)|^2dx-
\int_{\R^N}\dfrac{h_i\big(\frac{x}{|x|}\big)}{|x|^2}\,\phi^2(x)\,dx
+o(1),\quad\text{ as }\mu\to0.\nonumber
\end{align}
Letting $\mu\to0$, by \eqref{eq:prs} and \eqref{eq:lnh}, we obtain 
\begin{align*}
  \e\int_{\R^N}|\n\phi(x)|^2dx &\leq\int_{\R^N}|\n\phi(x)|^2dx
-\int_{\R^N}\dfrac{h_i\big(\frac{x}{|x|}\big)}{|x|^2}\,\phi^2(x)\,dx\\
&<\delta
\int_{\R^N}\dfrac{h_i\big(\frac{x}{|x|}\big)}{|x|^2}\,\phi^2(x)\,dx
\leq\delta\Lambda_N(h_i)\int_{\R^N}|\n \phi(x)|^2dx
\end{align*}
thus giving rise to a contradiction.

Suppose now that, $\mu_1(h_\infty)\leq-\frac{(N-2)^2}{4}$, or equivalently that 
$\La_N(h_\infty)\geq1$,
with $\La_N(h_\infty)$ defined in \eqref{eq:lnh}.  
Let $\delta\in(0,\e\Lambda_N(h_\infty)^{-1})$. By definition of $\La_N(h_\infty)$ and
density of $C^{\infty}_{\rm c}(\R^N\setminus\{0\})$ in $\Di$, 
there exists $\phi\in
C^{\infty}_{\rm c}(\R^N\setminus\{0\})$ such that  
$$
\int_{\R^N}|\n\phi(x)|^2dx-\int_{\R^N}\frac{h_\infty\big(\frac{x}{|x|}\big)}
{|x|^2}\,\phi^2(x)\,dx
<\delta\int_{\R^N}\frac{h_\infty\big(\frac{x}{|x|}\big)}
{|x|^2}\,\phi^2(x)\,dx.
$$
Since $\phi\in C^{\infty}_{\rm c}(\R^N\setminus\{0\})$,
the rescaled function $\phi_{\mu}(x)=\mu^{-(N-2)/2}\phi(x/\mu)$
satisfies
\begin{align*}
\int_{\R^N}|\n\phi_{\mu}(x)|^2dx&-{\displaystyle{\sum_{i=1}^k}}\int_{B(a_i,r_i)}
\dfrac{h_i\big(\frac{x-a_i}{|x-a_i|}\big)}{|x-a_i|^2}\,\phi_{\mu}^2(x)\,dx\\
&\quad\quad-\int_{\R^N\setminus B(0,R)}
\frac{h_\infty\big(\frac{x}{|x|}\big)}{|x|^2}\,\phi_{\mu}^2(x)\,dx-\int_{\R^N}
W(x)\,\phi_{\mu}^2(x)\,dx\nonumber\\
&=\int_{\R^N}|\n\phi(x)|^2dx-{\displaystyle{\sum_{i=1}^k}}
\int_{B\big(\frac{a_i}\mu,\frac{r_i}\mu\big)}
\dfrac{h_i\big(\frac{x-a_i/\mu}{|x-a_i/\mu|}\big)}{|x-a_i/\mu|^2}\,\phi^2(x)\,dx\\
&\quad\quad-\int_{\R^N\setminus B\big(0,\frac{R}\mu\big)}
\frac{h_\infty\big(\frac{x}{|x|}\big)}{|x|^2}\,\phi^2(x)\,dx-\mu^2\int_{\R^N}
W(\mu x)\,\phi^2(x)\,dx\nonumber\\
&=\int_{\R^N}|\n\phi(x)|^2dx-
\int_{\R^N}\frac{h_\infty\big(\frac{x}{|x|}\big)}
{|x|^2}\,\phi^2(x)\,dx+o(1),\quad\text{ as
}\mu\to\infty. 
\end{align*} 
Letting $\mu\to\infty$ and arguing as
above, we obtain easily a contradiction.~\end{pfn}

\section{The Shattering Lemma}\label{sec:bound}

The well-known Allegretto-Piepenbrink theory
\cite{allegretto,piepenbrink} suggests us a 
  criterion for establishing positivity of
Schr\"odinger operators with potentials in $\mathcal V$, by
 relating the existence of positive
solutions to a Schr\"odinger equation with the positivity of the
spectrum of the corresponding operator. For analogous criteria for
potentials in the Kato class we refer to \cite[Theorem 2.12]{CFKS}.

\begin{Lemma}\label{l:positivity_condition}
Let $V\in{\mathcal V}$. 
Then the two following
conditions are equivalent:
\begin{align*}
  (i)\quad &\mu(V):=\inf_{u\in\Di\setminus\{0\}}\frac{\int_{\R^N}
    \big(|\n u(x)|^2-V(x)u^2(x)\big)\,dx}{\int_{\R^N}|\n u(x)|^2\,dx}>0;\\[10pt]
  (ii)\quad &\text{there exist }\e>0 \text{ and }\varphi\in \Di,\
  \varphi>0\text{ in }\R^N\setminus\{a_1,\dots,a_k\},\text{ and
  }\varphi\text{
    continuous  }\\
  &\text{in }\R^N\setminus\{a_1,\dots,a_k\},\text{ such that
  }-\Delta\varphi-V\varphi\geq\e\,V\varphi\text{ in
  }\big(\Di\big)^\star, \text{ i.e. }\\[5pt]
  &
{}_{(\Di)^\star}\langle -\Delta\varphi-V\varphi-\e\,V\varphi,w\rangle_{\Di}
=\int_{\R^N}\big[\nabla\varphi\cdot\nabla w-(1+\e)V\varphi\, w\big]\,dx
\geq 0\\[5pt]
  & \text{for any }w\in\Di \text{ such that }w\geq0 \text{ a.e. in
  }\R^N.
 \end{align*}
Moreover, if (ii) holds, $\mu(V)\geq\frac{\e}{1+\e}$.
\end{Lemma}

\begin{pf}
Let $V(x)=\sum_{i=1}^k\alchi_{B(a_i,r_i)}(x)
\frac{h_i((x-a_i)/|x-a_i|)}{|x-a_i|^2}
+\alchi_{\R^N\setminus B(0,R)}(x)
\frac{h_\infty(x/|x|)}{|x|^2}+W(x)
\in{\mathcal V}$ and set  $\bar h:=\sum_{i=1}^{k}\|h_i\|_{L^\infty({\mathbb S}^{N-1})}+
\|h_\infty\|_{L^\infty({\mathbb S}^{N-1})}$. 
From Hardy's, H\"older's, and Sobolev's inequalities there holds  
\begin{equation}\label{eq:4}
\int_{\R^N}V(x)u^2(x)\,dx\leq\bigg[\frac{4\bar h}{(N-2)^2}
+S^{-1}\|W\|_{L^{N/2}(\R^N)}\bigg]\int_{\R^N}|\n u(x)|^2\,dx,
\end{equation}
for every $u\in\Di$.

Assume that $(i)$ holds. 
If $0<\e<\frac{\mu(V)}{2}\Big[\frac{4\bar h}{(N-2)^2}+S^{-1}\|W\|_{L^{N/2}(\R^N)}\Big]^{-1}$, 
from (\ref{eq:4}) it follows that 
$$
\int_{\R^N}\big(|\n u(x)|^2-(1+\e)V(x)u^2(x)\big)\,dx\geq\frac{\mu(V)}2\int_{\R^N}|\n u(x)|^2\,dx.
$$
As a consequence, for any fixed $p\in L^{N/2}(\R^N)\cap L^{\infty}(\R^N)$, 
$p(x)>0$ a.e. in $\R^N$, the infimum
$$
\nu_p(V+\e V)=\inf_{u\in\Di\setminus\{0\}}\dfrac{{{\int_{\R^N}\big(|\n u(x)|^2-
(1+\e) V(x)\,u^2(x)\big)\,dx}}}{{{\int_{\R^N}p(x)u^2(x)}}}
$$
is strictly positive and  attained by some function $\varphi\in\Di\setminus\{0\}$ satisfying
$$
-\D\varphi(x)-V(x)\varphi(x)=\e\,V(x)\varphi(x)+\nu_p(V+\e V)p(x)\varphi(x).
$$
By evenness we can assume $\varphi\geq0$. Since $V\in{\mathcal V}$, the
Strong Maximum Principle allows us to conclude that $\varphi>0$ in
$\R^N\setminus\{a_1,\dots,a_k\}$, while standard regularity theory
ensures regularity of $\varphi$ outside the poles. Hence $(ii)$
holds. 

Assume now that $(ii)$ holds. For any $u\in C^{\infty}_{\rm
c}(\R^N\setminus\{a_1,\dots,a_k\})$, testing the weak inequality
satisfied by $\varphi$ with $u^2/\varphi$ we get  
\begin{align*}
(1\!+\!\e)\!\int_{\R^N}\!\!V(x)u^2(x)\,dx
&\leq2\!\int_{\R^N}\!\frac{u(x)}{\varphi(x)}\n u(x)\cdot\n
\varphi(x)\,dx-\!\int_{\R^N}\!\frac{u^2(x)}{\varphi^2(x)}|\n
\varphi(x)|^2\,dx\leq\!\int_{\R^N}\!|\n u(x)|^2\,dx.
\end{align*}
By density of $C^{\infty}_{\rm c}(\R^N\setminus\{a_1,\dots,a_k\})$ in
$\Di$ we deduce that,  
for every $u\in\Di\setminus\{0\}$, 
$$\int_{\R^N}V(x)u^2(x)\leq\frac{1}{1+\e}\int_{\R^N}|\n u(x)|^2\,dx,$$
implying
\begin{equation}\label{eq:66}
\int_{\R^N}\big(|\n
u(x)|^2-V(x)u^2(x)\big)\,dx\geq\frac{\e}{1+\e}\int_{\R^N}|\n u(x)|^2\,dx,
\end{equation}
and hence $\mu(V)\geq \frac{\e}{\e+1}>0$.
\end{pf}

The above positivity criterion allows to extend to multiple
dipole Schr\"odinger operators the Shattering Lemma in 
\cite[Lemma 1.3]{FMT1} yielding positivity in the case of
 singularities  localized strictly near the poles. 

Let us start by observing that,
evaluating the quotient minimized in the definition of $\mu(V)$ at 
functions concentrating at the singularities,  
 $\mu(V)$ can be estimated from above as follows.
\begin{Lemma}\label{l:estmu}
For any 
\begin{equation}\label{eq:viinV}
V(x)=\displaystyle{\sum_{i=1}^k\alchi_{B(a_i,r_i)}(x)
\frac{h_i\big(\frac{x-a_i}{|x-a_i|}\big)}{|x-a_i|^2}
+\alchi_{\R^N\setminus B(0,R)}(x)
\frac{h_\infty\big(\frac{x}{|x|}\big)}{|x|^2}+W(x)}
\in{\mathcal V},
\end{equation}
there holds
\begin{align*}
\mu(V)\leq 1-\max\big\{
0,\Lambda_N(h_1),\dots,\Lambda_N(h_k),\Lambda_N(h_\infty)\big\},
\end{align*}
where $\Lambda_N(h_i)$ is defined in \eqref{eq:lnh}.
\end{Lemma}

\begin{pf}
Let us first consider the case $\Lambda_N(h_i)= 0$ for every $i=1,\dots,k,\infty$.
Let us fix
$u\in C^{\infty}_{\rm c}(\R^N)$ and
$P\in \R^N\setminus\{a_1,\dots,a_k\}$. Letting
$u_{\mu}(x)=\mu^{-\frac{N-2}2}u\big(\frac{x-P}\mu\big)$, for 
$\mu$ small there holds 
$$
\mu(V)\leq 1-\frac{\int_{\R^N}W(x)u_{\mu}^2(x)\,dx}
{\int_{\R^N}|\n u_{\mu}(x)|^2\,dx}=1+o(1)\quad\text{as }\mu\to 0^+.
$$
Letting $\mu\to 0^+$ we obtain that $\mu(V)\leq 1$.

Assume now that $\max_{i=1,\dots,k,\infty}\Lambda_N(h_i)>0$. 
  Suppose $\Lambda_N(h_1)\leq\Lambda_N(h_2)\leq\dots\Lambda_N(h_k)$ and let
  $\e>0$.  From (\ref{eq:lnh}) and by density of $C^{\infty}_{\rm
    c}(\R^N)$ in $\Di$, there exists
  $\phi\in C^{\infty}_{\rm c}(\R^N)$ such
  that
\begin{equation*}
\int_{\R^N}|\n\phi(x)|^2dx<\bigg[\Lambda_N(h_k)^{-1}+\e\bigg]
\int_{\R^N}\frac{h_k\big(\frac{x-a_k}{|x-a_k|}\big)}{|x-a_k|^2}\,\phi^2(x)\,dx.
\end{equation*}
Letting $\phi_{\mu}(x)=\mu^{-\frac{N-2}2}\phi\big(\frac{x-a_k}\mu\big)$, 
for any $\mu>0$ there holds
\begin{align*}
\mu(V)&\leq 1-
\frac{\int_{B(a_k,r_k)}h_k\big(\frac{x-a_k}{|x-a_k|}\big)|x-a_k|^{-2}
\phi_{\mu}^2(x)\,dx}{\int_{\R^N}|\n\phi_{\mu}(x)|^2\,dx}-\sum_{i=1}^{k-1}
\frac{\int_{B(a_i,r_i)}h_i\big(\frac{x-a_i}{|x-a_i|}\big)|x-a_i|^{-2}
\phi_{\mu}^2(x)\,dx}{\int_{\R^N}|\n\phi_{\mu}(x)|^2\,dx}\\
&\quad-\frac{\int_{\R^N\setminus B(0,R)}h_\infty\big(\frac{x}{|x|}\big)|x|^{-2}\phi_{\mu}^2(x)\,dx}{\int_{\R^N}|\n \phi_{\mu}(x)|^2\,dx}-\frac{\int_{\R^N}W(x)\phi_{\mu}^2(x)\,dx}
{\int_{\R^N}|\n \phi_{\mu}(x)|^2\,dx}
\\
&=1-\frac{\int_{\R^N}h_k\big(\frac{x}{|x|}\big)|x|^{-2}\phi^2(x)\,dx}
{\int_{\R^N}|\n\phi(x)|^2\,dx}+o(1)\quad\text{as}\quad \mu\to 0^+.
\end{align*}
Letting $\mu\to 0^+$, by the choice of $\phi$ we obtain
$$
\mu(V)\leq 1-\bigg[\Lambda_N(h_k)^{-1}+\e\bigg]^{-1}
$$
for any $\e>0$. Letting $\e\to0$ we derive that
$\mu(V)\leq 1-\Lambda_N(h_k)$. Repeating the same argument, we obtain 
a function 
$\psi_\infty\in C^{\infty}_{\rm
c}(\R^N\setminus\{0\})$ such that  
\begin{equation*}
\int_{\R^N}|\n\psi_\infty(x)|^2dx<\bigg[\Lambda_N(h_\infty)^{-1}+\e\bigg]
\int_{\R^N}\frac{h_\infty\big(\frac{x}{|x|}\big)}{|x|^2}\,\psi_\infty^2(x)\,dx.
\end{equation*}
Setting $\psi_{\mu}(x)=\mu^{-\frac{N-2}2}\psi_\infty\big(\frac{x}\mu\big)$
and letting $\mu\to+\infty$ we obtain also that 
$\mu(V)\leq 1-\Lambda_N(h_\infty)$. The required estimate is thereby proved.
\end{pf}

The proof of Lemma \ref{l:sep} is immediate if
$\Lambda_N(h_i)= 0$ for all $i=1,\dots,k,\infty$, i.e. if $h_i\leq 0$ a.e. 
in  ${\mathbb S}^{N-1}$ for all $i=1,\dots,k,\infty$, as the following
lemma states.

\begin{Lemma}\label{l:sep_parte1}
Let  $a_1,a_2,\dots,a_k\in\R^N$, $a_i\neq
a_j$ for $i\neq j$, and $h_1,\dots,h_k,k_\infty
\in L^{\infty}\big({\mathbb S}^{N-1}\big)$ with
$\mu_1(h_i)>-(N-2)^2/4$, for $i=1,\dots,k,\infty$. Then, if 
$\Lambda_N(h_i)= 0$ for all $i=1,\dots,k,\infty$, for every
${\mathcal U}_1,\dots,{\mathcal U}_k,{\mathcal U}_\infty\subset\R^N$
such that ${\mathcal U}_i$ is a neighborhood of $a_i$ for every 
$i=1,\dots,k$ and ${\mathcal U}_\infty$ is a neighborhood of $\infty$, there holds
\begin{align*}
\mu\bigg(
\sum_{i=1}^k\alchi_{{\mathcal U}_i}(x)
\frac{h_i\big(\frac{x-a_i}{|x-a_i|}\big)}{|x-a_i|^2}
+\alchi_{{\mathcal U}_\infty}(x)
\frac{h_\infty\big(\frac{x}{|x|}\big)}{|x|^2}
\bigg)=1.
\end{align*}
\end{Lemma}
\begin{pf}
The inequality $\mu\big(
\sum_{i=1}^k\alchi_{{\mathcal U}_i}(x)
\frac{h_i((x-a_i)/|x-a_i|)}{|x-a_i|^2}
+\alchi_{{\mathcal U}_\infty}(x)
\frac{h_\infty({x}/{|x|})}{|x|^2}
\big)\leq 1$ follows from Lemma \ref{l:estmu}, whereas the reverse inequality 
comes immediately from the assumption $h_i\leq 0$ a.e. 
in  ${\mathbb S}^{N-1}$ for all $i=1,\dots,k,\infty$.
\end{pf}

When dealing with isotropic potentials, it is quite easy to study the
positivity of potentials obtained as juxtaposition of potentials with
different singularity rate. If, for example,
$\Omega_1\subset\Omega_2$, then the operator
$-\Delta-\frac{\lambda_1}{|x|^2}
\alchi_{\Omega_1}-\frac{\lambda_2}{|x|^2} \alchi_{\Omega_2\setminus
  \Omega_1}$ is positive definite whenever
$\max\{\lambda_1,\lambda_2\}< (N-2)^2/4$. On the other hand, the
positivity of a potential obtained as juxtaposition of two potentials
with changing sign angular components is more delicate to be
established.  The analysis we are going to develop shows how
juxtaposition of potentials giving rise to positive quadratic forms
produces positive operators if their contact region has some
particular shape (which resembles a sphere deformed according to
dipole coefficients).

For any $h\in L^{\infty}\big({\mathbb S}^{N-1}\big)$, let  
$\psi_1^{h}$ denote the positive $L^2$-normalized
eigenfunction associated to the first eigenvalue
$\mu_1(h)$  of the
operator $-\D_{\mathbb S^{N-1}}-h(\theta)$ on~$\mathbb S^{N-1}$.
Let us notice that, for $h\equiv0$, $\psi_1^{0}\equiv 1/\sqrt{\omega_N}$,
 where
$\omega_N$ denotes the volume of the unit sphere ${\mathbb
    S}^{N-1}$, i.e. $\omega_N=\int_{{\mathbb S}^{N-1}}dV(\theta)$.

For $h_1,h_2\in L^{\infty}\big({\mathbb S}^{N-1}\big)$,
$\sigma>0$, and $R>0$, let us denote 
\begin{equation}\label{eq:26}
{\mathcal E}^{\sigma,R}_{h_1,h_2}:=
\bigg\{x\in\R^N:\ |x|<R\bigg(\frac{\psi_1^{h_1}(x/|x|)}
{\psi_1^{h_2}(x/|x|)}\bigg)^{\!\!\sigma}\bigg\}.
\end{equation}
Let us notice that 
$$
B\bigg (0,R\bigg[\min_{\theta\in\mathbb
  S^{N-1}}\frac{\psi_1^{h_1}(\theta)}
{\psi_1^{h_2}(\theta)}\bigg]^\sigma\bigg)\subset
{\mathcal E}^{\sigma,R}_{h_1,h_2} \subset B\bigg (0,R\bigg[\max_{\theta\in\mathbb
  S^{N-1}}\frac{\psi_1^{h_1}(\theta)}
{\psi_1^{h_2}(\theta)}\bigg]^\sigma\bigg).
$$
We also set, for $a\in\R^N$,
$$
{\mathcal E}^{\sigma,R}_{h_1,h_2}(a):=
\big\{x\in\R^N:\ x-a\in {\mathcal E}^{\sigma,R}_{h_1,h_2}\big\}.
$$

\begin{Lemma}\label{l:intercapedine}
Let $h_i\in L^{\infty}\big({\mathbb S}^{N-1}\big)$, 
$\mu_1(h_i)>-\big(\frac{N-2}2\big)^{\!2}$, $i=1,2$, 
$\sigma_1,\sigma_2>0$ such that
\begin{equation}\label{eq:7}
\frac{N-2}2<\frac1{\sigma_1}<\frac1{\sigma_1}+\frac1{\sigma_2}<
\frac{N-2}2+\min_{i=1,2}\sqrt{\frac{(N-2)^2}{4}+\mu_1(h_i)},
\end{equation}
and $R_1,R_2\in(0,+\infty)$ such that ${\mathcal
  E}^{\sigma_1,R_1}_{h_1,0}\subset {\mathcal
  E}^{\sigma_2,R_2}_{h_2,h_1}$.  Then 

\medskip
\begin{compactenum}[(i)]
\item the function
\begin{equation}\label{eq:10}
u(x)=
\begin{cases}
  \omega_N^{-1/2},&\text{in } {\mathcal
    E}^{\sigma_1,R_1}_{h_1,0},\\[7pt]
  R_1^{1/\sigma_1}|x|^{-\frac1{\sigma_1}}\psi_1^{h_1}(x/|x|),&\text{in
  } {\mathcal E}^{\sigma_2,R_2}_{h_2,h_1}
  \setminus {\mathcal E}^{\sigma_1,R_1}_{h_1,0},\\[7pt]
  R_1^{1/\sigma_1}R_2^{1/\sigma_2}
  |x|^{-(\frac1{\sigma_1}+\frac1{\sigma_2})}\psi_1^{h_2}(x/|x|),&\text{in
  } \R^N\setminus{\mathcal E}^{\sigma_2,R_2}_{h_2,h_1},
\end{cases}
\end{equation}
belongs to $\Di$;
\item the distribution 
$$
{\mathcal H}=-\Delta
  u-\frac{h_1(x/|x|)}{|x|^2} \alchi_{{\mathcal
      E}^{\sigma_2,R_2}_{h_2,h_1} \setminus {\mathcal E}^{\sigma_1,R_1}_{h_1,0}}
\,u-\frac{h_2(x/|x|)}{|x|^2}
  \alchi_{\R^N\setminus{\mathcal E}^{\sigma_2,R_2}_{h_2,h_1}}\,u
$$
belongs to the dual space $\big(\Di\big)^\star$;
\item there exists a
positive constant $C$ (depending only on $N$, $\sigma_1$, $\sigma_2$,
$h_1$, and $h_2$) such that ${\mathcal H}\geq C\,\alchi_{\R^N
  \setminus {\mathcal E}^{\sigma_1,R_1}_{h_1,0}}\,
\frac{u}{|x|^2} $ in $\big(\Di\big)^\star$, i.e.
$$
{}_{(\Di)^\star}\langle {\mathcal H},w\rangle_{\Di}\geq 
C\,\int_{\R^N \setminus {\mathcal E}^{\sigma_1,R_1}_{h_1,0}}\frac{u\,w}{|x|^2}
$$
for any $w\in\Di$ such that $w\geq0$ a.e. in $\R^N$.
\end{compactenum}
\end{Lemma}
\begin{pf}
Let us consider the function
$$
u(x)=
\begin{cases}
\omega_N^{-1/2},&\text{in } {\mathcal
    E}^{\sigma_1,R_1}_{h_1,0},\\[6pt]
\varphi_1(x),&\text{in } {\mathcal
    E}^{\sigma_2,R_2}_{h_2,h_1}
  \setminus {\mathcal E}^{\sigma_1,R_1}_{h_1,0},\\[6pt]
 \varphi_2(x),&\text{in } \R^N\setminus{\mathcal
    E}^{\sigma_2,R_2}_{h_2,h_1},
\end{cases}
$$
where
\begin{gather*}
\varphi_1(x)=
 R_1^{1/\sigma_1}|x|^{-\frac1{\sigma_1}}\psi_1^{h_1}(x/|x|),\quad
\text{and}\quad\varphi_2(x)=R_1^{1/\sigma_1}R_2^{1/\sigma_2}
  |x|^{-(\frac1{\sigma_1}+\frac1{\sigma_2})}\psi_1^{h_2}(x/|x|).
\end{gather*}
By definition of the sets ${\mathcal
  E}^{\sigma_2,R_2}_{h_2,h_1}$
and ${\mathcal E}^{\sigma_1,R_1}_{h_1,0}$, it
follows that $u\in \Di\cap C^0(\R^N)$. Moreover
\begin{align*}
&  -\Delta\varphi_1-{\textstyle{
\frac{h_1(x/|x|)}{|x|^2}}}\,
\varphi_1
=\left[\big({\textstyle{\frac{N-2}2}}\big)^{2}+\mu_1(h_1)-
\big({\textstyle{\frac1{\sigma_1}-\frac{N-2}2}}\big)^{\!2}\right]
{\textstyle{\frac{\varphi_1}{|x|^2}}},\quad\text{in }\R^N\setminus\{0\},\\[10pt]
&  -\Delta\varphi_2-{\textstyle{
\frac{h_2(x/|x|)}{|x|^2}}}\,
\varphi_2
=\left[\big({\textstyle{\frac{N-2}2}}\big)^{2}+\mu_1(h_2)-
\big({\textstyle{\frac1{\sigma_1}+\frac1{\sigma_2}-\frac{N-2}2}}\big)^{\!2}
\right]
{\textstyle{\frac{\varphi_2}{|x|^2}}},\quad\text{in }\R^N\setminus\{0\}.
\end{align*}
Let us denote as $\nu_1(x)$ the outward normal derivate to ${\mathcal
  E}^{\sigma_1,R_1}_{h_1,0}$ and as $\nu_2(x)$ the
outward normal derivate to ${\mathcal
    E}^{\sigma_2,R_2}_{h_2,h_1}$. 
A direct calculation shows that
$$
\partial {\mathcal
  E}^{\sigma_1,R_1}_{h_1,0}=\{x\in\R^N:\ 
\varphi_1(x)=\omega_N^{-1/2}\}\quad\text{and}\quad
\nabla\varphi_1\cdot \frac{x}{|x|}=-
\frac{\varphi_1(x)}{\sigma_1|x|}<0\quad\text{on }\partial {\mathcal
  E}^{\sigma_1,R_1}_{h_1,0},
$$
hence
$\nu_1(x)=-\frac{\nabla\varphi_1(x)}
{|\nabla\varphi_1(x)|}$ for all $x\in \partial {\mathcal
  E}^{\sigma_1,R_1}_{h_1,0}$. In a similar way
$$
\partial {\mathcal
    E}^{\sigma_2,R_2}_{h_2,h_1}
=\{x\in\R^N:\ \varphi_1(x)=
\varphi_2(x)\}\quad\text{and}\quad
\nabla(\varphi_1-\varphi_2)\cdot \frac{x}{|x|}=
\frac{\varphi_2(x)}{\sigma_2|x|}>0\quad\text{on }\partial {\mathcal
    E}^{\sigma_2,R_2}_{h_2,h_1},
$$
hence
$\nu_2(x)=\frac{\nabla(\varphi_1-\varphi_2)(x)}
{|\nabla(\varphi_1-\varphi_2)(x)|}$ for all $x\in \partial {\mathcal
    E}^{\sigma_2,R_2}_{h_2,h_1}$. 
Therefore, for any $w\in\Di$, $w\geq 0$ a.e. in~$\R^N$,  
there holds
\begin{align*} &{\phantom{\frac12}}_{(\Di)^*}\bigg\langle-\Delta
  u-\frac{h_1(x/|x|)}{|x|^2} \alchi_{{\mathcal
      E}^{\sigma_2,R_2}_{h_2,h_1} \setminus {\mathcal
      E}^{\sigma_1,R_1}_{h_1,0}} \,u-\frac{h_2(x/|x|)}{|x|^2}
  \alchi_{\R^N\setminus{\mathcal E}^{\sigma_2,R_2}_{h_2,h_1}}\,u,w\bigg\rangle_{\Di}\\
  &\quad=-\int_{\partial{\mathcal E}^{\sigma_1,R_1}_{h_1,0}}
  \big(\nabla\varphi_1\cdot\nu_1\big)w\,dS+ \int_{\partial{\mathcal
      E}^{\sigma_2,R_2}_{h_2,h_1}}
  \big(\nabla(\varphi_1-\varphi_2)\cdot\nu_2\big)w\,dS\\
  &\quad\quad + \left[\bigg(\frac{N-2}2\bigg)^{\!\!2}+\mu_1(h_1)-
    \bigg(\frac1{\sigma_1}-\frac{N-2}2\bigg)^{\!\!2}\right] \int_{
    {\mathcal E}^{\sigma_2,R_2}_{h_2,h_1} \setminus {\mathcal
      E}^{\sigma_1,R_1}_{h_1,0}}{{\frac{\varphi_1\,w}{|x|^2}}}\\
  & \quad\quad+ \left[\bigg(\frac{N-2}2\bigg)^{\!\!2}+\mu_1(h_2)-
    \bigg(\frac1{\sigma_1}+\frac1{\sigma_2}-\frac{N-2}2\bigg)^{\!\!2}
  \right] \int_{\R^N\setminus{\mathcal E}^{\sigma_2,R_2}_{h_2,h_1}}
  {{\frac{\varphi_2\,w}{|x|^2}}}\geq C\, \int_{\R^N \setminus
    {\mathcal E}^{\sigma_1,R_1}_{h_1,0}} \frac{u\,w}{|x|^2},
\end{align*}
where $C=\min\big\{
\big({\textstyle{\frac{N-2}2}}\big)^{2}+\mu_1(h_1)-
\big({\textstyle{\frac1{\sigma_1}-\frac{N-2}2}}\big)^{\!2},
\big({\textstyle{\frac{N-2}2}}\big)^{2}+\mu_1(h_2)-
\big({\textstyle{\frac1{\sigma_1}+\frac1{\sigma_2}-\frac{N-2}2}}\big)^{\!2}\big\}>0$.
The proof is  complete.
\end{pf}

The Kelvin's transform yields the counterpart of Lemma \ref{l:intercapedine}
in the case of singularities located at  finite dipoles, as we prove below.

\begin{Lemma}\label{l:sha1}
Let $h_i\in L^{\infty}\big({\mathbb S}^{N-1}\big)$, 
$\mu_1(h_i)>-\big(\frac{N-2}2\big)^{\!2}$, $i=1,2$, 
$\sigma_1,\sigma_2>0$ satisfying \eqref{eq:7},
and $r_1,r_2\in(0,+\infty)$ such that ${\mathcal
  E}^{\sigma_2,r_2}_{h_1,h_2}\subset {\mathcal
  E}^{\sigma_1,r_1}_{0,h_1}$.  Then 

\medskip
\begin{compactenum}[(i)]
\item the function
\begin{equation}\label{eq:11}
v(x)=
\begin{cases}
  |x|^{-(N-2)}\omega_N^{-1/2},&\text{in } \R^N\setminus{\mathcal
    E}^{\sigma_1,r_1}_{0,h_1},\\[7pt]
  r_1^{-1/\sigma_1}|x|^{-(N-2)+\frac1{\sigma_1}}\psi_1^{h_1}(x/|x|),&\text{in
  } {\mathcal E}^{\sigma_1,r_1}_{0,h_1}
  \setminus {\mathcal E}^{\sigma_2,r_2}_{h_1,h_2},\\[7pt]
  r_1^{-1/\sigma_1}r_2^{-1/\sigma_2}
  |x|^{-(N-2)+\frac1{\sigma_1}+\frac1{\sigma_2}}\psi_1^{h_2}(x/|x|),&\text{in
  } {\mathcal E}^{\sigma_2,r_2}_{h_1,h_2},
\end{cases}
\end{equation}
belongs to $\Di$;
\item the distribution 
$$
{\mathcal H}=-\Delta
  v-\frac{h_1(x/|x|)}{|x|^2} \alchi_{{\mathcal
      E}^{\sigma_1,r_1}_{0,h_1} \setminus {\mathcal E}^{\sigma_2,r_2}_{h_1,h_2}}
\,v-\frac{h_2(x/|x|)}{|x|^2}
  \alchi_{{\mathcal E}^{\sigma_2,r_2}_{h_1,h_2}}\,v
$$
belongs to the dual space $\big(\Di\big)^\star$;
\item there exists a
positive constant $C$ (depending only on $N$, $\sigma_1$, $\sigma_2$,
$h_1$, and $h_2$) such that ${\mathcal H}\geq C\,\alchi_{{\mathcal E}^{\sigma_1,r_1}_{0,h_1}}\,
\frac{v}{|x|^2} $ in $\big(\Di\big)^\star$, i.e.
$$
{}_{(\Di)^\star}\langle {\mathcal H},w\rangle_{\Di}\geq 
C\,\int_{
{\mathcal E}^{\sigma_1,r_1}_{0,h_1}}\frac{v\,w}{|x|^2}
$$
for any $w\in\Di$ such that $w\geq0$ a.e. in $\R^N$.
\end{compactenum}
\end{Lemma}
\begin{pf}
  Let us consider the function $u$ defined in (\ref{eq:10}) with
  $R_1=1/r_1$ and $R_2=1/r_2$.  Then the function $v$ defined in
  (\ref{eq:11}) is the  Kelvin's transformed of $u$, i.e. 
$$
u(x)=|x|^{-(N-2)}v(x/|x|^2).
$$
Since  $\Delta v(x)=|x|^{-N-2}\Delta u(x/|x|^2)$, the conclusion follows 
from Lemma \ref{l:intercapedine}.
\end{pf}

\begin{Lemma}\label{l:sep_variante}
  Let $h_1,\dots,h_k,h_{\infty},H_1,\dots,H_k,H_{\infty}\in L^{\infty}\big({\mathbb
    S}^{N-1}\big)$
 satisfying 
$$
0<\max_{i=1,\dots,k,\infty}\big\{\Lambda_N(h_i),\Lambda_N(H_i)\big\}<1,
$$
$0<\e<  \Big(
\max\limits_{i=1,\dots,k,\infty}\big\{\Lambda_N(h_i),\Lambda_N(H_i)\big\}
\Big)^{-1}\!\!-1$,
    $\tilde h_i=(1+\e)h_i$, $\tilde H_i=(1+\e)H_i$ for all   $i=1,\dots,k,\infty$,
$\sigma_1,\sigma_2>0$ such that
$$
{\textstyle{\frac{N-2}2<\frac1{\sigma_1}<\frac1{\sigma_1}+\frac1{\sigma_2}<
\frac{N-2}2+\min_{i=1,\dots,k,\infty}
\left\{\frac{N-2}2,\sqrt{\frac{(N-2)^2}{4}+\mu_1(\tilde h_i)},
\sqrt{\frac{(N-2)^2}{4}+\mu_1(\tilde H_i)}\right\}}},
$$
$\{a_1,a_2,\dots,a_k\}\subset\R^N$, $a_i\neq a_j$ for $i\neq j$, and $R_0>0$ such that
$$
\{a_1,a_2,\dots,a_k\}\subset
B\Big (0,R_0\Big[\min_{\theta\in\mathbb
  S^{N-1}}\sqrt{\omega_N}\,\psi_1^{\tilde h_\infty}(\theta)
\Big]^{\sigma_1}\Big)\subset
{\mathcal E}^{\sigma_1,R_0}_{\tilde h_\infty,0}.
$$
Then there exists $\bar\delta>0$ such that for all $0<\delta<\bar\delta$,  for any $R$ such that
${\mathcal E}^{\sigma_1,R_0}_{\tilde h_\infty,0}\subset{\mathcal
  E}^{\sigma_2,R}_{\tilde H_\infty,\tilde h_\infty}$, and for any 
$r>0$ such that ${\mathcal E}^{\sigma_2,r}_{\tilde h_i,\tilde H_i}\subset
{\mathcal E}^{\sigma_1,\delta}_{0,\tilde h_i}$ for any $i=1,\dots,k$, there holds
\begin{align*}
  \mu\Bigg(&
  \sum_{i=1}^k\alchi_{
    {\mathcal E}^{\sigma_1,\delta}_{0,\tilde h_i}(a_i)\setminus
{\mathcal E}^{\sigma_2,r}_{\tilde h_i,\tilde H_i}(a_i)
  }(x)\frac{h_i\big(\frac{x-a_i}{|x-a_i|}\big)}{|x-a_i|^2}+
\sum_{i=1}^k\alchi_{
    {\mathcal E}^{\sigma_2,r}_{\tilde h_i,\tilde H_i}(a_i)
  }(x)\frac{H_i\big(\frac{x-a_i}{|x-a_i|}\big)}{|x-a_i|^2}\\
&\hskip3cm+
  \alchi_{{\mathcal
      E}^{\sigma_2,R}_{\tilde H_\infty,\tilde h_\infty} 
    \setminus {\mathcal E}^{\sigma_1,R_0}_{\tilde h_\infty,0}}
  \frac{h_\infty\big(\frac{x}{|x|}\big)}{|x|^2}+
  \alchi_{\R^N\setminus {\mathcal
      E}^{\sigma_2,R}_{\tilde H_\infty,\tilde h_\infty} }(x)
  \frac{H_\infty\big(\frac{x}{|x|}\big)}{|x|^2}
  \Bigg)
  \geq \frac{\e}{\e+1}.
\end{align*}
\end{Lemma}
\begin{pf}
  By scaling properties of the operator and in view of Lemma
  \ref{l:positivity_condition}, to prove the statement it is enough
  to find $\varphi\in\Di$ positive and continuous outside the
  singularities such that
\begin{equation}\label{eq:80bis}
-\Delta\varphi(x)-\sum_{i=1}^kV_i(x)\,\varphi(x)- V_{\infty}(x)\,\varphi(x)\geq0
\quad\text{in }\big(\Di\big)^\star,
\end{equation}
where 
$$
V_{\infty}(x)=\alchi_{{\mathcal E}^{\sigma_2,{R}/{\delta}}_{\tilde
    H_\infty,\tilde h_\infty} \setminus {\mathcal
    E}^{\sigma_1,{R_0}/\delta}_{\tilde h_\infty,0}}\frac{\tilde
  h_\infty(x/|x|)}{|x|^2} +\alchi_{\R^N\setminus{\mathcal E}^{\sigma_2,
    R/\delta}_{\tilde H_\infty,\tilde h_\infty}}\frac{\tilde
  H_\infty(x/|x|)}{|x|^2}
$$
and 
$$
V_i(x)=\alchi_{ {\mathcal E}^{\sigma_1,1}_{0,\tilde h_i}(a_i/\delta)
\setminus
{\mathcal E}^{\sigma_2,r/\delta}_{\tilde h_i,\tilde H_i}(a_i/\delta)
}
\frac{\tilde h_i\big(\frac{x-(a_i/\delta)}{|x-(a_i/\delta)|}\big)}
{|x-\frac{a_i}{\delta}|^2}
+\alchi_{
    {\mathcal E}^{\sigma_2,r/\delta}_{\tilde h_i,\tilde H_i}(a_i/\delta)
  }(x)\frac{\tilde H_i\big(\frac{x-(a_i/\delta)}{|x-(a_i/\delta)|}\big)}
{|x-\frac{a_i}{\delta}|^2},
$$ 
and $\d>0$  depends neither on $R$ nor on $r$ (but could depend on $\e$). 

Let us consider the function $\varphi_\infty(x)=u(\delta x)$, where 
$$
u(x)=
\begin{cases}
  \omega_N^{-1/2},&\text{in } {\mathcal
    E}^{\sigma_1,R_0}_{\tilde h_\infty,0},\\[10pt]
  R_0^{1/\sigma_1}|x|^{-\frac{1}{\sigma_1}}\psi_1^{\tilde h_\infty}(x/|x|),&\text{in }
  {\mathcal E}^{\sigma_2,R}_{\tilde H_\infty,\tilde h_\infty} \setminus
  {\mathcal E}^{\sigma_1,R_0}_{\tilde h_\infty,0},\\[10pt]
  R_0^{1/\sigma_1}R^{1/\sigma_2}|x|^{-(\frac1{\sigma_1}+\frac1{\sigma_2})}
\psi_1^{\tilde H_\infty}
(x/|x|),&\text{in } \R^N\setminus{\mathcal
    E}^{\sigma_2,R}_{\tilde H_\infty,\tilde h_\infty}.
\end{cases}
$$ 
 From Lemma
\ref{l:intercapedine}, we have that, for some positive constant $C$
(depending on $N$, $\sigma_1$, $\sigma_2$, $h_\infty$, $H_\infty$, $\e$,
but independent of $R$ and $\delta$),
$$
-\D\varphi_\infty(x)-V_\infty(x)\,\varphi_\infty(x)\geq C
\,\alchi_{\R^N\setminus {\mathcal E}^{\sigma_1,{R_0}/\delta}_{\tilde
    h_\infty,0}}\frac{\varphi_\infty(x)} {|x|^2}\quad\text{in
}\big(\Di\big)^\star.
$$
Let us also consider the functions $\varphi_i(x)=u_i\big(x-
\frac{a_i}{\delta}\big)$, $i=1,\dots,k$, where 
$$
u_i(x)=
\begin{cases}
  |x|^{-(N-2)}\omega_N^{-1/2},&\text{in } \R^N\setminus{\mathcal
    E}^{\sigma_1,1}_{0,\tilde h_i},\\[7pt]
  |x|^{-(N-2)+\frac1{\sigma_1}}\psi_1^{\tilde h_i}(x/|x|),&\text{in
  } {\mathcal E}^{\sigma_1,1}_{0,\tilde h_i}
  \setminus {\mathcal E}^{\sigma_2,r/\delta}_{\tilde h_i,\tilde H_i},\\[7pt]
  \big(\frac{r}{\delta}\big)^{-1/\sigma_2}
  |x|^{-(N-2)+\frac1{\sigma_1}+\frac1{\sigma_2}}\psi_1^{\tilde H_i}(x/|x|),&\text{in
  } {\mathcal E}^{\sigma_2,r/\delta}_{\tilde h_i,\tilde H_i}.
\end{cases}
$$
 From Lemma
\ref{l:sha1}, we have that, for some positive constant $\bar C$
(depending on $N$, $\sigma_1$, $\sigma_2$, $h_i$, $H_i$, and $\e$,
but independent of $\delta$ and $r$) and for all $i=1,\dots,k$,
$$
-\D\varphi_i-V_i\,\varphi_i\geq \bar C\, \alchi_{ {\mathcal
    E}^{\sigma_1,1}_{0,\tilde h_i}(a_i/\delta)} \,\frac{\varphi_i}
{\big|x-\frac{a_i}{\delta}\big|^2}\quad\text{in }\big(\Di\big)^\star.
$$
We notice that, if 
${\mathcal E}^{\sigma_2,r}_{\tilde h_i,\tilde H_i}\subset
{\mathcal E}^{\sigma_1,\delta}_{0,\tilde h_i}$ , then
$\frac{r}{\delta}\leq\min_{\mathbb S^{N-1}}\big(\psi_1^0\big)^{\sigma_1}
\big(\psi_1^{\tilde H_i}\big)^{\sigma_2}\big(\psi_1^{\tilde h_i}\big)^{-\sigma_2-\sigma_1}$.
Hence, by definition of $\varphi_i$, there exists some 
constant $C_0>0$ (depending on $N$, $\sigma_1$, $\sigma_2$, $h_i$, $H_i$, and $\e$,
but independent of~$\delta$ and $r$) such that, for all $i=1,\dots,k$, 
\begin{align}\label{eq:32bi}
  \varphi_i(x)\geq\frac1{C_0}\Big|x-\frac{a_i}\d\Big|^{-(N-2)+\frac1{\sigma_1}
    +\frac1{\sigma_2}},\quad&\text{ for all }x\in {\mathcal
    E}^{\sigma_1,1}_{0,\tilde h_i}(a_i/\delta).
\end{align}
Moreover, by the definition of $\varphi_\infty$, 
\begin{equation}\label{eq:6}
\varphi_\infty(x)\geq \frac1{C_1}|\delta x|^{-(\frac1{\sigma_1}+\frac1{\sigma_2})},\quad\text{ for
  all }x\in \R^N\setminus {\mathcal
    E}^{\sigma_1,R_0/\delta}_{\tilde h_\infty,0},
\end{equation}
where $C_1>0$ depends on $N$, $\sigma_1$, $\sigma_2$, 
$h_\infty$, $H_\infty$, $\e$, and $R_0$, but is independent of $R$ and $\delta$.
Let $\varphi=\sum_{i=1}^k\varphi_i+\eta\varphi_{\infty}$ for some
$\eta$
such that
$$
0<\eta<
 \min_{ i=1,\dots,k}\left\{ \frac{\bar C\sqrt{\omega_N}}{2C_0
\max\big\{\|\tilde
    h_i\|_{L^{\infty}({\mathbb S}^{N-1})},\|\tilde
    H_i\|_{L^{\infty}({\mathbb S}^{N-1})}\big\}
} \bigg[\max_{\theta\in\mathbb
   S^{N-1}}\frac{\omega_N^{-1/2}} {\psi_1^{\tilde
     h_i}(\theta)}\bigg]^{1+\frac{\sigma_1}{\sigma_2}-\sigma_1(N-2)}\right\}.
$$
 Then we have 
\begin{align*}
-\D\varphi(x)-\sum_{i=1}^kV_i(x)\,\varphi(x)&-
 V_{\infty}(x)\,\varphi(x)\geq f(x)
\end{align*}
in $\big(\Di\big)^\star$, where
\begin{align*}
f(x)=&\bar C\,\sum_{i=1}^k  \alchi_{ {\mathcal
    E}^{\sigma_1,1}_{0,\tilde h_i}(a_i/\delta)} \,\frac{\varphi_i}
{\big|x-\frac{a_i}{\delta}\big|^2}+
C\,\eta
\,\alchi_{\R^N\setminus {\mathcal E}^{\sigma_1,{R_0}/\delta}_{\tilde
    h_\infty,0}}\frac{\varphi_\infty(x)} {|x|^2}-\sum_{\substack{i,j=1\\i\neq
    j}}^kV_i(x)\varphi_j(x)\\
&-\eta\sum_{i=1}^kV_i(x)\varphi_{\infty}(x)-V_{\infty}(x)
\sum_{i=1}^k\varphi_i(x).
\end{align*}
In particular a.e. in the set ${\mathcal
  E}^{\sigma_1,R_0/\delta}_{\tilde h_\infty,0}\setminus\bigcup_{i=1}^k
{\mathcal E}^{\sigma_1,1}_{0,\tilde h_i}(a_i/\delta)$, we have that
$f(x)\equiv 0$.  Let us consider ${\mathcal E}^{\sigma_1,1}_{0,\tilde
  h_i}(a_i/\delta)$. Since, for $\delta$ small,
$${\mathcal E}^{\sigma_1,1}_{0,\tilde
  h_i}(a_i/\delta)\subset  {\mathcal
    E}^{\sigma_1,R_0/\delta}_{\tilde h_\infty,0}\quad\text{and}\quad
{\mathcal E}^{\sigma_1,1}_{0,\tilde
  h_i}(a_i/\delta)\subset\R^N\setminus {\mathcal E}^{\sigma_1,1}_{0,\tilde
  h_j}(a_j/\delta)\text{ for }j\neq i,
$$
from
(\ref{eq:32bi})  it
follows that,
in  ${\mathcal E}^{\sigma_1,1}_{0,\tilde
  h_i}(a_i/\delta)$, 
\begin{align*}
  f(x)&\geq \bar C\, \alchi_{ {\mathcal E}^{\sigma_1,1}_{0,\tilde
      h_i}(a_i/\delta)} \,\frac{\varphi_i}
  {\big|x-\frac{a_i}{\delta}\big|^2}
  -V_i(x)\Big(\sum_{\substack{i,j=1\\i\neq j}}^k
  \varphi_j(x)+\eta\,\varphi_{\infty}(x)\Big)\\
  &\geq\Big|x-\frac{a_i}{\delta}\Big|^{-2} \Bigg[\frac{\bar
    C}{C_0}\Big|x-\frac{a_i}{\delta}\Big|^{-(N-2)+\frac1{\sigma_1}+\frac1{\sigma_2}}\\
&  \quad-\omega_N^{-1/2}
\max\Big\{\|\tilde
    h_i\|_{L^{\infty}({\mathbb S}^{N-1})},\|\tilde
    H_i\|_{L^{\infty}({\mathbb S}^{N-1})}\Big\}
  \Bigg(\sum_{\substack{j=1\\j\neq i}}^k
  \Big|x-\frac{a_j}{\delta}\Big|^{-(N-2)}+\eta 
  \Bigg)\Bigg].
\end{align*}
It is easy to see that, for $\delta$ small, in  ${\mathcal E}^{\sigma_1,1}_{0,\tilde
  h_i}(a_i/\delta)$ 
\begin{align*}
&\Big|x-\frac{a_i}{\delta}\Big|^{-(N-2)+\frac1{\sigma_1}+\frac1{\sigma_2}}\geq
\bigg[\max_{\theta\in\mathbb
  S^{N-1}}\frac{\omega_N^{-1/2}}
{\psi_1^{\tilde h_i}(\theta)}\bigg]^{1+\frac{\sigma_1}{\sigma_2}-\sigma_1(N-2)}\quad\text{and}\\
&\Big|x-\frac{a_j}{\delta}\Big|^{-(N-2)}\leq\Big(\frac{2}{|a_i-a_j|}\Big)^{N-2}
\d^{N-2}<\frac{\eta}{k-1}, 
\end{align*}
and hence the choice of $\eta$ ensures that $f\geq 0$ a.e. in
${\mathcal E}^{\sigma_1,1}_{0,\tilde h_i}(a_i/\delta)$, provided
$\delta$ is sufficiently small.
Let us finally consider $\R^N\setminus {\mathcal
    E}^{\sigma_1,R_0/\delta}_{\tilde h_\infty,0}$. 
From  (\ref{eq:6}), we deduce that in $\R^N\setminus {\mathcal
    E}^{\sigma_1,R_0/\delta}_{\tilde h_\infty,0}$
\begin{align*}
  f(x)&\geq \frac1{|x|^{2}}\bigg[ \frac{C\,\eta}{C_1}\big|\d
  x\big|^{-(\frac1{\sigma_1}+\frac1{\sigma_2})}-\omega_N^{-1/2}\max\{\|\tilde
    h_\infty\|_{L^{\infty}({\mathbb S}^{N-1})},\|\tilde H_\infty
    \|_{L^{\infty}({\mathbb S}^{N-1})}\}
\sum_{i=1}^k \Big|x-\frac{a_i}{\delta}\Big|^{-(N-2)}\bigg]. 
\end{align*}
It is easy to see that,  in $\R^N\setminus {\mathcal
    E}^{\sigma_1,R_0/\delta}_{\tilde h_\infty,0}$, 
$\big|x-\frac{a_i}{\delta}\big|\geq\big(1-\frac{\a}{\beta R_0}\big)|x|$
where $\a=\max_{j}\{|a_j|\}$ and $\beta=\big[\min_{\theta\in\mathbb
  S^{N-1}}\sqrt{\omega_N}\,\psi_1^{\tilde h_\infty}(\theta)
\big]^{\sigma_1}$, hence 
\begin{multline*}
f(x)\geq |x|^{-2-\frac{1}{\sigma_1}-\frac{1}{\sigma_2}} \bigg[
\frac{C}{C_1}\,\eta\,\delta^{-(\frac1{\sigma_1}+\frac1{\sigma_2})}\\
-\omega_N^{-1/2}
\max\{\|\tilde
    h_\infty\|_{L^{\infty}({\mathbb S}^{N-1})},\|\tilde
    H_\infty\|_{L^{\infty}({\mathbb S}^{N-1})}\}
k\bigg(1-\frac{\alpha}{\beta R_0}\bigg)^{-(N-2)}(R_0\beta)^{-(N-2)+
\frac1{\sigma_1}+\frac1{\sigma_2}}
\bigg]>0
\end{multline*}
a.e. in $\R^N\setminus {\mathcal
  E}^{\sigma_1,R_0/\delta}_{\tilde h,0}$, provided
$\delta$ is sufficiently small (notice that the choice of $\delta$ is
independent of~$R$ and $r$).  The proof is thereby complete.
\end{pf}

\begin{Lemma}\label{l:sep_parte2}
Let  $a_1,a_2,\dots,a_k\in\R^N$, $a_i\neq
a_j$ for $i\neq j$, and $h_i
\in L^{\infty}\big({\mathbb S}^{N-1}\big)$, $i=1,\dots,k,\infty$, with
$\mu_1(h_i)>-(N-2)^2/4$, for $i=1,\dots,k,\infty$, and 
$\max_{i=1,\dots,k,\infty}\Lambda_N(h_i)>0$. Then for every
$0<\alpha<1-\max_{i=1,\dots,k,\infty}\Lambda_N(h_i)$  there exist
${\mathcal U}_1,\dots,{\mathcal U}_k,{\mathcal U}_\infty\subset\R^N$
such that ${\mathcal U}_i$ is a neighborhood of $a_i$ for every 
$i=1,\dots,k$, ${\mathcal U}_\infty$ is a neighborhood of $\infty$, and 
\begin{align*}
\mu\bigg(
\sum_{i=1}^k\alchi_{{\mathcal U}_i}(x)
\frac{h_i\big(\frac{x-a_i}{|x-a_i|}\big)}{|x-a_i|^2}
+\alchi_{{\mathcal U}_\infty}(x)
\frac{h_\infty\big(\frac{x}{|x|}\big)}{|x|^2}
\bigg)\geq
1-\max_{i=1,\dots,k,\infty}\Lambda_N(h_i)-\alpha>0.
\end{align*}
\end{Lemma}
\begin{pf}
For any $0<\alpha<1-\max_{i=1,\dots,k,\infty}\Lambda_N(h_i)$, the statement follows 
from Lemma \ref{l:sep_variante} with $\e=\big(\alpha+
\max_{i=1,\dots,k,\infty}\Lambda_N(h_i)\big)^{-1}-1$ and $H_i=h_i$, $i=1,\dots,k,\infty$.
\end{pf}

\noindent 
\begin{pfn}{Lemma \ref{l:sep}}
It follows from Lemmas \ref{l:sep_parte1} and \ref{l:sep_parte2}.
\end{pfn}

\noindent 
\begin{pfn}{Theorem \ref{l:semibounded}}
It follows from Lemmas \ref{l:sep} and Proposition \ref{t:mainresult}.
\end{pfn}

\noindent Lemma \ref{l:sep} can be extended to infinitely many
dipole-type singularities distributed on
reticular structures.

\begin{Lemma}{\bf [\,Shattering of reticular singularities\,]}\label{l:ret}
Let 
 $\{h_n\}_{n\in\N}\subset L^{\infty}\big({\mathbb
    S}^{N-1}\big)$ such that 
\begin{equation}\label{eq:25}
\sup_{n\in\N}\|h_n\|_{L^{\infty}({\mathbb S}^{N-1})}<+\infty
\quad\text{and}\quad
0<\sup_{n\in\N}\Lambda_N(h_n)<1,
\end{equation}
and $\{a_n\}_n\subset\R^N$ satisfy 
\begin{equation}\label{eq:55}
\sum_{n=1}^{\infty}|a_n|^{-(N-2)}<+\infty,\quad 
\sum_{{\substack{m=1\\m\neq n}}}^{\infty}|a_{m}-a_n|^{-(N-2)}\quad\text{is bounded uniformly
  in }n,
\end{equation}
and $|a_n-a_m|\geq 1$
for all $n\neq m$. Then  there exist $\e>0$, $\sigma>0$, and $\bar \delta>0$,  
such that, for all $0<\delta<\bar \delta$,
$$
\inf_{{\substack{u\in \Di
\\u\not\equiv0}}}\frac{\displaystyle\int_{\R^N}\big(|\n
u(x)|^2-V(x)u^2(x)\big)\,dx}{\displaystyle\int_{\R^N}|\n
u(x)|^2\,dx}>0
$$
where
$$
V(x)=\sum_{n=1}^{\infty}\frac{h_n\big(\frac{x-a_n}{|x-a_n|}\big)}{|x-a_n|^2}
\alchi_{{\mathcal E}^{\sigma,\delta}_{0, (1+\e)h_n}(a_n)}.
$$
\end{Lemma}
\begin{pf}
  Taking into account the characterization of $\Lambda_N(h)$ given in
  \cite[Lemma 2.4]{FMT2}, a direct calculation yields that, for any
  $h\in L^{\infty}\big({\mathbb S}^{N-1}\big)$,
\begin{align}
\label{eq:24}&\Lambda_N(h)\geq -\frac4{(N-2)^2}\mu_1(h).
\end{align}
Letting $0<\e<\big(\sup_{n\in\N}\Lambda_N(h_n)\big)^{-1}-1$
and $\tilde h_n:=(1+\e)h_n$, from (\ref{eq:24}), it follows that
$$
\inf_{n\in\N}\mu_1(\tilde h_n)>-\bigg(\frac{N-2}2\bigg)^{\!\!2},
$$
hence there exists  
$\sigma>0$ such that
\begin{equation}\label{eq:12}
{{\frac{N-2}2<\frac1{\sigma}<
\frac{N-2}2+\min\left\{\frac{N-2}2,\sqrt{\frac{(N-2)^2}{4}+\inf_{n\in\N}
\mu_1(\tilde h_n)}\right\}}}.
\end{equation}
Moreover, classical elliptic regularity theory and bootstrap methods (see also 
\cite[Lemma 2.3]{FMT2}) yield that  $\psi_1^{\tilde h_n}$ is
bounded in $C^{0,\alpha}\big({\mathbb
    S}^{N-1}\big)$ uniformly with respect to $n$, thus implying that 
the sequence $\{\psi_1^{\tilde h_n}\}_{n\in\N}$ is equi-continuous and, 
by Ascoli-Arzel\`a's Theorem,  
compact in  $C^{0}\big({\mathbb
    S}^{N-1}\big)$. We deduce that, for positive constant $C>0$ independent on $n$,
\begin{equation}\label{eq:22}
\frac1C\leq \psi_1^{\tilde h_n}(\theta)\leq C\quad\text{for all }\theta\in 
{\mathbb S}^{N-1}.
\end{equation}
Let 
$$
\psi_n(x)=
\begin{cases}
|x|^{-(N-2)+\frac1\sigma}\,\psi_1^{\tilde h_n}(x/|x|),&\text{in } {\mathcal
    E}^{\sigma,1}_{0,\tilde h_n},\\[7pt]
|x|^{-(N-2)}\omega_N^{-1/2},&\text{in } \R^N\setminus{\mathcal
    E}^{\sigma,1}_{0,\tilde h_n}.
\end{cases}
$$
Arguing as in Lemma \ref{l:sha1}, one can easily verify that $\psi_n\in\Di$ and
$$
-\Delta
  \psi_n-\frac{\tilde h_n(x/|x|)}{|x|^2} \alchi_{{\mathcal
      E}^{\sigma,1}_{0,\tilde h_n}}\,\psi_n\geq \bar C\frac{\psi_n}{|x|^2}\, 
\alchi_{{\mathcal
      E}^{\sigma,1}_{0,\tilde h_n}}\quad \text{in }\big(\Di\big)^\star,
$$
for some $\bar C>0$ depending on $N$, $\inf_{n\in\N}
\mu_1(\tilde h_n)$,  and $\sigma$.

Let
$\varphi(x)=\sum_{n=1}^{\infty}\psi_n\big(x-\frac{a_n}\delta\big)$.
Due to estimate (\ref{eq:22}), for any open set $\Omega\subset\R^N$
such $\overline\Omega$ is compact, we have that there exists $\bar n$
such that $\big|\psi_n\big(x-\frac{a_n}\delta\big)\big|
=\omega_N^{-1/2} \big|x-\frac{a_n}\delta\big|^{-(N-2)}\leq{\rm
  const\,}\big|\frac{a_n}\delta\big|^{-(N-2)}$ 
and $\big|\nabla \psi_n\big(x-\frac{a_n}\delta\big)\big|
=\omega_N^{-1/2}(N-2) \big|x-\frac{a_n}\delta\big|^{-(N-1)}\leq{\rm
  const\,}\big|\frac{a_n}\delta\big|^{-(N-1)}$ 
for all $n\geq \bar n$
and $x\in \Omega$.  Then
$$
\varphi\big|_\Omega(x)=\sum_{n=1}^{\bar
  n-1}\psi_n\Big(x-\frac{a_n}\delta\Big)+\sum_{n=\bar
  n}^{\infty}\psi_n\Big(x-\frac{a_n}\delta\Big)\in H^1(\Omega)+L^{\infty}(\Omega)
$$
and $\nabla\big(\varphi\big|_\Omega)\in L^2(\Omega)$.  In particular
$\varphi\in L^{2^*}_{\rm loc}(\R^N)$, $\varphi\in H^1_{\rm
  loc}(\R^N)$, and $\sum_{n=1}^{\infty}\psi_n\big(x-\frac{a_n}\delta\big)$ converges 
to $\varphi$ in $H^1_0(\Omega)$ for every $\Omega\Subset\R^N$. Moreover
$$
-\Delta\varphi(x)-\sum_{n=1}^{\infty}\frac{\tilde h_n\big(\frac{x-\frac{a_n}{\delta}}
{|x-\frac{a_n}{\delta}|}\big)}{|x-\frac{a_n}{\delta}|^2}
\alchi_{{\mathcal E}^{\sigma,1}_{0,\tilde h_n}(a_n/\delta)}
\varphi(x)\geq f(x) \quad \text{in }\big(H^1_0(\Omega))^\star\quad\text{for all }
\Omega\Subset\R^N,
$$
where
$$
f(x)=
\bar C\sum_{n=1}^{\infty}\frac{
\psi_n\big(x-\frac{a_n}\delta\big)}{\big|x-\frac{a_n}\delta \big|^2}
\alchi_{{\mathcal
      E}^{\sigma,1}_{0,\tilde h_n}(a_n/\delta)}(x)
- \sum_{\substack{n,m=1\\ m\neq
  n}}^\infty\frac{\tilde h_n\big(\frac{x-({a_n}/{\delta})}
{|x-{a_n}/{\delta}|}\big)}{|x-\frac{a_n}{\delta}|^2}
\alchi_{{\mathcal E}^{\sigma,1}_{0,\tilde h_n}(a_n/\delta)}
\psi_m \Big(x-\frac{a_m}\delta\Big).
$$
A.e. in the set $\R^N\setminus\bigcup_{n=1}^{\infty} {\mathcal
      E}^{\sigma,1}_{0,\tilde h_n}(a_n/\delta)$, we have that $f(x)= 0$.
From the definition of $\psi_n$ and estimate (\ref{eq:22}), it follows that in each  ${\mathcal
      E}^{\sigma,1}_{0,\tilde h_n}(a_n/\delta)$
\begin{align*}
f(x)\geq C_1\,\Big|x-\frac{a_n}\delta\Big|^{-2}\bigg(
C_2\Big|x-\frac{a_n}\delta\Big|^{-(N-2)+1/\sigma}-\sum_{m\neq
  n}\Big|x-\frac{a_m}\delta\Big|^{-(N-2)}\bigg),
\end{align*}
for some constants $C_1,C_2$ independent of $n$.
Since, for small $\d$, $\big|x-\frac{a_m}\delta\big|\geq
\frac{|a_m-a_n|}\d-{\rm const}\geq \frac{|a_m-a_n|}{2\d}$ provided $\d$ small
enough, we deduce that
$$
\sum_{m\neq
  n}\Big|x-\frac{a_m}\delta\Big|^{-(N-2)}\leq (2\d)^{N-2}\sum_{m\neq
  n}|a_m-a_n|^{-(N-2)}\leq {\rm const\,}\d^{N-2}.
$$
Hence, we can choose $\delta$ small enough independently of $n$ such
that $f(x)\geq 0$ a.e. in ${\mathcal E}^{\sigma,1}_{0,\tilde
  h_n}(a_n/\delta)$.
Hence we have constructed  a supersolution  
$\varphi\in L^1_{\rm loc}(\R^N)$, $\varphi>0$ in
$\R^N\setminus\{a_n/\delta\}_{n\in\N}$ and $\varphi$ continuous in
$\R^N\setminus\{a_n/\delta\}_{n\in\N}$, such that
$$
-\Delta\varphi-\sum_{n=1}^{\infty}\frac{h_n\big(\frac{x-\frac{a_n}{\delta}}
{|x-\frac{a_n}{\delta}|}\big)}{|x-\frac{a_n}{\delta}|^2}
\alchi_{{\mathcal E}^{\sigma,1}_{0,\tilde h_n}(a_n/\delta)}
\varphi\geq \e \sum_{n=1}^{\infty}\frac{h_n\big(\frac{x-\frac{a_n}{\delta}}
{|x-\frac{a_n}{\delta}|}\big)}{|x-\frac{a_n}{\delta}|^2}
\alchi_{{\mathcal E}^{\sigma,1}_{0,\tilde h_n}(a_n/\delta)}
\varphi 
 \quad \text{in }\big(H^1_0(\Omega))^\star,
$$
for all $\Omega\Subset\R^N$.
Therefore, arguing as in Lemma \ref{l:positivity_condition} and taking into
account the scaling properties of the operator, we  obtain 
$$
\tilde\mu(V):=\inf_{{\substack{u\in C^{\infty}_{\rm
        c}\left(\R^N\setminus\{a_n\}_{n\in\N}\right) 
\\u\not\equiv0}}}\frac{\displaystyle\int_{\R^N}\big(|\n
u(x)|^2-V(x)u^2(x)\big)\,dx}{\displaystyle\int_{\R^N}|\n
u(x)|^2\,dx}>0,
$$
i.e.
\begin{equation}\label{eq:56}
\int_{\R^N}V(x)u^2(x)\,dx\leq(1-\tilde\mu(V))\int_{\R^N}|\n
u(x)|^2\,dx
\end{equation}
for all $u\in C^{\infty}_{\rm
        c}\left(\R^N\setminus\{a_n\}_{n\in\N}\right)$.
By density of $C^{\infty}_{\rm
        c}\left(\R^N\setminus\{a_n\}_{n\in\N}\right)$ in $\Di$ and the
      Fatou Lemma, we can easily prove that (\ref{eq:56}) holds for
      all $u\in\Di$.

\end{pf}

\begin{remark}\label{rem:ipo}
Taking into account the characterization of $\Lambda_N(h)$ given in 
\cite[Lemma2.4]{FMT2}, we can easily prove that, for any $h\in L^{\infty}\big({\mathbb
    S}^{N-1}\big)$,
\begin{align*}
 \mu_1(h)\leq
  -\bigg(\frac{N-2}2\bigg)^{\!\!2}+\left[(\Lambda_N(h))^{-1}-1\right]\|h^+\|_{L^{\infty}({\mathbb
      S}^{N-1})},
\end{align*}
which, together with (\ref{eq:24}), implies that, for $\{h_n\}_{n\in\N}$ bounded in 
$L^{\infty}({\mathbb
      S}^{N-1})$,  
$$\sup_{n\in\N}\Lambda_N(h_n)<1\quad\text{if and only if}
\quad \inf_{n\in\N}\mu_1(h_n)>-\bigg(\frac{N-2}2\bigg)^{\!\!2}.
$$
\end{remark}

\section{Spectral stability under perturbation at singularities}\label{sec:pertinf}

In this section we discuss the stability of positivity with respect to
perturbations of the potentials with a singularity sitting at 
dipolar-shaped neighborhoods either of a dipole or of infinity.  

In order to analyze the stability of the sign of $\mu(V)$ under
perturbations at singularities, it is useful to investigate its
attainability.  Due to inverse square homogeneity, it is easy to
verify that $\mu\big(\frac{h(x/|x|)}{|x|^2}\big)=1-\Lambda_N(h)$ is
not attained if $h\in L^{\infty}\big({\mathbb S}^{N-1}\big)$ is
positive somewhere, i.e. if $\Lambda_N(h)>0$. 
If $h\leq0$ a.e. in $\R^N$, then $\mu\big(\frac{h(x/|x|)}{|x|^2}\big)=1$ 
could be  achieved, as it happens for example  if $h$ vanishes in a nonempty open set.

The best constant in the Hardy-type inequality associated
to a multisingular potential $V\in{\mathcal V}$ is attained 
  if  $\mu(V)$ 
stays strictly below the bound provided in Lemma~\ref{l:estmu}.

\begin{Proposition}\label{p:attai}
Let $V\in{\mathcal V}$ be as in (\ref{eq:viinV}). If 
\begin{equation}\label{eq:77}
\mu(V)<1-\max\big\{
0,\Lambda_N(h_1),\dots,\Lambda_N(h_k),\Lambda_N(h_\infty)\big\}
\end{equation}
then $\mu(V)$ is attained.
\end{Proposition}
\begin{pf}
Let us assume that (\ref{eq:77}) holds. 
Hence there exists  $\alpha>0$ such that $\mu(V)=
1-\bar\Lambda-\alpha$ where $\bar \Lambda=
\max\big\{
0,\Lambda_N(h_1),\dots,\Lambda_N(h_k),\Lambda_N(h_\infty)\big\}$.
From Lemmas \ref{l:sep_parte1} and \ref{l:sep_parte2}, there exist
$\widetilde V\in{\mathcal V}$ and $\widetilde W\in L^{N/2}(\R^N)\cap
L^{\infty}(\R^N)$ such that $V= \widetilde V+\widetilde W$ and 
\begin{equation}\label{eq:76}
\mu(\widetilde V)\geq 1-\bar\Lambda-\frac{\alpha}2.
\end{equation}
 Let $u_n\in\Di$ a minimizing sequence
for $\mu(V)$, namely  
$$
\int_{\R^N}|\n u_n(x)|^2\,dx=1\quad\text{and}\quad 
\int_{\R^N}\big(|\n u_n(x)|^2-V(x)u_n^2(x)\big)\,dx=\mu(V)+o(1)
\quad\text{as }n\to\infty.
$$
Being $\{u_n\}_n$ bounded in $\Di$, we can assume that, up to a
subsequence still denoted as $u_n$, $u_n$ converges to
some $u$ a.e. and weakly in $\Di$.   Since
\begin{align*}
1-\bar\Lambda-\frac{\alpha}2\leq
\mu(\widetilde V)&\leq \int_{\R^N}\big(|\n
u_n(x)|^2-V(x)u_n^2(x)\big)\,dx+\int_{\R^N}\widetilde
W(x)u_n^2(x)\,dx\\
&
=\mu(V)+\int_{\R^N}\widetilde W(x)u^2(x)\,dx+o(1)
=1-\bar\Lambda-\alpha+\int_{\R^N}\widetilde W(x)u^2(x)\,dx
\end{align*}
as $n\to\infty$, we obtain that $\int_{\R^N}\widetilde
W(x)u^2(x)\,dx\geq\frac{\alpha}2>0$, thus implying $u\not\equiv 0$.
From weak convergence of $u_n$ to $u$, we deduce that
\begin{align*}
&\frac{\int_{\R^N}\big(|\n u(x)|^2-V(x)u^2(x)\big)\,dx}{\int_{\R^N}|\n
  u(x)|^2\,dx}\\[5pt]
&\ =\frac{\big[\int_{\R^N}\big(|\n
  u_n(x)|^2-V(x)u_n^2(x)\big)\,dx\big]-\big[ \int_{\R^N}\big(|\n
  (u_n-u)(x)|^2-V(x)(u_n-u)^2(x)\big)\,dx   \big]+o(1)}{\int_{\R^N}|\n
  u_n(x)|^2\,dx-\int_{\R^N}|\n
  (u_n-u)(x)|^2\,dx+o(1)}\\[5pt]
&\ \leq \mu(V)\,\frac{1-\int_{\R^N}|\n
  (u_n-u)(x)|^2\,dx+o(1)}{1-\int_{\R^N}|\n (u_n-u)(x)|^2\,dx+o(1)}=
\mu(V)\bigg[1+\frac{o(1)}{\int_{\R^N}|\n u(x)|^2\,dx+o(1)}\bigg]
\quad\text{as }n\to\infty.
\end{align*}
Letting $n\to\infty$, we  obtain that $u$ attains the infimum defining
$\mu(V)$.
\end{pf}

Let us now study the stability of the sign of $\mu(V)$ under perturbations at infinity.

\begin{Theorem}\label{t:pertinf}
  For $i=1,\dots,k$, let $r_i,R\in\R^+$, $a_i\in\R^N$, $a_i\neq a_j$
  for $i\neq j$, $h_i,h_\infty\in L^{\infty}\big({\mathbb S}^{N-1}\big)$ with
$\mu_1(h_i)>-(N-2)^2/4$, $\mu_1(h_\infty)>-(N-2)^2/4$,  and 
$$
V(x)=\sum_{i=1}^k\alchi_{B(a_i,r_i)}(x)\frac{
h_i\big(\frac{x-a_i}{|x-a_i|}\big)}{|x-a_i|^2}
      +\alchi_{\R^N\setminus B(0,R)}(x)
\frac{
h_\infty\big(\frac{x}{|x|}\big)}{|x|^2}+W(x)\in{\mathcal V}
$$
where $W\in L^{N/2}(\R^N)\cap L^{\infty}(\R^N)$. Assume that
$\mu(V)>0$ and let $h\in
L^{\infty}\big({\mathbb S}^{N-1}\big)$ such that 
$\mu_1(h+h_\infty)>-(N-2)^2/4$,
 $\e>0$ such that, setting $H:=h+h_\infty$,
$$
\e<
\begin{cases}
\frac1{\max_{i=1,\dots,k,\infty}\{
\Lambda_N(h_i),
\Lambda_N(H)\}}-1,&\text{if }
\max_{i=1,\dots,k,\infty}\big\{
\Lambda_N(h_i),
\Lambda_N(H)\big\}>0,\\
+\infty&\text{if }
\max_{i=1,\dots,k,\infty}\big\{
\Lambda_N(h_i),
\Lambda_N(H)\big\}= 0,
\end{cases}
$$
 and  $\sigma>0$ such that 
$
0<\frac1\sigma<\min_{i=1,\dots,k,\infty}
\Big\{\sqrt{(N-2)^2/4+\mu_1(\tilde h_i)},
\sqrt{(N-2)^2/4+\mu_1(\tilde H)}\Big\},
$
where $\tilde h_i:=(1+\e)h_i$ and $\tilde H:=(1+\e)H$.
  Then there exists
 $\bar R$ such that
$$
\mu\bigg(V+\frac{h(x/|x|)}{|x|^2}\alchi_{\R^N\setminus {\mathcal E}^{\sigma, \tilde
    R}_{\tilde H,\tilde h_\infty}}\bigg)>0,
$$
for all $\tilde R>\bar R$, where ${\mathcal E}^{\sigma, \tilde
    R}_{\tilde H,\tilde h_\infty}$ is defined in (\ref{eq:26}).
\end{Theorem}
\begin{pf}
Let $\sigma_1>0$ such that 
$$
\frac{N-2}2<\frac1{\sigma_1}<\frac1{\sigma_1}+\frac1{\sigma}<
\frac{N-2}2+
\min_{i=1,\dots,k,\infty}
\bigg\{\sqrt{\frac{(N-2)^2}{4}+\mu_1(\tilde h_i)},
\sqrt{\frac{(N-2)^2}{4}+\mu_1(\tilde H)}\bigg\}.
$$
Assume, by contradiction, that there exists a sequence $R_n\to+\infty$
such that, setting 
$$
V_n= 
V+\frac{h(x/|x|)}{|x|^2}\alchi_{\R^N\setminus {\mathcal E}^{\sigma, 
R_n}_{\tilde H,\tilde h_\infty}},
$$ 
there holds $\mu(V_n)\leq 0$.  By Lemmas \ref{l:sep_parte1} and \ref{l:sep_variante},
$V_n=\widetilde V_n+\widetilde W$, where $ \mu(\widetilde V_n)\geq \frac\e{\e+1}>0$
and
$$
\widetilde W(x)=
\sum_{i=1}^k\alchi_{B(a_i,r_i)\setminus 
 {\mathcal E}^{\sigma_1,\delta}_{0,\tilde h_i}(a_i)}(x)
\frac{h_i\big(\frac{x-a_i}{|x-a_i|}\big)}{|x-a_i|^2}
+ \alchi_{
{\mathcal E}^{\sigma_1,{R_0}}_{\tilde h_\infty,0}
\setminus B(0,R)}(x)\frac{
h_\infty\big(\frac{x}{|x|}\big)}{|x|^2}+W(x)
$$
for some  $\d>0$ and $R_0>0$ independent of $n$.
By Proposition \ref{p:attai}, $\mu(V_n)$ is attained by some $\varphi_n\in\Di$ 
such that 
$$
\int_{\R^N}|\nabla \varphi_n(x)|^2\,dx=1,\quad 
\int_{\R^N}|\nabla \varphi_n(x)|^2\,dx-\int_{\R^N}V_n(x)\varphi_n^2(x)\,dx=\mu(V_n),
$$
and
\begin{equation}\label{eq:2}
-\Delta\varphi_n-V_n\varphi_n=-\mu(V_n)\Delta\varphi_n\quad\text{in }\R^N.
\end{equation}
Up to a subsequence, $\varphi_n\weakly \varphi$ weakly in $\Di$
for some $\varphi\in\Di$, and hence 
\begin{align*}
\frac\e{\e+1}\leq\mu(\widetilde V_n)&\leq\int_{\R^N}|\nabla
\varphi_n(x)|^2\,dx-\int_{\R^N}V_n(x)\varphi_n^2(x)\,dx+\int_{\R^N}\widetilde
W(x)\varphi_n^2(x)\,dx\\
&=\mu(V_n)+\int_{\R^N}\widetilde
W(x)\varphi_n^2(x)\,dx\leq \int_{\R^N}\widetilde
W(x)\varphi^2(x)\,dx+o(1)\quad\text{as }n\to+\infty.
\end{align*}
Therefore $\int_{\R^N}\widetilde
W(x)\varphi^2(x)\,dx>0$ and we conclude that $\varphi\not\equiv 0$.
We claim that 
\begin{equation}\label{eq:82bi}
\lim_{j\to\infty}\int_{\R^N} V_{n}(x)\varphi_n(x)\varphi(x)\,dx=
\int_{\R^N}V(x)\varphi^2(x)\,dx.
\end{equation}
Indeed for any $\eta>0$, by density there exists $\psi\in
C^{\infty}_{\rm c}(\R^N)$ such that $\|\varphi-\psi\|_{\Di}<\eta$.
Since $\psi$ has compact support, from Hardy's inequality we have
that, for large $n$,
\begin{align*}
\bigg|&\int_{\R^N}V_{n}(x)\varphi_n(x)\varphi(x)\,dx-
\int_{\R^N} V(x)\varphi^2(x)\,dx\bigg|\\
&\leq\bigg|\int_{\R^N}(V_{n}(x)-V(x))\varphi_n(x)(\varphi(x)-\psi(x))\,dx\bigg|+
\bigg|\int_{\R^N}V(x)(\varphi_n(x)-\varphi(x))\varphi(x)\,dx\bigg|\\
&\leq {\rm const\,}\eta+
\bigg|\int_{\R^N}V(x)(\varphi_n(x)-\varphi(x))\varphi(x)\,dx\bigg|
={\rm const\,}\eta+o(1)
\quad\text{as }n\to\infty.
\end{align*}
(\ref{eq:82bi}) is thereby proved. From (\ref{eq:82bi}),  multiplying
 (\ref{eq:2}) by $\varphi$ and passing to
limit as $n\to\infty$, we obtain
$$
\mu(V)\int_{\R^N}|\n
\varphi(x)|^2\,dx\leq
\int_{\R^N}\big(|\n
\varphi(x)|^2-V(x)\varphi^2(x)\big)\,dx\leq 0,
$$
 a contradiction.
\end{pf}

The following theorem is the counterpart of Theorem \ref{t:pertinf} as far as
the possibility of perturbing singularities at dipoles is concerned.
\begin{Theorem}\label{t:pertpole}
  For $i=1,\dots,k$, let $r_i,R\in\R^+$, $a_i\in\R^N$, $a_i\neq a_j$
  for $i\neq j$, $h_i,h_\infty\in L^{\infty}\big({\mathbb S}^{N-1}\big)$ with
$\mu_1(h_i)>-(N-2)^2/4$, $\mu_1(h_\infty)>-(N-2)^2/4$,  and 
$$
V(x)=\sum_{i=1}^k\alchi_{B(a_i,r_i)}(x)\frac{
h_i\big(\frac{x-a_i}{|x-a_i|}\big)}{|x-a_i|^2}
      +\alchi_{\R^N\setminus B(0,R)}(x)
\frac{
h_\infty\big(\frac{x}{|x|}\big)}{|x|^2}+W(x)\in{\mathcal V}
$$
where $W\in L^{N/2}(\R^N)\cap L^{\infty}(\R^N)$. Assume that
$\mu(V)>0$ and let $h\in
L^{\infty}\big({\mathbb S}^{N-1}\big)$ such that, for some $i_0\in\{1,\dots,k\}$, 
$
\mu_1(h+h_{i_0})>-(N-2)^2/4,
$
 $\e>0$ such that, setting $H:=h+h_{i_0}$,
$$
\e<
\begin{cases}
\frac1{\max_{i=1,\dots,k,\infty}\{
\Lambda_N(h_i),
\Lambda_N(H)\}}-1,&\text{if }
\max_{i=1,\dots,k,\infty}\big\{
\Lambda_N(h_i),
\Lambda_N(H)\big\}>0,\\
+\infty&\text{if }
\max_{i=1,\dots,k,\infty}\big\{
\Lambda_N(h_i),
\Lambda_N(H)\big\}= 0,
\end{cases}
$$
 and  $\sigma>0$ such that 
$0<\frac1\sigma<\min_{i=1,\dots,k,\infty}
\Big\{\sqrt{(N-2)^2/4+\mu_1(\tilde h_i)},
\sqrt{(N-2)^2/4+\mu_1(\tilde H)}\Big\}$,
where $\tilde h_i:=(1+\e)h_i$ and $\tilde H:=(1+\e)H$.
  Then, for every $i\in\{1,\dots,k\}$, there exists
 $\bar r$ such that
$$
\mu\bigg(V+\frac{h\big(\frac{x-a_{i_0}}{|x-a_{i_0}|}\big)}
{|x-a_{i_0}|^2}\alchi_{{\mathcal E}^{\sigma, 
r}_{\tilde h_{i_0},\tilde H}(a_{i_0})}\bigg)>0,
$$
for all $0<r<\bar r$.
\end{Theorem}
\begin{pf}
Let $\sigma_1>0$ such that 
$$
\frac{N-2}2<\frac1{\sigma_1}<\frac1{\sigma_1}+\frac1{\sigma}<
\frac{N-2}2+
\min_{i=1,\dots,k,\infty}
\bigg\{\sqrt{\frac{(N-2)^2}{4}+\mu_1(\tilde h_i)},
\sqrt{\frac{(N-2)^2}{4}+\mu_1(\tilde H)}\bigg\}.
$$
Assume, by contradiction, that there exists a sequence $r_n\to0^+$
such that, setting 
$$
V_n= 
V+\frac{h\big(\frac{x-a_{i_0}}{|x-a_{i_0}|}\big)}
{|x-a_{i_0}|^2}
\alchi_{
{{\mathcal E}^{\sigma, 
r_n}_{\tilde h_{i_0},\tilde H}(a_{i_0})}},
$$ 
there holds $\mu(V_n)\leq 0$.  
By Lemmas \ref{l:sep_parte1} and \ref{l:sep_variante},
$V_n=\widetilde V_n+\widetilde W$, where $ \mu(\widetilde V_n)\geq \frac\e{\e+1}>0$
and
$$
\widetilde W(x)=
\sum_{i=1}^k\alchi_{B(a_i,r_i)\setminus 
 {\mathcal E}^{\sigma_1,\delta}_{0,\tilde h_i}(a_i)}(x)
\frac{h_i\big(\frac{x-a_i}{|x-a_i|}\big)}{|x-a_i|^2}
+ \alchi_{
{\mathcal E}^{\sigma_1,{R_0}}_{\tilde h_\infty,0}
\setminus B(0,R)}(x)\frac{
h_\infty\big(\frac{x}{|x|}\big)}{|x|^2}+W(x)
$$
for some  $\d>0$ and $R_0>0$ independent of $n$.
By Proposition \ref{p:attai}, $\mu(V_n)$ is attained by some $\varphi_n\in\Di$ 
such that $\int_{\R^N}|\nabla \varphi_n(x)|^2\,dx=1$, 
$\int_{\R^N}|\nabla \varphi_n(x)|^2\,dx-\int_{\R^N}V_n(x)\varphi_n^2(x)\,dx=\mu(V_n)$,
and \eqref{eq:2} is satisfied.
Up to a subsequence, $\varphi_n\weakly \varphi$ weakly in $\Di$
for some $\varphi\in\Di$, and hence 
\begin{align*}
\frac\e{\e+1}\leq\mu(V_n)+\int_{\R^N}\widetilde
W(x)\varphi_n^2(x)\,dx\leq \int_{\R^N}\widetilde
W(x)\varphi^2(x)\,dx+o(1)\quad\text{as }n\to+\infty.
\end{align*}
Therefore $\int_{\R^N}\widetilde
W(x)\varphi^2(x)\,dx>0$ and we conclude that $\varphi\not\equiv 0$.
We claim that 
\begin{equation}\label{eq:82tri}
\lim_{j\to\infty}\int_{\R^N} V_{n}(x)\varphi_n(x)\varphi(x)\,dx=
\int_{\R^N}V(x)\varphi^2(x)\,dx.
\end{equation}
Indeed for any $\eta>0$, by density there exists $\psi\in
C^{\infty}_{\rm c}(\R^N\setminus\{a_{i_0}\})$ such that $\|\varphi-\psi\|_{\Di}<\eta$.
Since the support of $\psi$ is detouched from $a_{i_0}$, from Hardy's inequality we have
that, for large $n$,
\begin{align*}
\bigg|&\int_{\R^N}V_{n}(x)\varphi_n(x)\varphi(x)\,dx-
\int_{\R^N} V(x)\varphi^2(x)\,dx\bigg|\\
&\leq\bigg|\int_{\R^N}(V_{n}(x)-V(x))\varphi_n(x)(\varphi(x)-\psi(x))\,dx\bigg|+
\bigg|\int_{\R^N}V(x)(\varphi_n(x)-\varphi(x))\varphi(x)\,dx\bigg|\\
&\leq {\rm const\,}\eta+
\bigg|\int_{\R^N}V(x)(\varphi_n(x)-\varphi(x))\varphi(x)\,dx\bigg|
={\rm const\,}\eta+o(1)
\quad\text{as }n\to\infty.
\end{align*}
(\ref{eq:82tri}) is thereby proved. From (\ref{eq:82tri}),  multiplying
 (\ref{eq:2}) by $\varphi$ and passing to
limit as $n\to\infty$, we obtain
$$
\mu(V)\int_{\R^N}|\n
\varphi(x)|^2\,dx\leq
\int_{\R^N}\big(|\n
\varphi(x)|^2-V(x)\varphi^2(x)\big)\,dx\leq 0,
$$
 a contradiction.
\end{pf}

\section{Localization of binding}\label{sec:localization-binding}

This section deals with the \emph{localization of binding} of
Schr\"odinger operators with potentials in the class ${\mathcal V}$.
The theory of localization of binding was first developed by Simon
\cite{Simon} who proved that if $V_1$ and $V_2$ are compactly supported
potentials such that the corresponding Schr\"odinger operators are
positive, then also the operator $-\Delta-V_1-V_2(\cdot-y)$ is positive definite  
provided $|y|$ is sufficiently large. Pinchover \cite{pinchover95}
extended the above result to the case of potentials belonging to the Kato class.
 
We notice that  both  Simon and Pinchover consider 
potentials which are lower order perturbations of the Laplacian, thus
excluding  the case of potentials with inverse
square type singularities. Such a case presents an additional difficulty due to the 
interaction of singularities at infinity. In \cite{FMT1}, 
 the following localization of binding result was proved for locally isotropic 
 inverse square potentials under some additional assumptions to control the
singularity at infinity:

\begin{Theorem}\label{t:scattering} 
\cite[Theorem 1.5]{FMT1} 
For $j=1,2$, let
$$
V_j=
\sum_{i=1}^{k_j}\alchi_{B(a_i^j,r_i^j)}(x)\frac{
h_i^j\big(\frac{x-a_i^j}{|x-a_i^j|}\big)}{|x-a_i^j|^2}
      +\alchi_{\R^N\setminus B(0,R_j)}(x)
\frac{
h_\infty^j\big(\frac{x}{|x|}\big)}{|x|^2}+W_j(x)\in{\mathcal V}.
$$
Assume that $\mu(V_j)>0$, the functions $h_i^j$ are constant for
all $j=1,2$, $i=1,\dots,k_j$, and that \eqref{eq:locbind1} is
satisfied.  Then, there exists $R>0$ such that, for every $y\in\R^N$
with $|y|\geq R$, the quadratic form associated to the operator
$-\Delta-\left(V_1+V_2(\cdot-y)\right)$ is positive definite.
\end{Theorem}

For general anisotropic inverse-square potentials, condition  (\ref{eq:locbind1})
is necessary to have a localization of binding type result, as the following 
proposition states.

\begin{Proposition}\label{p:nclb}
For $j=1,2$, let
$$
V_j=
\sum_{i=1}^{k_j}\alchi_{B(a_i^j,r_i^j)}(x)\frac{
h_i^j\big(\frac{x-a_i^j}{|x-a_i^j|}\big)}{|x-a_i^j|^2}
      +\alchi_{\R^N\setminus B(0,R_j)}(x)
\frac{
h_\infty^j\big(\frac{x}{|x|}\big)}{|x|^2}+W_j(x)\in{\mathcal V}.
$$
If there exists 
$y\in\R^N$ such that  $\mu(V_1+V_2(\cdot-y))>0$, then 
$$
\mu_1(h_\infty^1+h_\infty^2)>
-\bigg(\frac{N-2}2\bigg)^{\!\!2}.
$$
\end{Proposition}

\begin{pf}
  Assume that, 
  for some $y\in\R^N$, for some $\e>0$, and for any $u\in\Di$,
\begin{align}\label{eq:5}
\int_{\R^N}\big(|\n u(x)|^2-V_1(x)u^2(x)-V_2(x-y)u^2(x)\big)\,dx
\geq\e\int_{\R^N}|\n u(x)|^2\,dx.
\end{align}
Arguing by contradiction, suppose that
$\mu_1(h_\infty^1+h_\infty^2)\leq -(N-2)^2/4$,
 or, equivalently, that
$\La_N(h_\infty^1+h_\infty^2)\geq1$.  Let
$0<\delta<\e\Lambda_N(h_\infty^1+h_\infty^2)^{-1}$.  By \eqref{eq:lnh}
 and
density of $C^{\infty}_{\rm c}(\R^N\setminus\{0\})$ in $\Di$, 
there exists $\phi\in
C^{\infty}_{\rm c}(\R^N\setminus\{0\})$ such that  
$$
\int_{\R^N}|\n\phi(x)|^2dx-\int_{\R^N}\frac{(h_\infty^1+h_\infty^2)\big(\frac{x}{|x|}\big)}
{|x|^2}\,\phi^2(x)\,dx
<\delta\int_{\R^N}\frac{(h_\infty^1+h_\infty^2)\big(\frac{x}{|x|}\big)}
{|x|^2}\,\phi^2(x)\,dx.
$$
Let $\phi_{\mu}(x)=\mu^{-(N-2)/2}\phi(x/\mu)$. A direct calculation (similar
to that performed in the proof of Proposition \ref{t:mainresult}) yields
\begin{align*}
\int_{\R^N}&|\n\phi_{\mu}(x)|^2dx-\int_{\R^N}
V_1(x)\phi_{\mu}^2(x)\,dx-\int_{\R^N}
V_2(x-y)\phi_{\mu}^2(x)\,dx\\
&= \int_{\R^N}|\n\phi(x)|^2dx-
\int_{\R^N}\frac{h_\infty^1\big(\frac{x}{|x|}\big)}
{|x|^2}\,\phi^2(x)\,dx-
\int_{\R^N\setminus B\big(\frac{y}{\mu},\frac{R_2}{\mu}\big)}
\frac{h_\infty^2\big(\frac{x-(y/\mu)}{|x-(y/\mu)|}\big)}
{|x-(y/\mu)|^2}\,\phi^2(x)\,dx+o(1), 
\end{align*} 
 as $\mu\to\infty$. From continuity of $\varphi$ and the Dominated Convergence Theorem, 
 we deduce that
\begin{align*}
\int_{\R^N\setminus B\big(\frac{y}{\mu},\frac{R_2}{\mu}\big)}&
\frac{h_\infty^2\big(\frac{x-(y/\mu)}{|x-(y/\mu)|}\big)}
{|x-(y/\mu)|^2}\,\phi^2(x)\,dx=
\int_{\R^N\setminus B\big(0,\frac{R_2}{\mu}\big)}
\frac{h_\infty^2\big(\frac{x}{|x|}\big)}
{|x|^2}\,\phi^2\big(x+{\textstyle{\frac{y}{\mu}}}\big)\,dx\\
&=\int_{\R^N}\frac{h_\infty^2\big(\frac{x}{|x|}\big)}
{|x|^2}\,\phi^2(x)\,dx+o(1), \quad\text{as } \mu\to\infty,
\end{align*}
hence 
\begin{align*}
  \int_{\R^N}&|\n\phi_{\mu}(x)|^2dx-\int_{\R^N}
V_1(x)\phi_{\mu}^2(x)\,dx-\int_{\R^N}
V_2(x-y)\phi_{\mu}^2(x)\,dx\\
&= \int_{\R^N}|\n\phi(x)|^2dx-
\int_{\R^N}\frac{(h_\infty^1+h_\infty^2)\big(\frac{x}{|x|}\big)}
{|x|^2}\,\phi^2(x)\,dx+o(1), 
\end{align*} 
 as $\mu\to\infty$.
Letting $\mu\to\infty$ we obtain that
\begin{align*}
 \e \int_{\R^N}|\n\phi(x)|^2dx&\leq 
\int_{\R^N}|\n\phi(x)|^2dx-
\int_{\R^N}\frac{(h_\infty^1+h_\infty^2)\big(\frac{x}{|x|}\big)}
{|x|^2}\,\phi^2(x)\,dx\\
&<\delta \Lambda_N(h_\infty^1+h_\infty^2)\int_{\R^N}|\n\phi(x)|^2dx
<\e\int_{\R^N}|\n\phi(x)|^2dx,
\end{align*}
thus giving a contradiction.~\end{pf}

On the other hand, the control at infinity required in (\ref{eq:locbind1}) is no 
 sufficient to obtain a localization of binding result for potentials which are
 anisotropic at infinity, as enlightened by Example \ref{ex:controes} stated in the introduction: 

\smallskip
{\bf Example \ref{ex:controes}.\, }{\it For $N\geq 4$, there exist $V_1,V_2\in{\mathcal V}$ such that 
$\mu(V_1),\mu(V_2)>0$,  \eqref{eq:locbind1} holds, and for every $R>0$ there
exists $y_R\in\R^N$ such that $|y_R|>R$ and the quadratic form associated to 
the operator
$-\Delta-\left(V_1+V_2(\cdot-y_R)\right)$ is not positive semidefinite, i.e.
$\mu(V_1+V_2(\cdot-y_R))<0$.  }
\smallskip

\begin{pf}
For $N\geq 4$, let $\bar
y=(0,\dots,0,1)\in\R^N$ and
$\frac12\big(\frac{N-3}2\big)^{\!2}<\lambda<\big(\frac{N-3}2\big)^{\!2}$.  
From \cite{SSW}, scaling
invariance, translation invariance in the $\bar y$-direction, and 
density of $C^{\infty}_{\rm
    c}\big((\R^{N-1}\setminus\{0\})\times\R\big)$ in $\Di$, it
follows that, for
any $0<\e<2\lambda-\big(\frac{N-3}2\big)^2$,
there exists $\psi\in
C^{\infty}_{\rm
    c}\big((\R^{N-1}\setminus\{0\})\times\R\big)$
such that
$$
\mathop{\rm supp}\psi\subset Q:=\bigg\{x\in\R^N:\ \frac{x}{|x|}\cdot\bar
  y>\frac{\sqrt3}2\ \text{and}\ \frac{x-\bar y}{|x-\bar y|}\cdot\bar
  y<-\frac{\sqrt3}2\bigg\}
$$
and
$$
\bigg(\frac{N-3}2\bigg)^{\!\!2}<\frac{{\displaystyle{\int_{\R^N}|\nabla
      \psi(x)|^2\,dx}} } {{\displaystyle{\int_{\R^N}\frac{\psi^2(x)}{|x'|^2}\,dx}}}
<\bigg(\frac{N-3}2\bigg)^{\!\!2}+\e,
$$
where a generic point $x\in \R^N$  is denoted, from now on, as
$x=(x',x_N)\in\R^{N-1}\times\R$.
Since $\psi\in C^{\infty}_{\rm
    c}\big(\R^N\setminus\big\{x'=0\big\}\big)$, 
$$
\delta=\min\left\{\sqrt{1-\big({\textstyle{\frac{x}{|x|}}}\cdot\bar
  y\big)^{2}}, \sqrt{1-\big({\textstyle{\frac{x-\bar y}{|x-\bar y|}}}\cdot\bar
  y\big)^{2}}:\ x\in\mathop{\rm supp}\psi \right\}>0.
$$
Hence 
$$
\mathop{\rm supp}\psi\subset Q':=
\bigg\{x\in\R^N:\ \frac{\sqrt3}2<\frac{x}{|x|}\cdot\bar
  y<\sqrt{1-\delta^2}\ \text{and}\ -\sqrt{1-\delta^2}<
\frac{x-\bar y}{|x-\bar y|}\cdot\bar
  y<-\frac{\sqrt3}2\bigg\}
$$
Let 
$$
V_1(x)=\frac{\lambda\,\alchi_{C^+}(x)}{|x|^2\Big(1-(\frac{x}{|x|}\cdot \bar y)^2\Big)}
=\frac{\lambda\,\alchi_{C^+}(x)}{|x'|^2}\quad\text{and}\quad
V_2(x)=\frac{\lambda\,\alchi_{C^-}(x)}{|x|^2\Big(1-(\frac{x}{|x|}\cdot \bar y)^2\Big)}
=\frac{\lambda\,\alchi_{C^-}(x)}{|x'|^2},
$$
where
\begin{align*}
&C^+=\bigg\{x\in\R^N:\ \frac{\sqrt3}2<\frac{x}{|x|}\cdot\bar
  y<\sqrt{1-\delta^2}\bigg\}\quad\text{and}\\
&
C^-=\bigg\{x\in\R^N:\ -\sqrt{1-\delta^2}<
\frac{x}{|x|}\cdot\bar
  y<-\frac{\sqrt3}2\bigg\}=-C^+.
\end{align*}
We notice that $C^+\cap \big(\bar y+C^-\big)=Q'$. 
Moreover, we can write $V_1$, $V_2$ as
\begin{align*}
  & V_1(x)=\frac{h_1(\frac{x}{|x|})}{|x|^2}\text{\quad with\quad }h_1(\theta)=
  \frac{\lambda\,\alchi_{\{\theta\in{\mathbb
        S}^{N-1}:\sqrt{3}/2<\theta\cdot\bar
      y<\sqrt{1-\delta^2}\}}(\theta)} {\Big(1-(\theta\cdot \bar
    y)^2\Big)},\\
  & V_2(x)=\frac{h_2(\frac{x}{|x|})}{|x|^2}\text{\quad with\quad }h_2(\theta)=
  \frac{\lambda\,\alchi_{\{\theta\in{\mathbb
        S}^{N-1}:-\sqrt{1-\delta^2}<\theta\cdot\bar
      y<-\sqrt{3}/2\}}(\theta)} {\Big(1-(\theta\cdot \bar
    y)^2\Big)}.
\end{align*}
Being $\lambda<\big(\frac{N-3}2\big)^{\!2}$, from \cite{BT, SSW} it follows 
easily that $\mu(V_1),\mu(V_2)>0$. Moreover, $C^+\cap C^-=\emptyset$ and
$$
V_1(x)+V_2(x)=\frac{h_1(\frac{x}{|x|})+h_2(\frac{x}{|x|})}{|x|^2}
=\frac{\lambda\,\big(\alchi_{C^+}(x)+\alchi_{C^-}(x)\big)}{|x'|^2}
=\frac{\lambda\,\alchi_{C^+\cup C^-}(x)}{|x'|^2}\leq
\frac{\lambda}{|x'|^2},
$$
hence, from $\lambda<\big(\frac{N-3}2\big)^{\!2}$ and \cite{BT, SSW}, it follows 
 that $\mu(V_1+V_2)>0$, and thus 
$\mu_1(h_1+h_2)>-\big(\frac{N-2}2\big)^{\!2}$.

For any $\mu>0$, let $\psi_{\mu}(x):=\mu^{-{\frac{N-2}2}}\psi(x/\mu)$. 
Since
$$
V_1(x)+V_2(x-\bar
    y)=\frac{\lambda\,\big(\alchi_{C^+}(x)+\alchi_{\bar y+C^-}(x)\big)}{|x'|^2},
$$
a direct calculation yields, for all $\mu>0$,
\begin{align*}
 & \frac{\int_{\R^N}\big(V_1(x)+V_2(x-\mu\,\bar
    y)\big)\psi_\mu^2(x)\,dx}{\int_{\R^N}|\nabla\psi_\mu(x)|^2\,dx}
=\frac{\int_{\R^N}\big(V_1(x)+V_2(x-\bar
    y)\big)\psi^2(x)\,dx}{\int_{\R^N}|\nabla\psi(x)|^2\,dx}\\
&\qquad\geq 2\lambda\, \frac{\int_{Q'}\frac{\psi^2(x)}{|x'|^2}\,dx}
{\int_{\R^N}|\nabla\psi(x)|^2\,dx}=
2\lambda\, \frac{\int_{\R^N}\frac{\psi^2(x)}{|x'|^2}\,dx}
{\int_{\R^N}|\nabla\psi(x)|^2\,dx}>\frac{2\lambda}{\big(\frac{N-3}2\big)^2+\e}>1,
\end{align*}
thus implying that $\mu\big(V_1+V_2(\cdot-\mu\,\bar
    y)\big)<0$ for all $\mu>0$. Hence, for every $R>0$, it is enough to choose
$\mu>R$ and $y_R=\mu\,\bar y$ to obtain the example we are looking for.
\end{pf}

The above example justifies the need of the stronger control of the singularity
 at infinity required in the following theorem, which provides a positive
supersolution for the Schr\"odinger operator  $-\Delta-(V_1+V_2(\cdot-y))$
with $|y|$ sufficiently large.

\begin{Theorem}\label{t:separation} 
For $j=1,2$, let
$$
V_j=
\sum_{i=1}^{k_j}\alchi_{B(a_i^j,r_i^j)}(x)\frac{
h_i^j\big(\frac{x-a_i^j}{|x-a_i^j|}\big)}{|x-a_i^j|^2}
      +\alchi_{\R^N\setminus B(0,R_j)}(x)
\frac{
h_\infty^j\big(\frac{x}{|x|}\big)}{|x|^2}+W_j(x)\in{\mathcal V},
$$
where $W_j\in L^{\infty}(\R^N)$, $W_j(x)=O(|x|^{-2-\delta})$, with
$\delta>0$, as $|x|\to\infty$. 
Assume that $\mu(V_1),\mu(V_2)>0$,  and that 
\begin{equation}\label{eq:1}
\supess_{{\mathbb
        S}^{N-1}}(h_1^\infty)^+
+\supess_{{\mathbb
        S}^{N-1}}(h_2^\infty)^+
<\frac{(N-2)^2}4.
\end{equation}   
Then, there exists $R>0$ such that, for every $y\in\R^N$ with $|y|\geq
R$, 
there exists $\Phi_y\in\Di$,  $\Phi_y> 0$ and continuous in
$\R^N\setminus\{a_i^1,a_i^2+y\}_{i=1,\dots k_j,j=1,2}$, such that 
$$
-\Delta\Phi_y-\left(V_1+V_2(\cdot-y)\right)\Phi_y>0\quad\text{ in
}(\Di)^\star.
$$ 
\end{Theorem}
\begin{pf}
  Let us set $h_j=(\supess_{{\mathbb S}^{N-1}}h_j^\infty)-h_j^\infty$,
  $j=1,2$ and, for any $j=1,2$ and $i=1,\dots,k_j$, let us choose
  $\max\{0,\mu_1(h_i^j)\}<c_i^j< \mu_1(h_i^j)+(\frac{N-2}2)^2$.

In view of Theorems \ref{t:pertinf} and \ref{t:pertpole} and by assumption
  (\ref{eq:1}), there exist ${\mathcal N}_1$ and ${\mathcal N}_2$
  neighborhoods of infinity and ${\mathcal N}_i^j$ neighborhoods of $a_i^j$, 
$j=1,2$ and $i=1,\dots,k_j$,
 such that $\mu(\widetilde
  V_j)>0$, $j=1,2$, where  $\widetilde
  V_j(x):=V_j(x)+
\sum_{i=1}^{k_j}\alchi_{{\mathcal N}_i^j}(x)\frac{
c_i^j}{|x-a_i^j|^2}+
|x|^{-2}h_j(x/|x|)\alchi_{{\mathcal N}_j}(x)$.
We notice that $\widetilde
  V_j\geq V_j$ a.e. in $\R^N$.

Let $0<\e<(\frac{N-2}{2})^2-
\supess_{{\mathbb
        S}^{N-1}}(h_1^\infty)^+-\supess_{{\mathbb
        S}^{N-1}}(h_2^\infty)^+$  and, for $j=1,2$, set
$$
\Lambda=\bigg(\frac{N-2}{2}\bigg)^{\!\!2}-\e\quad\text{ and }
\quad\gamma_j=\Lambda-\supess_{{\mathbb
        S}^{N-1}}h_j^\infty.
$$
Let us also choose $0<\eta<\!\!<1$ such that 
\begin{equation}\label{eq:13}
\supess_{{\mathbb
        S}^{N-1}}h_2^\infty
<\gamma_1(1-2\eta)\quad\text{and}\quad
\supess_{{\mathbb
        S}^{N-1}}h_1^\infty
<\gamma_2(1-2\eta).
\end{equation}
We can choose $\bar R>0$ such that, for $j=1,2$,
$\bigcup_{i=1}^{k_j}B(a_i^j,r_i^j)\subset B(0,\bar R)$, and define
\begin{equation}\label{eq:scatt9}
p_j(x):=\begin{cases}
|x-a_i^j|^{-2+\tau}&\text{in }B(a_i^j,r_i^j),\\
1&\text{in }B(0,\bar R)\setminus\bigcup_{i=1}^{k_j}B(a_i^j,r_i^j),\\
0&\text{in }\R^N\setminus B(0,\bar R),
\end{cases}
\end{equation}
with $0<\tau<\min_{i,j}
\Big\{1,\frac{N-2}2-\sqrt{\big(\frac{N-2}2\big)^2+\mu_1(h_i^j)-c_i^j}\Big\}$.
 In view of Theorem \ref{t:pertinf}, there exist $\tilde
R_j$ such that the quadratic forms associated to the operators
$-\Delta-\widetilde V_j-\frac{\gamma_j}{|x|^2}\alchi_{\R^N\setminus
  B(0,\tilde R_j)}$ are positive definite.  Therefore, since $p_j\in
L^{N/2}$, the infima
$$
\mu_j=\inf_{u\in\Di\setminus\{0\}}\frac{\int_{\R^N}\big[|\n
  u(x)|^2-\widetilde
  V_j(x)u^2(x)-\gamma_j|x|^{-2}\alchi_{\R^N\setminus
    B(0,\tilde R_j)}u^2(x)\big]\,dx}{\int_{\R^N}p_j(x)u^2(x)\,dx}
$$
are achieved by some $\psi_j\in\Di$,  $\psi_j>0$ and continuous in
$\R^N\setminus\{a_1^j,\dots,a_{k_j}^j\}$, $\Di$-weakly solving 
\begin{equation}\label{eq:8}
  -\Delta\psi_j(x)-\widetilde V_j(x)\psi_j(x)=
 \mu_{j}p_j(x)\psi_j(x)+\frac{\gamma_j}{|x|^2}\alchi_{\R^N\setminus
    B(0,\tilde R_j)}\psi_j(x).
\end{equation}
From the asymptotic analysis of the exact behavior near the singularity 
of solutions to Schr\"odinger equations with inverse square potentials
proved in  \cite{FS3} (see also \cite{murata} and \cite{FMT2}),  there holds
$$
\lim_{|x|\to+\infty} \psi_j(x)|x|^{N-2+\sigma_{\Lambda}}=\ell_j>0, \quad
\text{where}\quad
\sigma_{\Lambda}:=-\frac{N-2}2+\sqrt{\bigg(\frac{N-2}2\bigg)^{\!\!2}-\Lambda},
$$
hence the function $\varphi_j:=\frac{\psi_j }{\ell_j}$ solves (\ref{eq:8})
and $\varphi_j(x)\sim|x|^{-(N-2+\sigma_{\Lambda})}$ at $\infty$. 
Then there exists $\rho>\max\{\tilde R_1,\tilde R_2,\bar R\}$ such that, in  $\R^N\setminus B(0,\rho)$, 
\begin{equation}\label{eq:14}
(1-\eta^2)|x|^{-(N-2+\sigma_{\Lambda})}\leq\varphi_j(x)\leq
(1+\eta)|x|^{-(N-2+\sigma_{\Lambda})} 
\end{equation}
and that 
\begin{equation}\label{eq:15}
|W_1(x)|\leq\eta\gamma_2|x|^{-2}\quad\mbox{ and }
\quad|W_2(x)|\leq\eta\gamma_1|x|^{-2}.
\end{equation}
Moreover, from \cite[Theorem 1.1]{FMT2},  we can deduce 
that for some positive constant $C$
\begin{align}\label{eq:16}
  \frac1C |x-a_i^j|^{\sigma_i^j}&\leq\varphi_j(x)\leq
  C|x-a_i^j|^{\sigma_i^j} \text{ in }B(a_i^j,r_i^j),\quad
  i=1,\dots,k_j,
\end{align}
where $\sigma_i^j=-\frac{N-2}2+\sqrt{\big(\frac{N-2}2\big)^{2}
+\mu_1(h_i^j)-c_i^j}$. 

For any $y\in\R^N$, let us consider the function
$$
\Phi_y(x):=\gamma_2\varphi_1(x)+\gamma_1\varphi_2(x-y)\in\Di.
$$
Then $\Phi_y$ satisfies, in the weak $\Di$-sense,
$$
  -\D\Phi_y-\big(\widetilde V_1+\widetilde
  V_2(\cdot-y)\big)\Phi_y=f
$$
where
\begin{align*}
f(x)=&\mu_{1}\gamma_2 p_1(x)\varphi_1(x)+
  \frac{\gamma_1\gamma_2}{|x|^2}\alchi_{\R^N\setminus
    B(0,\tilde R_1)}\varphi_1(x)
  +\mu_{2}\gamma_1 p_2(x-y)\varphi_2(x-y)\\
&+ \frac{\gamma_1
    \gamma_2}{|x-y|^2}\alchi_{\R^N\setminus B(y,\tilde R_2)}\varphi_2(x-y)
  -\gamma_1\widetilde V_1(x)\varphi_2(x-y)- \gamma_2\widetilde
  V_2(x-y)\varphi_1(x).
\end{align*}
From (\ref{eq:13}), (\ref{eq:14}), and (\ref{eq:15}), it follows that in
$\R^N\setminus\big(B(0,\rho)\cup B(y,\rho)\big)$
\begin{align*}
f(x)
&\geq\frac{\gamma_1\gamma_2}{|x|^2}\varphi_1(x)+
\frac{\gamma_1\gamma_2}{|x-y|^2}\varphi_2(x-y)\\ 
&\quad-\gamma_1\bigg(\frac{\supess_{{\mathbb
        S}^{N-1}}h_1^\infty}{|x|^2}+W_1(x)\bigg)
\varphi_2(x-y)-\gamma_2\bigg(\frac{\supess_{{\mathbb
        S}^{N-1}}h_2^\infty}{|x-y|^2}+
W_2(x-y)\bigg)\varphi_1(x)\\
&\geq\frac{\gamma_1\gamma_2(1-\eta^2)}{|x|^{N+\sigma_{\Lambda}}}+
\frac{\gamma_1\gamma_2(1-\eta^2)}{|x-y|^{N+\sigma_{\Lambda}}}-
\frac{\gamma_1\gamma_2(1-\eta^2)}{|x|^2|x-y|^{N-2+\sigma_{\Lambda}}}-
\frac{\gamma_1\gamma_2(1-\eta^2)}{|x-y|^2|x|^{N-2+\sigma_{\Lambda}}}\\
&=\gamma_1\gamma_2(1-\eta^2)\left(\frac{1}{|x|^{N+\sigma_{\Lambda}-2}}-
  \frac{1}{|x-y|^{N+\sigma_{\Lambda}-2}}\right)\left(\frac{1}{|x|^2}-
  \frac{1}{|x-y|^2}\right)\geq0. 
\end{align*}
For $|y|$ sufficiently large, $B(0,\rho)\cap B(y,\rho)=\emptyset$. In
$B(a_i^1,r_i^1)$, for $i=1,\dots,k_1$, from (\ref{eq:scatt9}),
(\ref{eq:14}),  (\ref{eq:16}) and the fact that 
$\sigma_i^j<-\tau$ for all $j=1,2$ and $i=1,\dots,k_j$,
 we have that
\begin{align*}
  f(x)
  &\geq\mu_{1}\gamma_2p_1(x)\varphi_1(x)+ \frac{\gamma_1\gamma_2}{|x-y|^2}\varphi_2(x-y)-\gamma_1\widetilde V_1(x)\varphi_2(x-y)- \gamma_2\widetilde
  V_2(x-y)\varphi_1(x)\\
  &\geq|x-a_i^1|^{-2+\tau+\sigma_i^j}\left[\frac{\mu_{1}\gamma_2}C+o(1)\right],
  \quad\mbox{ as } |y|\to\infty.
\end{align*}
In $B(0,\rho)\setminus\bigcup_{i=1}^{k_1}B(a_i^1,r_i^1)$, from
(\ref{eq:scatt9}), (\ref{eq:14}), (\ref{eq:15}) and since
$\varphi_1>c>0$, we obtain that
$$
f(x)\geq
\mu_{1}\gamma_c+o(1),\quad\mbox{ as }|y|\to\infty.
$$
In a similar way we can prove that, if $|y|$ is large enough, 
$f(x)>0$, a.e. in $B(y,\rho)$.
Since $\widetilde V_1(x)+\widetilde
V_2(x-y)\geq  V_1(x)+V_2(x-y)$, we conclude that 
$$
-\Delta\Phi_y-\left(V_1+V_2(\cdot-y)\right)\Phi_y>0\quad\text{  in }(\Di)^\star.
$$ 
\end{pf}

\begin{pfn}{Theorem \ref{t:localization}}
Let us fix $\e\in(0,1)$ such that, for $j=1,2$, 
\begin{align}
\e<\min\Bigg\{
2S\,\mu(V_j),  \frac{\mu(V_j)}{4}\Bigg[
\frac{4
\Big(\sum_{i=1}^{k_j}\|h_i^j\|_{L^\infty({\mathbb S}^{N-1})}+
\|h_\infty^j\|_{L^\infty({\mathbb S}^{N-1})}\Big)}{(N-2)^2}
+\frac{\|W_j\|_{L^{N/2}(\R^N)}}S
\Bigg]^{-1}\Bigg\}\label{eq:18}
\end{align}
and
\begin{align*}
\mu((1+\e)V_j)>0,\quad
(1+\e)\Big[\supess_{{\mathbb
        S}^{N-1}}(h_1^\infty)^+
+\supess_{{\mathbb
        S}^{N-1}}(h_2^\infty)^+\Big]
<\frac{(N-2)^2}4,
\end{align*}
see the proof of Lemma \ref{l:positivity_condition}. Fix $\widetilde R>0$ such
that
\begin{align}\label{eq:23}
\|W_j\alchi_{\R^N\setminus
  B(0,\widetilde R)}\|_{L^{N/2}(\R^N)}<\min\Big\{\frac\e4,\frac\e{4}\,S\Big\}.
\end{align}
Denoting $V_{j,\widetilde R}:=V_j-W_j\alchi_{\R^N\setminus B(0,\widetilde R)}$, from
(\ref{eq:18}) and (\ref{eq:23}), there results
\begin{multline*}
\int_{\R^N}\big(|\n u(x)|^2-V_{j,\widetilde R}(x)u^2(x)\big)\,dx\\
\geq \Big[\mu(V_j)-\|W_j\alchi_{\R^N\setminus
  B(0,\widetilde R)}\|_{L^{N/2}(\R^N)}S^{-1}\Big]
\int_{\R^N}|\n u(x)|^2\,dx\geq \frac{\mu(V_j)}2\int_{\R^N}|\n
u(x)|^2\,dx,
\end{multline*}
therefore, from (\ref{eq:18}), it follows
\begin{align*}
&\int_{\R^N}\big(|\n u(x)|^2-(1+\e)V_{j,\widetilde R}(x)
u^2(x)\big)\,dx\\
&\geq \left[\frac{\mu(V_j)}2-\e \Bigg(
\frac{4
\Big(\sum_{i=1}^{k_j}\|h_i^j\|_{L^\infty({\mathbb S}^{N-1})}+
\|h_\infty^j\|_{L^\infty({\mathbb S}^{N-1})}\Big)}{(N-2)^2}
+\frac{\|W_j\|_{L^{N/2}(\R^N)}}S
\Bigg)\right]\int_{\R^N}|\n
u(x)|^2\,dx\\
&\geq \frac{\mu(V_j)}4\int_{\R^N}|\n
u(x)|^2\,dx.
\end{align*}
Hence the potentials $(1+\e)V_{j,\widetilde R}$ satisfy the assumptions of
Lemma \ref{t:separation}, which yields, for $|y|$ sufficiently large,
the existence of  $\Phi_y\in\Di$,  $\Phi_y> 0$ and continuous in
$\R^N\setminus\{a_i^1,a_i^2+y\}_{i=1,\dots k_j,j=1,2}$, such that 
\begin{equation*}
-\Delta\Phi_y-V_{\widetilde R,y}\Phi_y>
\e\, V_{\widetilde R,y}\Phi_y(x)
\quad\text{  in }(\Di)^\star,
\end{equation*}
where $V_{\widetilde R,y}(x):=V_{1,\widetilde R}(x)+V_{2,\widetilde
  R}(x-y)$. From Lemma \ref{l:positivity_condition}, we deduce that
$\mu(V_{\widetilde R,y})\geq\frac{\e}{\e+1}$.  Hence, from (\ref{eq:23}), for any
$u\in\Di$, there holds
\begin{align*}
&\int_{\R^N}\big(|\n u(x)|^2-(V_1(x)+V_2(x-y))
u^2(x)\big)\,dx\\
&=\int_{\R^N}\!\!\big(|\n u(x)|^2-V_{\widetilde R,y}(x)u^2(x)\big)\,dx
-\int_{\R^N\setminus B(0,\widetilde R)}\!\!W_1(x)u^2(x)\,dx
-\int_{\R^N\setminus B(y,\widetilde R)}\!\!W_2(x-y)u^2(x)\,dx\\
&\geq \Big[\frac{\e}{\e+1}-\|W_1\alchi_{\R^N\setminus
  B(0,\widetilde R)}\|_{L^{N/2}(\R^N)}S^{-1}-\|W_2\alchi_{\R^N\setminus
  B(0,\widetilde R)}\|_{L^{N/2}(\R^N)}S^{-1}\Big]\int_{\R^N}|\n u(x)|^2\,dx\\
&\geq \Big[\frac{\e}{\e+1}-\frac{\e}2\Big]\int_{\R^N}|\n u(x)|^2\,dx
=\frac{\e(1-\e)}{2(\e+1)}\int_{\R^N}|\n u(x)|^2\,dx.
\end{align*}
Therefore $\mu(V_1+V_2(\cdot-y))\geq \frac{\e(1-\e)}{2(\e+1)}>0$ for $|y|$ sufficiently large.
 The theorem is thereby proved.~\end{pfn}

\section{Essential self-adjointness}\label{sec:semi-bound-essent}

By the Shattering Lemma \ref{l:sep}, 
Schr\"odinger operators with potentials in  ${\mathcal
V}$ are  compact perturbations of positive operators,
see Theorem \ref{l:semibounded}. As a consequence, they are {\em semi-bounded
symmetric operators} and their $L^2(\R^N)$-spectrum  is bounded from below. 

In the present section, we discuss 
 essential self-adjointness of operators $-\D-V$ on the domain
$C^{\infty}_{\rm c}(\R^N\setminus\{a_1,\dots,a_k\})$, for every $V\in\mathcal V$.

From Theorem \ref{l:semibounded} (see also Lemmas \ref{l:sep_parte1} and
\ref{l:sep_variante}), 
we can split any $V\in{\mathcal V}$ as  $V(x)=\widetilde
V(x)+\widetilde W(x)$ where 
\begin{equation}\label{eq:37}
\widetilde
V(x)=\sum_{i=1}^k\alchi_{{\mathcal U}_i}(x)
\frac{h_i\big(\frac{x-a_i}{|x-a_i|}\big)}{|x-a_i|^2}
+\alchi_{{\mathcal U}_\infty}(x)
\frac{h_\infty\big(\frac{x}{|x|}\big)}{|x|^2},
\qquad \mu(\widetilde V)>0, 
\end{equation}
${\mathcal U}_i\subset B(a_i,1/2)$ is a neighborhood of $a_i$ for every 
$i=1,\dots,k$, ${\mathcal U}_\infty$ is a neighborhood of $\infty$, 
and 
$$
\widetilde W\in L^{N/2}(\R^N)\cap L^{\infty}(\R^N).
$$
  The
following self-adjointness criterion in ${\mathcal V}$ is an easy
consequence of  the Kato-Rellich Theorem (see e.g. \cite[Theorem 4.4]{kato})  and
well-known self-adjointness criteria for positive operators.

\begin{Lemma}{\bf [\,Self-adjointness  criterion in 
${\mathcal V}$\,]}\label{l:sac} 
Let $V\in{\mathcal V}$ and $V=\widetilde V+\widetilde W$, with
$\widetilde V$ as in~(\ref{eq:37}) and $\widetilde W\in
L^{N/2}(\R^N)\cap L^{\infty}(\R^N)$. Then the operator $-\D-V$ is
essentially  
self-adjoint in $C^{\infty}_{\rm c}(\R^N\setminus\{a_1,\dots,a_k\})$
if and only if $\mathop{\rm Range}(-\D-\widetilde V+b)$ is dense in
$L^2(\R^N)$ for some $b\in L^{\infty}(\R^N)$ with $\infess_{\R^N} b>0$.
\end{Lemma}

\noindent As a consequence of the above criterion, the  following {\em{non
self-adjointness  condition}} in ${\mathcal V}$ holds. We refer to \cite{FMT1} for 
more details.
  
\begin{Corollary}\label{l:nsac} 
Let $V\in{\mathcal V}$ and $V=\widetilde V+\widetilde W$, with
$\widetilde V$ as in (\ref{eq:37}) and $\widetilde W\in
L^{N/2}(\R^N)\cap L^{\infty}(\R^N)$. Assume that there exist $v\in
L^2(\R^N)$, $v(x)\geq0$ a.e. in $\R^N$, $\int_{\R^N}v^2>0$, a distribution 
$h\in H^{-1}(\R^N)$, and $\beta>0$ such that 
\begin{equation}\label{eq:47}
{}_{H^{-1}(\R^N)}\big\langle
h,u\big\rangle_{H^{1}(\R^N)}\leq 0\quad\text{for all }u\in H^1(\R^N),
\ u\geq 0\ \text{a.e in }\R^N,
\end{equation}
and
\begin{equation*}
-\D v-\tilde V\,v+\beta \, v=h\quad\text{in }{\mathcal D}'(\R^N\setminus\{a_1,\dots,a_k\}).
\end{equation*}
Then the operator $-\D-V$ is not
essentially  
self-adjoint in $C^{\infty}_{\rm c}(\R^N\setminus\{a_1,\dots,a_k\})$.
\end{Corollary}

\noindent We now prove Theorem \ref{t:self-adjo}.

\bigskip\noindent

\begin{pfn}{Theorem \ref{t:self-adjo}}

\medskip\noindent {\bf Step 1: } if
$\mu_1(h_i)\geq-\big(\frac{N-2}2\big)^{\!2}+1$ for all
$i\in\{1,\dots,k\}$, then $-\D-V$ is essentially self-adjoint.

\smallskip\noindent
For $i=1,\dots,k$, let $\phi_i\in C^\infty_{\rm c}(\R^N)$ such that $\phi_i\equiv 1$ in 
${\mathcal U}_i$, $\phi_i\equiv 0$ in $\R^N\setminus B(a_i,1)$, and $0\leq\phi_i\leq 1$, and
let $\phi_\infty\in C^\infty(\R^N)$ such that $\mathop{\rm supp}\phi_{\infty}
\subset \R^N\setminus\{0\}$,
$\phi_\infty\equiv 1$ in 
${\mathcal U}_\infty$, and  $0\leq\phi_\infty\leq 1$. Then, setting
$$
\overline{V}(x)=\sum_{i=1}^k
\frac{h_i\big(\frac{x-a_i}{|x-a_i|}\big)}{|x-a_i|^2}\phi_i(x)
+\frac{h_\infty\big(\frac{x}{|x|}\big)}{|x|^2}\phi_\infty(x),
$$
we have that 
$\widetilde V-\overline{V}\in L^\infty(\R^N)$. 

For $\alpha>\max\{0,-
\infess_{\R^N}(\widetilde V-\overline{V})\}$ and $b(x):=\widetilde V-\overline{V}+\alpha$,
there holds $b\in L^{\infty}(\R^N)$ and $\infess_{\R^N}b>0$.
Let us fix  $f\in C^{\infty}_{\rm c}(\R^N\setminus\{a_1,\dots,a_k\})$. 
Since $C^{\infty}_{\rm
  c}(\R^N\setminus\{a_1,\dots,a_k\})$ is dense in $L^2(\R^N)$, in
view of Lemma \ref{l:sac}, 
to prove essential self-adjointness it is enough to find some $g\in \mathop{\rm
Range}(-\D-\widetilde V+b)=\mathop{\rm
Range}(-\D-\overline{V}+\alpha)$ such that $g$ is arbitrarily close to
$f$ in $L^2(\R^N)$. To
this aim, we fix $\e>0$ and  claim that
there exist $\tilde h_i
\in C^{\infty}({\mathbb S}^{N-1})$,  $i=1,\dots,k,\infty$, and   and $\gamma_i\in\R$,
  $i=1,\dots,k$,
 such that, setting 
\begin{align*}
\widehat V(x)&:=
\sum_{i=1}^k
\frac{\tilde h_i\big(\frac{x-a_i}{|x-a_i|}\big)-\gamma_i}{|x-a_i|^2}
\phi_i(x)
+\phi_\infty(x)\frac{\tilde h_\infty\big(\frac{x}{|x|}\big)}{|x|^2}, 
\end{align*}
there holds
\begin{align}
  \label{eq:20}& \int_{\R^N}\big(|\nabla u|^2-\widehat
  V\,u^2+\alpha \,u^2\big)\,dx \geq \min\Big\{\frac12\mu(\tilde
  V), \infess_{\R^N}b\Big\}\|u\|^2_{H^1(\R^N)},\quad\text{for all }
  u\in H^1(\R^N),\\
\label{eq:27}& \mu_1(\tilde h_i)+\gamma_i>
-\big({\textstyle{\frac{N-2}2}}\big)^2+1, \quad\text{for all }i=1,\dots,k,\\
\label{eq:53}&
\|(\widehat V-\overline{V})u\|_{L^2(\R^N)}<\e,
\end{align}
for every $u\in H^1(\R^N)$ solving 
\begin{equation}\label{eq:54}
-\Delta u(x)-\widehat V(x)\, u(x)+\alpha\,u(x)=f(x).
\end{equation}
In order to prove the claim, for every $i=1,\dots,k,\infty$ let 
$h_n^i\in C^{\infty}({\mathbb S}^{N-1})$ such that  $h_n^i\to
h_i$ a.e. in ${\mathbb S}^{N-1}$, and $|h_n^i(\theta)|\leq
\|h_i\|_{L^{\infty}({\mathbb S}^{N-1})}$ for a.e. $\theta\in {\mathbb
  S}^{N-1}$. We notice that the existence of such approximating sequences 
can be proved using convolution methods in local charts.
From Lemma \ref{l:app}, it follows that
$\lim_{n\to+\infty} \mu_1(h_n^i)=\mu_1(h_i)$, while,
from the Dominated Convergence Theorem, $\lim_{n\to+\infty}
\int_{{\mathbb S}^{N-1}}(h_n^i-h_i)^2=0$, 
hence, for any
$i\in\{1,\dots,k\}$ such that
$\mu_1(h_i)=-\big(\frac{N-2}2\big)^{\!2}+1$, there exists a sequence
$\{\delta_n^i\}_{n\in\N}\subset(0,\infty)$ such that
$$
\delta^i_n>0,\quad \lim_{n\to+\infty}\delta_n^i=0, \quad
\text{and }\delta_n^i>\max\Big\{0,\Big|1-\Big(\frac{N-2}2\Big)^{\!2}-\mu_1(h_n^i)\Big|,
\sqrt{\textstyle{\int_{{\mathbb S}^{N-1}}(h_n^i-h_i)^2}}
\Big\}.
$$
Hence, setting, for $i=1,\dots,k$,
$$
\gamma^i_n:=
\begin{cases}
  0,&\text{if }\mu_1(h_i)>-\big(\frac{N-2}2\big)^{\!2}+1,\\[7pt]
  -\mu_1(h_n^i)-\big(\frac{N-2}2\big)^{\!2}+1+\delta_n^i,&\text{if
  }\mu_1(h_i)=-\big(\frac{N-2}2\big)^{\!2}+1,
\end{cases}
$$
and
\begin{align*}
\sigma_n^i&=-\frac{N-2}{2}+\sqrt{\frac{(N-2)^2}{4}+\mu_1(h_n^i)+\gamma_n^i}\\[7pt]
&=
\begin{cases}
-\frac{N-2}{2}+\sqrt{1+\delta_n^i},
&\text{if }\mu_1(h_i)=-\left(\frac{N-2}2\right)^2+1,\\[7pt]
-\frac{N-2}{2}+\sqrt{\frac{(N-2)^2}{4}+\mu_1(h_n^i)},
&\text{if }\mu_1(h_i)>-\left(\frac{N-2}2\right)^2+1,
\end{cases}
\end{align*}
there holds $\sigma_n^i+\frac{N-2}2>1$ for $n$ sufficiently large.
Let us set 
\begin{align}
\widetilde V_n(x)&:=
\sum_{i=1}^k
\frac{\tilde h_n^i\big(\frac{x-a_i}{|x-a_i|}\big)-\gamma_n^i}{|x-a_i|^2}
\phi_i(x)
+\phi_\infty(x)\frac{\tilde h_n^\infty\big(\frac{x}{|x|}\big)}{|x|^2}.
\end{align}
From \cite[Theorem 1.2]{FMT2} and taking into account i--ii) of Lemma \ref{l:app}, 
there exists a positive constant $C$
independent on $n$, such that every  solution $u\in H^1(\R^N)$ of equation
\begin{equation}\label{eq:28}
-\Delta u(x)-\widetilde V_n(x)\, u(x)+\alpha\,u(x)=f(x)
\end{equation}
can be estimated as 
$$
|u(x)|\leq C\,|x-a_i|^{\sigma_n^i}\|u\|_{H^1(\R^N)}\quad
\text{in }B(a_i,1),\quad\text{for all }i=1,\dots,k.
$$
 From iii) of Lemma \ref{l:app}, it follows that $\widetilde V_n\to \overline{V}$ 
in $L^{N/2,\infty}$, hence
it follows that  
\begin{align*}
\int_{\R^N}\big(|\nabla u|^2-\widetilde
  V_n\,u^2+\alpha \,u^2\big)\,dx &\geq 
\left[\min\Big\{\mu(\widetilde
  V), \infess_{\R^N}b\Big\}+o(1)\right]\|u\|_{H^1(\R^N)}^2\\
&\geq \frac12 \min\Big\{\mu(\widetilde
  V), \infess_{\R^N}b\Big\}\|u\|_{H^1(\R^N)}^2
\end{align*}
for all $u\in H^1(\R^N)$. Moreover, by testing (\ref{eq:28}) with $u$, 
every  solution $u\in H^1(\R^N)$ of (\ref{eq:28})  satisfies
$$
\|u\|_{H^1(\R^N)}\leq\frac{2\|f\|_{L^2(\R^N)}}{\min\{\mu(\widetilde
  V),\infess_{\R^N} b\}}.
$$
Then for every 
solution $u\in H^1(\R^N)$ of (\ref{eq:28}),
there holds 
\begin{align*}
&\left\|\frac{(h_n^\infty-h_\infty)\big(\frac{x}{|x|}\big)\phi_\infty
u}{|x|^2}\right\|_{L^2(\R^N)}^2
\leq 
\bigg(\int_{\R^N}|u(x)|^{2^*}\,dx\bigg)^{2/2^*}
\bigg(\int_{\R^N}
\frac{|h_n^\infty-h_\infty|^N\big(\frac{x}{|x|}\big)\phi_\infty^N(x)}{|x|^{2N}}\,dx
\bigg)^{2/N}\\
&\leq S^{-1}\frac{4\|f\|^2_{L^2(\R^N)}}{(\min\{\mu(\widetilde
  V),\infess_{\R^N} b\})^2}
\bigg(\int_{\R^N}
\frac{|h_n^\infty-h_\infty|^N\big(\frac{x}{|x|}\big)\phi_\infty^N(x)}{|x|^{2N}}\,dx
\bigg)^{2/N}=o(1)
\end{align*}
as $n\to+\infty$.
Furthermore for all $i=1,\dots,k,\infty$ and for any solution $u\in H^1(\R^N)$ of (\ref{eq:28}),
there holds 
\begin{align*}
&\left\|\frac{(h_n^i-h_i)\big(\frac{x-a_i}{|x-a_i|}\big)\phi_i
u}{|x-a_i|^2}\right\|_{L^2(\R^N)}^2
\leq C^2\|u\|_{H^1(\R^N)}^2\bigg(\int_{{\mathbb
  S}^{N-1}}(h_n^i(\theta)-h_i(\theta))^2\,dV(\theta)\bigg)
\int_0^1
r^{N-5+2\sigma_n^i}\,dr.
\end{align*}
Hence, if $\mu_1(h_i)>-\left(\frac{N-2}2\right)^2+1$, we have that
\begin{multline*}
\left\|\frac{(h_n^i-h_i)\big(\frac{x-a_i}{|x-a_i|}\big)\phi_i
u}{|x-a_i|^2}\right\|_{L^2(\R^N)}^2 \leq
C^2
\frac{4\|f\|^2_{L^2(\R^N)}}{(\min\{\mu(\widetilde
  V),\infess_{\R^N} b\})^2}\times\\
\times
\bigg(\int_{{\mathbb
  S}^{N-1}}(h_n^i(\theta)-h_i(\theta))^2\,dV(\theta)\bigg)
\bigg[\frac1{2\Big(\sqrt{\frac{(N-2)^2}{4}+\mu_1(h_i)}-1\Big)}+o(1)\bigg]=
o(1),
\end{multline*}
as $n\to+\infty$, while,  if $\mu_1(h_i)=-\left(\frac{N-2}2\right)^2+1$, by the choice of
$\delta_n^i$, there holds
\begin{multline*}
\left\|\frac{(h_n^i-h_i)\big(\frac{x-a_i}{|x-a_i|}\big)\phi_i
u}{|x-a_i|^2}\right\|_{L^2(\R^N)}^2\leq
C^2\|u\|_{H^1(\R^N)}^2\bigg(\int_{{\mathbb
  S}^{N-1}}(h_n^i(\theta)-h_i(\theta))^2\,dV(\theta)\bigg)
\frac{\sqrt{1+\delta_n^i}+1}{2\delta_n^i}\\
\leq
\frac12 C^2\frac{4\|f\|^2_{L^2(\R^N)}}{(\min\{\mu(\widetilde
  V),\infess_{\R^N} b\})^2}
\|h_n^i-h_i\|_{L^2({\mathbb S}^{N-1})}
\big(\sqrt{1+\delta_n^i}+1)=o(1)
\end{multline*}
as $n\to+\infty$. 
Furthermore, for all $i\in\{1,\dots,k\}$
 such that $\mu(h_i)=-\left(\!\frac{N-2}2\!\right)^2+1$, we have
that
\begin{align*}
\left\|\frac{\gamma_n^i\phi_iu}{|x-a_i|^2}\right\|_{L^2(\R^N)}^2
&\leq C^2(\gamma_n^i)^2\|u\|_{H^1(\R^N)}^2\int_0^1
r^{N-5+2\sigma_n^i}\,dr\\
&\leq\frac{2\,C^2\,\|f\|^2_{L^2(\R^N)}}{\big(\min\{\mu(\widetilde
  V),\infess_{\R^N}b\}\big)^2}\,\frac{(\gamma_n^i)^2\big(\sqrt{1+\delta_n^i}+1\big)}
{\delta_n^i}=
o(1)
\end{align*}
as $n\to+\infty$. Therefore, it is possible to choose $n$ large  enough in order to
ensure that every solution $u$ of (\ref{eq:28}) satisfies 
$$
\|(\widetilde V_n-\overline{V})u\|_{L^2(\R^N)}<\e.
$$
For such an $n$, let us set $\widehat V=\widetilde V_n$, $\tilde h_i=
\tilde h_n^i$,   $i=1,\dots,k,\infty$,  and $\gamma_i=\gamma_n^i$,
  $i=1,\dots,k$, so that conditions (\ref{eq:20}--\ref{eq:53}) 
are satisfied and the claim is proved.

Hence, by the Lax-Milgram Theorem, there exists  $w\in H^1(\R^N)$ satisfying 
\begin{equation}\label{eq:9}
-\D w(x)-\widehat V(x)w(x)+\alpha\,w(x)=f(x). 
\end{equation} 
Since $\widehat V$ is smooth outside poles,
 by classical regularity
theory $w\in C^{\infty}(\R^N\setminus\{a_1,\dots,a_k\})$.
From  \cite[Theorem 1.1]{FMT2} we deduce the following asymptotic behavior of
$w$ at poles  
\begin{align}\label{eq:38}
  w(x)\sim |x-a_i|^{\sigma_{i}} \psi_1^{\hat h_i}\Big(\frac{x-a_i}{|x-a_i|}\Big),\text{ as }
x\to a_i,
\end{align}
where 
$$
\hat h_i:=\tilde h_i-\gamma_i\quad\text{and}\quad 
\sigma_i=-\frac{N-2}{2}+\sqrt{\frac{(N-2)^2}{4}+\mu_1(\tilde h_i)+\gamma_i}.
$$
Hence the function
$g(x):=\widehat V(x)w(x)-\alpha\,w(x)+f(x)$ satisfies
$$
|g(x)|\leq {\rm const}\,\frac{\psi_1^{\hat h_i}
\big(\frac{x-a_i}{|x-a_i|}\big)}{|x-a_i|^{2-\sigma_{i}}}\,
\quad\text{ in }{\mathcal U}_i\quad\text{for all }i=1,\dots,k.
$$ 
Since, by (\ref{eq:27}), $\sigma_i+\frac{N-2}2>1$, it turns out that  $g\in
L^2(\R^N)$. 
Writing $w$ by its Green's representation formula  and 
using well known properties of differentiability of {\em{Newtonian potentials}}, 
the following asymptotic estimate for the gradient of $w$ near the poles 
can be deduced
\begin{equation}\label{eq:39}
\n w(x)=
\begin{cases}
O\big(|x-a_i|^{\sigma_{i}-1}\big),&\text{if }\mu_1(\hat h_i)<N-1,\\
O\big(|x-a_i|^{-\tau}\big),&\text{if }\mu_1(\hat h_i)\geq N-1,
\end{cases}
\qquad \text{ as }x\to a_i,
\end{equation}
for some $0<\tau<\frac{N-2}2$, see \cite{FMT1} for more details in the Hardy case.

For all $n\in\N$ let $\eta_n$ be a cut-off function such
that $\eta_n\in  C^{\infty}_{\rm c}(\R^N\setminus\{a_1,\dots,a_k\})$,
$0\leq\eta_n\leq 1$, and
\begin{align*}
&\eta_n(x)\equiv 0\text{ in }
\bigcup_{i=1}^kB\Big(a_i, \frac1{2n}\Big)\cup \big(\R^N\setminus B(0,2n)\big),\quad
\eta_n(x)\equiv 1\text{ in }B(0,n)\setminus \bigcup_{i=1}^k B\Big(a_i,
\frac1{n}\Big),\\
&|\n \eta_n(x)|\leq C\,n\text{ in } \bigcup_{i=1}^k \Big(B\Big(a_i,
\frac1{n}\Big)\setminus B\Big(a_i, \frac1{2n}\Big)\Big),\quad |\n
\eta_n(x)|\leq \frac{C}n \text{ in }B(0,2n)\setminus B(0,n),\\
&|\D \eta_n(x)|\leq C\,n^2\text{ in } \bigcup_{i=1}^k \Big(B\Big(a_i,
\frac1{n}\Big)\setminus B\Big(a_i, \frac1{2n}\Big)\Big),\quad |\D
\eta_n(x)|\leq \frac{ C}{n^2} \text{ in }B(0,2n)\setminus B(0,n),
\end{align*}
for some positive constant $C$ independent of $n$.
Setting $f_n:=\eta_n f-2\n \eta_n\cdot\n w-w\,\D\eta_n$,
we notice that $\eta_n f\to f$ in $L^2(\R^N)$, while (\ref{eq:38})
and (\ref{eq:39}) yield
\begin{align*}
  &\int_{\R^N}|\n \eta_n(x)|^2|\n w(x)|^2\,dx\\
  &\leq {\rm const\,} n^2\sum_{i=1}^k\int_{B(a_i, \frac1{n})\setminus
    B(a_i,\frac1{2n})}|\n w(x)|^2\,dx+\frac{\rm const\,}{n^2}
  \int_{B(0,2n)\setminus B(0,n)}|\n w(x)|^2\,dx\\
  &\leq{\rm const\,} n^2\left[
    \sum_{\mu_1(\hat h_i)<N-1}\int_{B(0,\frac1n)}|x|^{2\sigma_{i}-2}\,dx 
+ \sum_{\mu_1(\hat h_i)\geq N-1}\int_{B(0,\frac1n)}|x|^{-2\tau}\,dx\right]
+\frac{\rm
      const\,}{n^2}\|w\|_{H^1(\R^N)}\\
&\leq {\rm
      const\,}\bigg[\sum_{\mu_1(\hat h_i)<N-1} n^{-2\sigma_{i}+4-N}
+\sum_{\mu_1(\hat h_i)\geq N-1} n^{2+2\tau-N}
+n^{-2}\bigg]
\end{align*}
and
\begin{align*}
\int_{\R^N}&|\D \eta_n(x)|^2|w(x)|^2\,dx\\
&\leq {\rm const\,}n^4\sum_{i=1}^k\int_{B(a_i,
\frac1{n})\setminus B(a_i, \frac1{2n})}|w(x)|^2\,dx+\frac{\rm const\,}{n^4}
\int_{B(0,2n)\setminus B(0,n)}|w(x)|^2\,dx\\
&\leq  {\rm const\,}
n^4\sum_{i=1}^k\int_{B(0,1/n)}|x|^{2\sigma_{i}}\,dx+
\frac{\rm const\,}{n^4}\|w\|_{H^1(\R^N)}
\leq{\rm const\,}\bigg[\sum_{i=1}^k n^{-2\sigma_{i}+4-N}+n^{-4}\bigg].
\end{align*}
Since   $-2\sigma_{i}+4-N<0$, we
conclude that $f_n\to f$ in $L^2(\R^N)$.
 Hence, for $n$ large enough,
$\|f_n-f\|_{L^2(\R^N)}<\e$. 
The functions $w_n:=\eta_n w\in C^{\infty}_{\rm c}(\R^N\setminus\{a_1,\dots,a_k\})$
solve  
$$
-\D w_n(x)-\widetilde V(x) w_n(x)+b\,w_n(x)=g_n(x),
$$
where $g_n(x):=f_n(x)+(\widehat
V(x)-\overline{V}(x))w_n(x)$, i.e. $g_n\in  \mathop{\rm
  Range}(-\D-\widetilde V+b)$.
Moreover, from (\ref{eq:53}) and  $|\eta_n|\leq 1$, we deduce that
$$
\|g_n-
f_n\|_{L^2(\R^N)}=\|(\widehat V-\overline{V})w_n\|_{L^2(\R^N)}\leq
\|(\widehat V-\overline{V})w\|_{L^2(\R^N)}<\e,
$$ 
hence $\|g_n-
f\|_{L^2(\R^N)}<2\e$ for 
large $n$. The proof of step 1 is thereby complete.

\medskip\noindent {\bf Step 2: } if $\mu(h_i)<-\big(\frac{N-2}2\big)^{\!2}+1$ for some
$i\in\{1,\dots,k\}$, then $-\D-V$  is not essentially self-adjoint.

\smallskip\noindent Let $V=\widetilde V+\widetilde W$, with
$\widetilde V$ as in (\ref{eq:37}) and $\widetilde W\in
L^{N/2}(\R^N)\cap L^{\infty}(\R^N)$. Let us fix $i\in\{1,\dots,k\}$
such that $\mu(h_i)<-\big(\frac{N-2}2\big)^{\!2}+1$, $\beta>0$, and $\a<0$, and
consider the solution $\varphi\in C^1\big((-\infty,\ln\d]\big)$ of the
Cauchy problem
$$
\begin{cases}
\varphi''(s)-\omega_{h_i}^2\,\varphi(s)=\beta\, e^{2s}\,\varphi(s),\\
\varphi(\ln \d)=0,\quad
\varphi'(\ln \d)=\a,
\end{cases}
$$
where $\omega_{h_i}:=\sqrt{\big(\frac{N-2}2\big)^2+\mu_1(h_i)}$ and
$\delta$ is such that $B(a_i,\delta)\subset{\mathcal U}_i$. From the
Gronwall's inequality (see \cite[Lemma A.2]{FMT1} for more details),
we can estimate $\varphi$ as
\begin{equation}\label{eq:48}
0\leq \varphi(s)\leq C\,
\,e^{-\omega_{h_i} s}\quad\text{for all }s\leq\ln\d,
\end{equation}
for some positive constant $C=C(h_i,\d,\a,\beta)$. Let us set
$$
v(x):=
\begin{cases}
  |x-a_i|^{-\frac{N-2}2}\,\varphi(\ln|x-a_i|)
\,\psi_1^{h_i}\big(\frac{x-a_i}{|x-a_i|}\big)
,&\text{if }x\in B(a_i,\d)\setminus\{a_i\},\\
0,&\text{if }x\in \R^N\setminus\overline{B(a_i,\d)}.
\end{cases}
$$
From (\ref{eq:48}) we infer that 
\begin{equation}\label{eq:49}
0\leq v(x)\leq
C|x-a_i|^{-\frac{N-2}2-\sqrt{\left(\frac{N-2}2\right)^2+\mu(h_i)}}\quad\text{in }
B(a_i,\d). 
\end{equation}
The assumption  $\mu(h_i)<-\big(\frac{N-2}2\big)^{\!2}+1$ 
and estimate (\ref{eq:49}) ensure
that $v\in L^2(\R^N)$. Moreover the restriction of $v$ to $B(a_i,\d)$
satisfies 
$$
\begin{cases}
-\D v(x)-\dfrac{h_i\big(\frac{x-a_i}{|x-a_i|}\big)}
{|x-a_i|^2}\,v(x)+\beta\,v(x)=0,&\text{in } B(a_i,\d),\\[10pt]
v=0\quad\text{and}\quad \dfrac{\partial v}{\partial\nu}=\a\,\d^{-\frac
  N2} \psi_1^{h_i}\big(\frac{x-a_i}{|x-a_i|}\big)
,&\text{on } \partial B(a_i,\d).
\end{cases}
$$
As a consequence the distribution $-\D v-\widetilde
V\,v+\beta\,v\in{\mathcal D}'(\R^N\setminus\{a_1,\dots,a_k\})$ acts as follows: 
$$ 
{}_{{\mathcal D}'(\R^N\setminus\{a_1,\dots,a_k\})}\big\langle
\!\!-\!\D v-\widetilde
V\,v+\beta\,v\,,\,\varphi\big\rangle_{C^{\infty}_{\rm
  c}(\R^N\setminus\{a_1,\dots,a_k\})}=\d^{-\frac
  N2}\a \int_{\partial B(a_i,\d)} \!
\!\psi_1^{h_i}\bigg(\frac{x-a_i}{|x-a_i|}\bigg)
\varphi(x)\,ds.
$$
Hence $h=-\D v-\widetilde V\,v+\beta\,v\in H^{-1}(\R^N)$ and  satisfies
(\ref{eq:47}) as $\a<0$. From Corollary \ref{l:nsac}, we finally deduce that 
the operator $-\D-V$ is not
essentially  
self-adjoint in $C^{\infty}_{\rm c}(\R^N\setminus\{a_1,\dots,a_k\})$.
\end{pfn}

The following theorem analyzes  essential self-adjointness of 
anisotropic Schr\"odinger operators with potentials carrying infinitely many
singularities  on reticular structures.  

\begin{Theorem}\label{t:self-adjo-ret}
Assume that
 $\{h_n\}_{n\in\N}\subset C^{\infty}\big({\mathbb
    S}^{N-1}\big)$ satisfy \eqref{eq:25}
and $\{a_n\}_n\subset\R^N$ satisfy \eqref{eq:55}
and $|a_n-a_m|\geq 1$
for all $n\neq m$.  Let $\e>0$, $\sigma>0$, and $\bar \delta>0$ 
as in Lemma \ref{l:ret},  $0<\delta<\min\{1/4,\bar \delta\}$,
 and  
$$
V(x)=\sum_{n=1}^{\infty}\frac{h_n\big(\frac{x-a_n}{|x-a_n|}\big)}{|x-a_n|^2}
\alchi_{{\mathcal E}^{\sigma,\delta}_{0,\tilde h_n}(a_n)},
\quad\text{where }\tilde h_n:=(1+\e)h_n.
$$
Then $-\D-V$  is essentially self-adjoint
in $C^{\infty}_{\rm c}\left(\R^N\setminus\{a_n\}_{n\in\N}\right)$ if
and only if 
$$
\mu_1(h_n)\geq-\bigg(\frac{N-2}2\bigg)^{\!\!2}+1\quad\text{for all }n\in\N.
$$
\end{Theorem}

\begin{pf}
For $n\in\N$, let $\phi_n\in C^\infty_{\rm c}(\R^N)$ such that $\phi_n\equiv 1$ in 
${\mathcal E}^{\sigma,\delta}_{0,\tilde h_n}(a_n)$, $\varphi\equiv 0$ in $\R^N\setminus 
{\mathcal E}^{\sigma,2\delta}_{0,\tilde h_n}(a_n)$, and $0\leq\phi_n\leq 1$. Then, setting
$$
\overline{V}(x)=
\sum_{n=1}^{\infty}\frac{h_n\big(\frac{x-a_n}{|x-a_n|}\big)}{|x-a_n|^2}
\phi_n(x),
$$
we have that $\overline{V}\in C^{\infty}(\R^N\setminus\{a_n\}_{n\in\N})$ and 
$V-\overline{V}\in L^\infty(\R^N)$. 

For $\alpha>\max\{0,-
\infess_{\R^N}(V-\overline{V})$ and $b(x):=V-\overline{V}+\alpha$,
there holds $b\in L^{\infty}(\R^N)$ and $\infess_{\R^N}b>0$.
From the Kato-Rellich Theorem 
the operator $-\D-V$ is essentially self-adjoint
in $C^{\infty}_{\rm c}\left(\R^N\setminus\{a_n\}_{n\in\N}\right)$ if
and only if 
$-\D- V+b=-\D-\overline{V}+\alpha$ is essentially self-adjoint.
In view of Lemma \ref{l:ret},  $-\D- V+b$ is
positive. Hence essential self-adjointness is equivalent to density of
$\mathop{\rm Range}(-\D- V+b)$ in $L^2(\R^N)$.

Let us first prove that,  if  $\inf_{n\in\N}\mu_1(h_n)\geq-\big(\frac{N-2}2\big)^2+1$, then
$-\D-V$  is essentially self-adjoint.  
Let $f\in C^{\infty}_{\rm
  c}\left(\R^N\setminus\{a_n\}_{n\in\N}\right)$ and $\e>0$. Then there
exists $0\leq\gamma<1$  satisfying  
$$
\inf_{n\in\N}\mu_1(h_n-\gamma)>-\bigg(\frac{N-2}2\bigg)^{\!\!2}+1
$$ 
and such that if  $u\in H^1(\R^N)$ solves
\begin{equation}\label{eq:57}
-\D u(x)- V_{\gamma}(x)u(x)+\alpha\,u(x)=f(x),\quad\text{where }
V_{\gamma}(x):=\sum_{n=1}^{\infty}\frac{h_n\big(\frac{x-a_n}{|x-a_n|}\big)-\gamma}{|x-a_n|^2}
\phi_n(x),
\end{equation} 
then 
\begin{equation}\label{eq:58}
\|(V_{\gamma}- \overline{V})u\|_{L^2(\R^N)}<\e.
\end{equation}
If $\inf_{n\in\N}\mu_1(h_n)>-\big(\frac{N-2}2\big)^2+1$ it is enough to 
choose $\gamma=0$.
If $\inf_{n\in\N}\mu_1(h_n)=-\big(\frac{N-2}2\big)^2+1$, 
from \cite[Theorem 1.2]{FMT2}, 
there exists a positive constant $\tilde C$
independent on $\gamma\in(0,1)$, such that all solutions of (\ref{eq:57})
can be estimated as 
$$
|u(x)|\leq
\tilde C\,|x-a_n|^{-\frac{N-2}{2}+\sqrt{1+\gamma}}\|u\|_{H^1({\mathcal
    E}^{\sigma,\delta'}_{0,\tilde h_n}(a_n))}\quad \text{in }{\mathcal
  E}^{\sigma,2\delta}_{0,\tilde h_n}(a_n),
$$
for some $2\delta<\delta'<1/2$ and for all $n\in\N$. 
We emphasize that, thanks to assumption (\ref{eq:25}), the constant $\tilde C$ in the above
estimate can be taken to be independent of $n$, as
one can easily obtain by scanning through the proof of 
\cite[Theorem 1.2]{FMT2} and checking the dependence of the estimate 
constant of the angular coefficient of the dipole.
Consequently
\begin{align*}
\left\| \frac{\gamma  \,  
\phi_nu}{|x-a_n|^2}\right\|_{L^2(\R^N)}^2
&\leq {\rm const\,}(\tilde C,N,\e,\sigma)\|u\|_{H^1({\mathcal
    E}^{\sigma,\delta'}_{0,\tilde h_n}(a_n))}^2\,\frac{\gamma^2}{\sqrt{1+\gamma}-1},
\end{align*}
and hence
\begin{align*}
\|(V_{\gamma}- \overline{V})u\|_{L^2(\R^N)}&\leq
{\rm const\,}(\tilde C,N,\e,\sigma)\frac{\gamma}{\sqrt{\sqrt{1+\gamma}-1}}\|u\|_{H^1(\R^N)}
\\
&\leq {\rm
  const\,}\|f\|_{L^2(\R^N)}\frac{\gamma}{\sqrt{\sqrt{1+\gamma}-1}}
\to 0\quad\text{as }\gamma\to0. 
\end{align*}
Therefore it is possible to choose $\gamma$ small enough in order to
ensure that all solutions of (\ref{eq:57}) satisfy~(\ref{eq:58}).
For such a $\gamma$,  the
Lax-Milgram Theorem provides a unique $w\in
H^1(\R^N)$  weakly solving
$$
-\D w(x)-V_\gamma(x)w(x)+\alpha \,w(x)=f(x) \quad\text{in }\R^N.
$$
Being $V_\gamma\in C^{\infty}(\R^N\setminus\{a_n\}_{n\in\N})$, from classical elliptic 
regularity theory, $w\in  C^{\infty}(\R^N\setminus\{a_n\}_{n\in\N})$.
From \cite{FMT2} and arguing as in the proof of Theorem
\ref{t:self-adjo}, we deduce that
\begin{align}\label{eq:100}
& w(x)\sim
|x-a_n|^{ \sigma_n},\quad\text{and}\\
&\label{eq:21}
\n w(x)
=
\begin{cases}
O\big(|x-a_n|^{\sigma_n-1}\big),
&\text{if }\mu_1(h_n)+\gamma<N-1,\\[5pt]
O\big(|x-a_n|^{-\tau}\big),&\text{if }\mu_1(h_n)+\gamma\geq N-1,
\end{cases}
\end{align}
as $x\to a_n$, where 
$\sigma_n=-\frac{N-2}{2}+\sqrt{(\frac{N-2}2)^2+\mu_1(h_n)+\gamma}$ and 
$0<\tau<\frac{N-2}{2}$. Since 
$$
\tilde\sigma:=\inf_{n\in\N}\sigma_n>2-\frac N2,
$$
 for all $j\in\N$, we can choose
$N_j\in\N$ such that $N_j\to+\infty$ as $j\to\infty$,
$N_jj^{4-N-2\tilde\sigma}\to 0$, and $N_j j^{2\tau-N+2}\to 0$, and let $R_j>0$
such that $R_j\to+\infty$ as $j\to\infty$ and $B(a_n,1/j)\subset
B(0,R_j)$ for all $n=1,\dots, N_j$.  Let $\eta_j$ be a cut-off
function such that $\eta_j\in C^{\infty}_{\rm
  c}\left(\R^N\setminus\{a_n\}_{n\in\N}\right)$, $0\leq\eta_j\leq 1$,
and
\begin{align*}
&\eta_j(x)\equiv 0\text{ in }
\bigcup_{n=1}^{N_j}B\Big(a_n, \frac1{2j}\Big)\cup (\R^N\setminus
B(0,2R_j)),\quad 
\eta_j(x)\equiv 1\text{ in }B(0,R_j)\setminus \bigcup_{n=1}^{N_j} B\Big(a_n,
\frac1{j}\Big),\\
&|\n \eta_j(x)|\leq C\,j\text{ in } \bigcup_{n=1}^{N_j} \Big(B\Big(a_n,
\frac1{j}\Big)\setminus B\Big(a_n, \frac1{2j}\Big)\Big),\quad
|\n \eta_j(x)|\leq \frac{C}{R_j}\text{ in }B(0,2 R_j) \setminus B(0,R_j),
\\
&|\D \eta_j(x)|\leq C\,j^2\text{ in } \bigcup_{n=1}^{N_j} \Big(B\Big(a_n,
\frac1{j}\Big)\setminus B\Big(a_n, \frac1{2j}\Big)\Big),\quad
|\D \eta_j(x)|\leq \frac{C}{R_j^2}\text{ in }B(0,2 R_j) \setminus B(0,R_j),
\end{align*}
for some positive constant $C$ independent of $j$ and $n$.
Setting  $f_j:=\eta_j f-2\n \eta_j\cdot\n w-w\,\D\eta_j$, we have that 
  $\eta_j f\to f$ in $L^2(\R^N)$, and,
from  (\ref{eq:100}), 
\begin{align*}
&\int_{\R^N}|\n \eta_j(x)|^2|\n w(x)|^2\,dx
\\[5pt]
&\leq {\rm const\,} j^2\sum_{\substack{n=1\\\mu_1(h_n)+\gamma<N-1}}^{N_j}\int_{B(a_n,
\frac1{j})\setminus B(a_n, \frac1{2j})}|x-a_n|^{2(\sigma_n-1)}\,dx\\
&\quad+
{\rm const\,} j^2\sum_{\substack{n=1\\\mu_1(h_n)+\gamma\geq N-1}}^{N_j}\int_{B(a_n,
\frac1{j})\setminus B(a_n, \frac1{2j})}|x-a_n|^{-2\tau}\,dx+
\frac{\rm
const\,}{R_j^2} 
\int_{B(0,2R_j)\setminus B(0,R_j)}|\n w(x)|^2\,dx
\\[5pt]
&\leq 
  {\rm const\,}\left[
  N_j\, j^{4-N-2\tilde\sigma}+N_j\,
  j^{2\tau-N+2}+R_j^{-2}\|w\|_{H^1(\R^N)}\right]
\end{align*}
and, in a similar way,
\begin{align*}
  \int_{\R^N}|\D \eta_j(x)|^2|w(x)|^2\,dx\leq {\rm const\,}\big[N_j\,
  j^{4-N-2\tilde\sigma}+R_j^{-4}\big].
\end{align*}
By the choice of $N_j$, we deduce that $f_j\to f$ in $L^2(\R^N)$. 
 Hence, for $j$ large,
$\|f_j-f\|_{L^2(\R^N)}<\e$. 
The functions $w_j:=\eta_j w\in C^{\infty}_{\rm c}(\R^N\setminus\{a_n\}_{n\in\N})$
solve  
$$
-\D w_j(x)- V(x) w_j(x)+b\,w_j(x)=g_j(x),
$$
where $g_j(x):=f_j(x)+(
V_{\gamma}(x)- \overline{V}(x))w_j(x)$, i.e. $g_j\in  \mathop{\rm
  Range}(-\D-V+b)$.
Moreover, from (\ref{eq:58}) and  $|\eta_j|\leq 1$, we deduce that
$$
\|g_j-
f_j\|_{L^2(\R^N)}=\|(V_{\gamma}-\overline{V})w_j\|_{L^2(\R^N)}\leq
\|(V_{\gamma}- \overline{V})w\|_{L^2(\R^N)}<\e,
$$ 
hence $\|g_j-
f\|_{L^2(\R^N)}<2\e$ for 
large $n$. The density of 
$\mathop{\rm Range}(-\D-\
V+b)$  in $C^{\infty}_{\rm
  c}\left(\R^N\setminus\{a_n\}_{n\in\N}\right)$ and consequently  in $L^2(\R^N)$ 
is proved.

\medskip\noindent The proof of non essential self-adjointness in the
case $\mu_1(h_n)<-\big(\frac{N-2}2\big)^2+1$ for some $n\in\N$
 can be obtained just by mimicking the arguments
of the proof of Theorem \ref{t:self-adjo} and Corollary \ref{l:nsac}.\end{pf}

\begin{remark}\label{rem:retessselfadj}
If $\{h_n: n\in\N\}$ is finite, i.e. if only a finite number of possible
angular coefficients is repeated in the reticle, then the assumption 
of $C^\infty$-smoothness of $h_n$ required in Theorem \ref{t:self-adjo-ret}
can be removed, as one can easily check arguing 
by approximation
as in the proof of Theorem \ref{t:self-adjo}.
\end{remark}

\appendix
\section*{Appendix}
\setcounter{section}{1}
\setcounter{equation}{0}
\setcounter{Theorem}{0}

\begin{Lemma}\label{l:app}
  Let $\{h_n\}_{n\in\N}\subset L^{\infty}({\mathbb S}^{N-1})$ and
  $h\in L^{\infty}({\mathbb S}^{N-1})$ such that
$$
h_n\to h\text{ a.e. in }{\mathbb S}^{N-1}\quad\text{and}\quad
\sup_n\|h_n\|_{L^{\infty}({\mathbb S}^{N-1})}<+\infty.
$$
Then 
\begin{align*}
  {\rm i)}\quad &\lim_{n\to+\infty}\mu_1(h_n)=\mu_1(h);\\[5pt]
{\rm ii)} \quad &\psi_1^{h_n} \text{ converge to }\psi_1^h
\text{ uniformly in ${\mathbb S}^{N-1}$ and in $H^1({\mathbb S}^{N-1})$};\\[5pt]
{\rm iii)} \quad & \frac{h_n\big(\frac{x}{|x|}\big)}{|x|^2}\text{ converge to }
\frac{h\big(\frac{x}{|x|}\big)}{|x|^2}\text{ in the Marcinkiewicz space }L^{N/2,\infty}.
\end{align*}
\end{Lemma}
\begin{pf}
For all $n\in\N$, let  $\psi_n=\psi_1^{h_n}$ be the positive $L^2$-normalized
eigenfunction associated to the first eigenvalue
$\mu_1(h_n)$, i.e.
\begin{equation}\label{eq:29}
-\Delta \psi_n-h_n\psi_n=\mu_1(h_n)\psi_n, \quad\text{in }{\mathbb S}^{N-1},\quad
\text{and}\quad
\int_{{\mathbb S}^{N-1}}\psi_n^2=1.
\end{equation}
Since $\{h_n\}_n$ is bounded in $L^{\infty}({\mathbb S}^{N-1})$, it is easy to verify that 
$\{\mu_1(h_n)\}_n$ is bounded in $\R$, hence it admits a subsequence, still denoted as 
$\{\mu_1(h_n)\}_n$, such that $\mu_1(h_n)\to \bar \mu$ as $n\to+\infty$ for some 
$\bar\mu\in\R$. By a standard bootstrap argument, it  follows 
that $\{\psi_n\}_n$ is relatively compact in $H^1({\mathbb S}^{N-1})$ and 
bounded in $C^{0,\alpha}({\mathbb S}^{N-1})$ for some positive $\alpha$.
In particular the sequence $\{\psi_n\}_n$ is equicontinuous and hence, by the Ascoli-Arzel\`a 
Theorem, there exists $\bar\psi\in H^1({\mathbb S}^{N-1})\cap C^0({\mathbb S}^{N-1})$
such that $\psi_n\to\bar\psi$ in $H^1({\mathbb S}^{N-1})$ and uniformly in
${\mathbb S}^{N-1}$. 

Passing to the limit in \eqref{eq:29}, strong $H^1$-convergence of $\psi_n$ to
$\psi$ and Dominated Convergence's Theorem
 yield  that $\bar\psi$ satisfies
$$
-\Delta \bar\psi-h\bar\psi=\bar \mu\bar\psi, \quad\text{in }{\mathbb S}^{N-1},\quad
\text{and}\quad
\int_{{\mathbb S}^{N-1}}\psi^2=1,
$$
and therefore
$$
\bar \mu\geq \mu_1(h).
$$
On the other hand, for all $\varphi\in H^1({\mathbb S}^{N-1})\setminus\{0\}$
$$
\mu_1(h_n)\leq
\frac{\int_{\mathbb S^{N-1}}|\n_{\mathbb
S^{N-1}}\varphi(\theta)|^2\,dV(\theta)- \int_{\mathbb S^{N-1}}h_n(\theta)
\varphi^2(\theta)\,dV(\theta)}{\int_{\mathbb S^{N-1}}\varphi^2(\theta)\,dV(\theta)},
$$
hence, letting $n\to+\infty$, 
$$
\bar\mu\leq
\frac{\int_{\mathbb S^{N-1}}|\n_{\mathbb
S^{N-1}}\varphi(\theta)|^2\,dV(\theta)- \int_{\mathbb S^{N-1}}h(\theta)
\varphi^2(\theta)\,dV(\theta)}{\int_{\mathbb S^{N-1}}\varphi^2(\theta)\,dV(\theta)},
$$
which implies that $\bar\mu\leq\mu_1(h)$. Then $\mu_1(h)=\bar\mu$ and $\bar\psi=
\psi_1^h$. Statements i--ii) follow from above and  the Uryson
property.

A direct calculation shows that
$$
\left\|\frac{(h_n-h)\big(\frac{x}{|x|}\big)}{|x|^2}\right\|_{L^{N/2,\infty}(\R^N)}
=N^{-2/N}\|h_n-h\|_{L^{N/2}({\mathbb S}^{N-1})}
$$
and hence statement iii) follows from the assumption on $\{h_n\}_n$ and the 
Dominated Convergence Theorem.
\end{pf}

\end{document}